\documentclass[final,narrowdisplay]{elsart1p} 
\usepackage{graphicx}
\usepackage{amssymb}
\usepackage{amsmath} 
\include{TeXFigs} 
\usepackage{amsfonts}
\usepackage{amsmath} 
\usepackage{amssymb}
\usepackage{psfig}
\usepackage{verbatim}
\usepackage{amsmath}
\usepackage[mathcal]{euscript}
\numberwithin{equation}{section}
\newcommand \RR    {{I\!\! R}} 
\newcommand \eps   {\epsilon} 
\newcommand \Sbold    {{\mathbf S}} 
\newcommand \del   {\partial}
\newcommand \lam   {\lambda}
\newcommand \vfa   {\varphi^\flat_\alpha}
\newcommand \psifa   {\psi^\flat_\alpha}
\newcommand \psifaq   {\psi^{\flat,q}_\alpha}
\newcommand \talpha {{\tilde \alpha}}
\newcommand \be   {\begin{equation}}
\newcommand \ee   {\end{equation}}
\newtheorem{definition}{Definition}[section]

\newtheorem{conjecture}[definition]{Conjecture}

\begin{document}

\begin{frontmatter}

\title{Why many theories of shock waves are necessary.
Kinetic functions, equivalent equations, and fourth-order models}  

\author{Philippe G. LeFloch}

\address{Laboratoire J.-L. Lions \& Centre National de la Recherche Scientifique, 
\\
Ê      Universit\'e de Paris VI, 4 Place Jussieu, 
 Ê     75252 Paris,
 Ê     France.
\\
E-mail: \texttt{LeFloch@ann.jussieu.fr}}

\author{Majid Mohammadian}

\address{Courant Institute of Mathematical Sciences, New York University, 
\\
251 Mercer Street, New York, NY 10012, USA. 
\\
E-mail: \texttt{Majid@cims.nyu.edu}}

\begin{abstract}
We consider several systems of nonlinear hyperbolic conservation laws describing the dynamics of nonlinear waves
in presence of phase transition phenomena.  
These models admit under-compressive shock waves which are not uniquely determined by a standard 
entropy criterion but must be characterized by a kinetic relation.   
Building on earlier work by LeFloch and collaborators, 
we investigate the numerical approximation of these models by {\sl high-order} finite difference schemes, 
and uncover several new features of the kinetic function associated with 
with physically motivated second and third-order regularization terms, especially
viscosity and capillarity terms. 
 
On one hand, the role of the equivalent equation associated with a finite difference scheme is 
discussed. We conjecture here and demonstrate numerically that the (numerical) kinetic function
associated with a scheme approaches the (analytic) kinetic function associated with the given model  
---~especially since its equivalent equation approaches the regularized model at a higher order.  
On the other hand, we demonstrate numerically that a kinetic function can be associated
with the thin liquid film model and the generalized Camassa-Holm model.
Finally, we investigate to what extent a kinetic function can be associated with the equations of 
van der Waals fluids, whose flux-function admits two inflection points. 
\end{abstract}

\begin{keyword}
hyperbolic equation \sep  conservation law \sep  shock wave \sep
kinetic relation \sep  viscosity \sep  capillarity \sep  equivalent equation \sep thin liquid film 
\sep Camassa-Holm. 

\PACS 35L65, 76L05

{\bf To appear in the Journal of Computational Physics.}  

\end{keyword}
\end{frontmatter}


\section{Introduction}
\label{IN-0}

\subsection{Background} 

In this paper we study the numerical approximation of several first-order 
nonlinear hyperbolic systems of conservation laws, and
we consider discontinuous solutions generated by supplementing the hyperbolic equations with 
higher-order, physically motivated, vanishing regularization terms. Specifically, we consider complex
fluid flows when physical features such as viscosity and capillarity effects 
can not be neglected even at the hyperbolic level of modeling, and 
need to be taken into account. This gives rise to {\sl many shock wave theories} 
associated with any given nonlinear hyperbolic system: depending on the underlying small-scale physics
of the problem under consideration, one need a different selection of ``entropy solutions''. 

It is well known that solutions of nonlinear hyperbolic equations become discontinuous in finite time,
and may, therefore, exhibit shock waves. {\sl Classical} shock waves satisfy standard entropy criteria
(due to Lax, Oleinik, Wendroff, Liu, Dafermos, etc); they are compressible and stable under
perturbation and approximation.

On the other hand, discontinuous solutions of hyperbolic problems may also exhibit
{\sl non-classical, under-compressive} shock waves -- also referred to as 
{\sl subsonic phase boundaries} in the context of phase transition theory.
To uniquely characterize under-compressive shocks one need to impose a jump condition that is
not implied by the given set of conservation laws and is called a {\sl kinetic relation.} 
The selection of physically meaningful shock waves of a first-order hyperbolic system is
determined by traveling waves associated with an augmented system
that includes viscosity and capillarity effects. That is, one searches for scale-invariant 
solutions depending only on the variable $y :=x- \lam \, t$ for some speed $\lam$ and connecting 
two constant states $(\tau_-, u_-)$ and $(\tau_+, u_+)$ at infinity. The characterization of these 
``admissible'' discontinuities is based on kinetic relations. (For background, 
see \cite{LeFloch2} and the references cited therein.)  

For instance, 
one important model of interest in fluid dynamics describes liquid-vapor flows governed by van der Waals's equation 
of state. For pioneering mathematical works on van der Waals fluids we refer to Slemrod et al.
 \cite{Slemrod1,Slemrod2,FanSlemrod}, who investigated self-similar approximations to the Riemann problem.  
The concept of a kinetic relation associated with (undercompressive) nonclassical shocks or phase boundaries 
was introduced by Abeyaratne and Knowles \cite{AbeyaratneKnowles1,AbeyaratneKnowles2}, 
Truskinovsky \cite{Truskinovsky1,Truskinovsky2}, and first analyzed mathematically  
by LeFloch \cite{LeFloch1}. 
Kinetic relations and nonclassical shocks were later 
studied extensively by Shearer et al. \cite{JMS,SS} (phase transition), 
LeFloch et al. \cite{HayesLeFloch1,HayesLeFloch2,HayesLeFloch3,LeFloch3,BedjaouiLeFloch} 
(traveling waves, Riemann problem, general hyperbolic systems),  
and by Bertozzi, Shearer, et al. \cite{BertozziShearer,BertozziMunchShearer} (thin film model), 
as well as Colombo, Corli, Fan, and others (see  \cite{ColomboCorli,CorliFan,FanLiu,FanSlemrod} 
and the references cited therein).

\subsection{Approximation of undercompressive waves} 

The present paper built on existing numerical work done by the first author and his collaborators
\cite{HayesLeFloch1,HayesLeFloch2,LeFlochRohde,ChalonsLeFloch1,ChalonsLeFloch2}
and devoted to the numerical investigation of non-classical shocks and phase boundaries
generated by diffusive and dispersive terms kept in balance. These papers cover 
scalar conservation laws, and systems of two or three conservation laws arising in fluid dynamics
(Euler equations) and material science. Entropy stable schemes were constructed, 
and kinetic relations were computed numerically. 
Numerical kinetic functions turn out to be very useful to 
to evaluate the accuracy and efficiency of the schemes under consideration.

Numerical experiments were also performed on undercompressive shocks for the
models of thin films \cite{BertozziMunchShearer,LS,Munch}
and phase dynamics \cite{NganTruskinovsky}. These authors did not investigate the role of 
the kinetic relation for these models, and one of our aims in the present paper is to tackle this issue. 

Recall that one effective numerical strategy to compute under-compressive waves is provided by 
the Glimm and front tracking schemes, 
for which both theoretical and numerical results are now available 
\cite{LeFloch2,ChalonsLeFloch2}. 
Both schemes converge to the correct non-classical solutions to any given Cauchy problem. 
The main feature of these schemes is to completely avoid spurious numerical dissipation or dispersion, 
which  
this guarantees their convergence to the physically meaningful solution selected by the kinetic relation. 
However, this class of schemes has some drawbacks: they are limited to first-order accuracy 
and require the precise knowledge of the Riemann solver, while 
numerical solutions may exhibit a noisy behavior. 

In contrast, in the last twenty years, modern, high-order accurate, finite difference techniques of 
(classical) shock capturing have been developed to deal with discontinuous solutions of hyperbolic problems. 
Extending these techniques to compute nonclassical shocks turned out to be quite difficult, however. 
Our aim in the present paper is to pursue this investigation of 
the interplay between dissipative and dispersive mechanisms in hyperbolic models, 
and to better understand how they compare with similar mechanisms taking place in finite difference schemes. 

This problem was first tackled by Hayes and LeFloch \cite{HayesLeFloch1,HayesLeFloch2} where the
importance of the equivalent equation associated with a scheme was pointed out and kinetic functions for schemes
were computed. Next, the role of high-order, entropy conservative schemes was emphasized and 
further kinetic functions were determined numerically \cite{LeFlochRohde}.
In the situation just described, the numerical solutions usually contain mild oscillations,
which vanish out when the mesh is refined. This feature is entirely consistent with the behavior 
of traveling wave solutions to the underlying dispersive model. In consequence, 
total variation diminishing (TVD) techniques should not be used to compute nonclassical 
shocks. Observe also that the diffusive-dispersive model serves only to provide a mechanism to select 
``admissible'' solutions of associated hyperbolic equations. It has been established, for many hyperbolic models, 
that physically meaningful solutions can be characterized uniquely by pointwise conditions on shocks 
(Rankine-Hugoniot jump conditions, an entropy inequality, a kinetic relation, and, possibly, 
a nucleation criterion). 

Finite difference schemes can also be studied for their own sake, and one important
issue is whether criteria can be found for a given scheme to generate nonclassical shocks or not.
Some heuristics were put forward in \cite{HayesLeFloch2} and allow one to distinguish between 
the following cases:  
\begin{itemize}
\item[(1) ] Schemes satisfying (a discrete version of) all of the entropy inequalities. 
In the context of scalar conservation laws, this is the case of 
the schemes with monotone numerical flux functions. For systems of two conservation laws
(which admit a large family of (convex) mathematical entropies) 
this class includes the Godunov and the Lax-Friedrichs schemes (for small data, at least).
These schemes, if convergent, must converge to a weak solution satisfying all entropy inequalities which, 
consequently, coincides with the classical solution selected by 
the Oleinik or Kruzkov entropy conditions (scalar equations) 
and the Wendroff or Liu entropy conditions (systems of two or more conservation laws). 

\item[(2) ] Schemes satisfying a single (discrete) entropy inequality but applied to a system
 with genuinely nonlinear 
characteristic fields. The classification of Riemann solutions given in \cite{HayesLeFloch3} 
shows that, again, only the scheme can converge to the classical entropy solution, only. 

\item[(3) ] Schemes satisfying a single (discrete) entropy inequality but applied to a system
with non-genuinely nonlinear characteristic fields. To decide whether such schemes are expected 
to generate nonclassical shocks, one should determine the equivalent equation. 
One may truncate the equivalent equation and keep only the first two terms. 
An analysis of the properties of traveling waves for this continuous model provides an indication 
of the expected behavior of the scheme. Interestingly, the behavior depends 
on the sign of the dispersion coefficient and the sign of the third-order derivative of the flux.  
For instance, for first-order schemes applied to scalar 
equations one typically obtains, after further linearization in the neighborhood of 
the origin $0$ in the phase space,   
$$
v_t \pm (v^3)_x = h \, v_{xx} + \alpha h^2 \, v_{xxx}. 
$$
Nonclassical shocks have been observed when the flux is concave-convex (that is, $v^3$) and $\alpha$ is positive, 
or else when the flux is convex-concave (that is, $-v^3$) and the coefficient $\alpha$ is negative.  
\end{itemize}

\subsection{Purpose of this paper} 
 
We will demonstrate here that the existence of under-compressive shocks and several 
typical behaviors of these nonlinear waves, especially the existence of a kinetic function, 
are properties shared by {\sl many examples} arising in continuum physics.  
To make our point, we present a number of physical models describing nonlinear wave dynamics dynamics: 
cubic flux, thin liquid films, generalized Camassa-Holm, van der Waals fluids. 
Specifically, we prove that kinetic functions can be associated 
to each of these models and we study their monotonicity, dependence upon (viscosity, capillarity, mesh) parameters,
and behavior in the large. We uncover several new features of the kinetic function that have not been
observed theoretically via analytical methods yet. It is our hope that the conclusions reached here numerically
will motivate further theoretical developments in the mathematical theory of non-classical shocks.
The work also provides further ground that not a single theory of entropy solutions but, rather, 
{\sl many theories of shock waves}
are required to accurately describe singular limits of hyperbolic equations, as supported by the
framework developed in \cite{LeFloch2,LeFloch3,LeFlochShearer}.

The outline of the paper is as follows. In Section~\ref{MO-0}, we present the physical models of interest
and discuss briefly their analytical properties. In Section~\ref{EQ-0}, after introducing some background on
non-classical shocks and kinetic relations, we investigate the role of the equivalent equation.
In Section~\ref{KI-0}, we establish the existence of kinetic functions associated with each of the models
and investigate their properties. In Section~\ref{MU-0}, we investigate van der waals fluids.
Finally, Section~\ref{CO-0} contains concluding remarks.


\section{Models of interest}
\label{MO-0}

We begin with a brief presentation of a few nonlinear hyperbolic models arising in continuum physics.
\subsection{Cubic conservation law} 
It will be convenient to start with an academic example consisting 
of a conservation law whose flux-function admits a non-degenerate inflection point. 
For simplicity and with little loss of generality as far as the local behavior near the inflection point is concerned, 
we can assume that the flux is a cubic function. After normalization, 
we arrive at the {\sl cubic conservation law} 
\be
u_t + (u^3)_x = 0, \quad u=u(t,x) \in \RR,  \, t \geq 0.
\label{EQ.1}
\ee
We are interested in (discontinuous) solutions that can be realized as limits of diffusive-dispersive solutions of
\be
u_t + (u^3)_x =  \eps \, u_{xx} + \alpha \, \eps^2 \, u_{xxx}, \quad u=u_\alpha^\eps(t,x),
\label{EQ.2}
\ee
where $\alpha$ is a fixed parameter and $\eps \to 0$. Recall that shock waves of (\ref{EQ.1})
are solutions containing a single propagating discontinuity connecting two states $u_-, u_+$
at the speed $\lambda$. These constants must satisfy the Rankine-Hugoniot relation
$$
\lambda = {u_+^3 - u_-^3 \over u_+ - u_-} = u_-^2 + u_- \, u_+ + u_+^2,
$$
as well as the entropy inequality associated with the quadratic entropy ${U(u) := u^2}$
$$
-\lam \, (u_+^2 - u_-^2 ) + {2 \over 3} \, (u_+^4 - u_-^4) \leq 0.
$$
As is now well-known, the solutions
$u_\alpha := \lim_{\eps \to 0} u^\eps_\alpha$
may contain both {\sl classical shock waves} satisfying the standard compressibility condition
(equivalent here to the entropy criteria introduced by Lax)
\be
3 \, u_-^2 = f'(u_-) \geq \lambda \geq f'(u_+) = 3 \, u_+^2,
\label{EQ.3}
\ee
as well as {\sl non-classical shock waves} which turn out to be {\sl under-compressive}
\be
\lambda < f'(u_\pm).
\label{EQ.34}
\ee

The characterization of the solutions generated by (\ref{EQ.2}) as $\eps \to 0$ is provided 
as follows.
Based on an analysis of all possible traveling wave solutions of  (\ref{EQ.2}), one can see that,
for a given viscosity/capillarity ratio $\alpha$ and for every left-hand state $u_-$,
there exists a single right-hand state
\be
u_+ =\vfa(u_-),
\label{EQ.5}
\ee
that can be attained by a non-classical shock. The function $\vfa$
is called the {\sl kinetic function} associated with the model (\ref{EQ.2}).
The existence of the kinetic function has been established theoretically for a large class of flux-functions
and nonlinear diffusion-dispersion operators, including (\ref{EQ.2}). The results are often stated in terms
of the shock set, $\Sbold_\alpha(u_-)$, 
consisting of all right-hand states $u_+$ that can be attained from a given left-hand state $u_-$
by a classical or by a non-classical shock. Note also that instead of the relation (\ref{EQ.5}), one
can equivalently prescribe the entropy dissipation of a non-classical shock, that is
the kinetic relation can be expressed in the following form observed in \cite{LeFloch1}: 
\be
\aligned
&  - \lambda \, (U(u_+) - U(u_-)) + F(u_+) - F(u_-) 
\\
& = - \int_\RR U''(v(y)) \, v_y(y)^2 \, dy \leq 0, 
\endaligned 
\label{EQ.6}
\ee
where $y \mapsto v(y)$ denotes the traveling wave trajectory
connecting $u_-$ to $u_+$. In orther words, the entropy dissipation must be prescribed on
a non-classical shock. 
As far as the specific model with cubic flux and linear diffusion
and dispersion is concerned, the kinetic function and the shock set
can be expressed explicitly by analytical formulas as follows. The
kinetic function associated with (\ref{EQ.2}) reads
 \be
 \vfa(u_-)=
 \begin{cases}
  - u_- - \talpha/2,
     &         u_- \leq - \talpha,
    \\
    - u_- /2,
    &         |u_-| \leq \talpha,
    \\
  - u_- + \talpha/2,
    &         u_- \geq \talpha,
 \end{cases}
 \label{EQ.7} \ee
with $\talpha := \sqrt{(8 /3\alpha)}$, while the corresponding shock
set is
 $$
 \Sbold_\alpha(u_-) \, = \, \begin{cases} (u_-, \talpha/2] \cup
 \big\{-u_- - \talpha/2 \big\},
               &   u_-  \leq -\talpha,
 \\
 [- u_- /2, u_-),     &  - \talpha \leq u_- \leq \talpha,
 \\
 \big\{-u_- + \talpha/2 \big\} \cup [- \talpha/2, u_-),
               &   u_-  \geq \talpha.
 \end{cases}
 $$
Observe that $\vfa$ converges to $-u/2$ when $\alpha \to 0$, and
that the shock set converges to the standard interval $[-u/2, u]$
determined by the Oleinik entropy inequalities. See
\cite{LeFloch2} for details.
%
%
\subsection{Thin liquid films}
The {\sl thin liquid film model} ($\eps, \alpha$ being positive, scaling parameters) 
\be
u_t+(u^2-u^3)_x = \eps \, (u^3 \, u_x)_x  - \alpha \eps^2 \, (u^3 \, u_{xxx})_x
\label{AP.1}
\ee
describes the dynamics of a thin film with height $u=u(t,x)$ moving on an inclined flat solid surface. 
The (non-convex) flux-function
$$
f(u) = u^2 - u^3, \quad 0 \leq u \leq 1,
$$
represents the competing effects of the gravity and a surface stress
known as the Marangoni stress. The latter arises in experiments due
to an imposed thermal gradient along the solid surface. The
fourth-order diffusion is due to surface tension, and the
second-order diffusion represents a contribution of the gravity to
the pressure. 
This model was introduced and extensively studied by Bertozzi, M\"unch, and Shearer \cite{BertozziMunchShearer}.
The existence of non-classical traveling waves was established analytically in \cite{BertozziShearer}, 
and various numerical studies were performed \cite{Munch,LS} which exhibited subtil
properties of stability and instability 
of these waves. 
For more analytical background, see also the convergence theory developed by Otto and Westdickenberg \cite{OW}. 

From the standpoint of the general well-posedness theory the kinetic relation 
is an important object which is necessary to uniquely characterized the physically meaningful solutions. 
A kinetic function was not exhibited in the work \cite{BertozziShearer} which, instead, relied on non-constructive arguments.
From the existing literature, it is not clear whether
a concept of a kinetic relation could be associated with the equation \eqref{AP.1} and, if so, 
and whether such a kinetic function would enjoy the same monotonicity properties, 
as the ones established earlier  
for conservation laws regularized by viscosity and capillarity. This issue will be addressed in Section~\ref{KI-0}.

Observe that the flux $f(u)$ here is positive on the interval $(0,1)$, increasing on $(0,2/3)$ and decreasing on $(2/3,1)$. 
It admits a single inflection point at $u=1/3$.  To every point $u \in (0,1/3)$ we can associate the ``tangent point'' 
$\varphi(u) \in (1/3, 1/2)$ characterized by 
$$
f'(\varphi^\natural(u)) = {f(u) - f(\varphi^\natural(u)) \over u - \varphi^\natural(u)}, 
$$
or 
$$
\varphi^\natural(u) := {1 - u \over 2}. 
$$
The same formula maps also the interval $(1/3,1)$ onto $(0,1/3)$. This function $\varphi^\natural$ 
allows us to define 
{\sl classical} shock waves associated with the equation. When $u \in (0,1/3)$ the left-hand state 
$u$ can be connected to the right-hand state
$\varphi^\natural(u)$ by a contact discontinuity. When $u \in (1/3,1)$, $\varphi^\natural(u)$ is the left-hand state
and $u$ is the right-hand state.  

Another important function is provided by considering the {\sl entropy dissipation} 
$$
D(u_-,u_+) := - {\lambda \over 2} (u_+^2 - u_-^2) + {2 \over 3} (u_+^3 - u_-^3) + {3 \over 4} (u_+^4 - u_-^4). 
$$
where the shock speed $\lambda$ is given by 
$$
\lambda := {f(u_+) - f(u_-) \over u_+ - u_-} = u_+ + u_- - (u_+^2 + u_+ u_- + u_-^2).  
$$
The {\sl zero dissipation function}  $\varphi^\sharp$ is by definition the {\sl non-trivial} root of $D$, i.e.
$$
D(u, \varphi^\sharp(u)) = 0, 
$$
or
$$
\varphi^\sharp(u) := {2 \over 3} - u. 
$$
The function $\varphi^\sharp$ maps the interval $(0,2/3)$ onto itself.

According to the theory in \cite{LeFloch2} the range of the nonclassical shocks
 is limited by the functions $\varphi^\natural$ and $\varphi^\sharp$. 
Precisely, for a nonclassical shock connecting $u_- < 1/3$ to $u_+ >1/3$, the right-hand state must have
$$
\varphi^\natural (u_-)= {1 - u_- \over 2} \leq u_+ < \varphi^\sharp(u_-)= {2 \over 3} - u_-.  
$$
The sign are reversed for a decreasing nonclassical shock. 


\subsection{Generalized Camassa-Holm model}

We consider a generalized version of the {\sl Camassa-Holm equation} 
\be
\label{AP.camassa} 
\aligned 
u_t + f(u)_x 
& = \eps \, u_{xx} + \alpha \eps^2 \big( u_{xxt} + 2 \, u_x u_{xx}+ u u_{xxx}\big),  
\\
& = \eps \, u_{xx} + {\alpha \over 2} \eps^2 \big( 2 u_{xt} + (u^2)_{xx} - (u_x)^2\big)_x,  
\endaligned
\ee
which arises as an asymptotic higher-order model of wave dynamics in shallow water.  
The second-order and third-order terms are related to the viscosity and the capillarity 
of the fluid. When the flux $f$ is nonconvex, for instance
$$
f(u) = u^3
$$
nonclassical shocks may in principle arise. We will demonstrate in this paper that indeed 
solutions may exhibit nonclassical shocks and establish the existence of an associated kinetic function.


\subsection{Van der Waals fluids}

Compressible fluids are governed by the following two conservation laws: 
\be
\begin{split}
& \del_t \tau - \del_x u = 0, 
\\
& \del_t u + \del_x p(\tau) = \eps \del_{xx} u  
- \alpha \, \eps^2 \del_{xxx} \tau.
\end{split}
\label{VD.1} 
\ee 
Here, $u$ and $\tau$  represent the velocity and the specific volume of the fluid, respectively, while
$\alpha$ is a non-negative parameter
representing the strength of the viscosity. 
The pressure law $p=p(\tau)$ is a positive
function defined for all $\tau \in ( 0, +\infty)$ and of the following {\sl van der Waals} type: 
there exist $0 < a < c$ such that 
\be
\begin{split}
& p''(\tau) > 0,  \quad \tau \in (0,a) \cup (c,+\infty),
\\
& p''(\tau) < 0,  \quad \tau \in ( a, c),
\\
& p'(a) >0, 
\end{split}
\label{VD.3} 
\ee 
and 
 \be 
 \lim_{\tau \to 0} p(\tau) = +
\infty, \qquad \lim_{\tau \to +\infty} p(\tau) = 0. 
\label{VD.2}
\ee 
The left-hand side of \eqref{VD.1} forms a first-order system of partial differential equations, which is of
elliptic type when $\tau$ belongs in the interval $(d, e)$ characterized by the conditions $0 < d < a < e < c$ and $p'(d) =
p'(e) = 0$. It is of hyperbolic type when $\tau \in (0, d) \cup (e, + \infty)$ and admits the two (distinct, real) wave speeds $
\pm \sqrt{ - p'(\tau)}$.

\

\section{A conjecture on the equivalent equation}
\label{EQ-0}

\subsection{Kinetic functions associated with difference schemes}

We now begin the discussion of the numerical approximation of the solution $u_\alpha$ generated by (\ref{EQ.2}).
The discussion applies to general conservation laws of the form
$$
\del_t u + \del_x f(u) = 0,
$$
where the flux-function $f$ admits a single inflection point.
To any finite difference scheme associated with (\ref{EQ.1}) one can in principle associate a kinetic function,
which we denote by $\psifa$. It may seem natural to request that $\psifa$ coincides with the kinetic function
$\vfa$ associated with the given model. However, it has been observed by Hayes and LeFloch \cite{HayesLeFloch2}
that, at least for all finite difference schemes that have been considered so far,
\be
\psifa \neq \vfa.
\label{EQ.8}
\ee
This discrepancy is due to the fact that the dynamics of non-classical shocks is determined by small-scale features of the
continuous model (\ref{EQ.2}) which can never be fully mimicked by a discrete model.

Given this perspective, the next natural question is to compare the kinetic functions $\psifa$ and $\vfa$.
We require that a scheme be a good approximation of the given, continuous model in the
following sense.  Its equivalent equation obtained by formal Taylor expansion should have the ``correct'' form
\be
v_t + (v^3)_x =  \eps \, v_{xx} + \alpha \, \eps^2 \, u_{xxx} + O(h^q), \quad h = c \, \eps,
\label{EQ.9}
\ee
where $c$ is a constant and $q$ represent the order of accuracy of the scheme.
Here, $v$ denotes the numerical solution and $h$ denotes the discretization parameter.
It is important to observe that all of the schemes considered in the present paper are first-order accurate, only,
as far as the hyperbolic equation (\ref{EQ.1}) is concerned.
They are, however, high-order approximations of the augmented model (\ref{EQ.2}).
For clarity, we emphasize the dependence in $q$ and denote by $\psifaq$ the kinetic function
associated with a scheme whose equivalent equation is (\ref{EQ.9}).

Let us introduce in this paper the following:

 \begin{conjecture}
 \label{conjecture}
As $q \to \infty$ the kinetic function $\psifaq$ associated with a
scheme with equivalent equation (\ref{EQ.7}) converges to the exact
kinetic function $\vfa$, \be \lim_{q \to \infty} \psifaq = \vfa.
\label{EQ.10} \ee
 \end{conjecture}

\

Rigorous results pointing toward the validity of this conjecture can be found in \cite{BKL}
which studies the role of relaxation terms in traveling wave solutions of conservation laws. A 
closely related problem
was tackled by Hou and LeFloch \cite{HouLeFloch} who considered nonconservative schemes for the
computation of nonlinear hyperbolic problems.
In these problems, small-scale features are critical in selecting shock waves, and the equivalent equation
have been found to provide a guide to designing difference schemes. We will not try here
to support the above conjecture on theoretical grounds, but we propose to investigate it numerically.
As we will see, very careful experiments are necessary.
We consider a large class of schemes based on standard differences, and obtained by approximating
the spatial derivatives arising in (\ref{EQ.2}) --with $\eps$ replaced by the mesh size $h$
(up to a constant multiplicative factor)--
$$
f(u)_x,  \quad   h \, u_{xx}, \quad    h^2 \, u_{xxx},
$$
by high-order finite differences,
so that the overall scheme is {\sl of order $q$ at least.} For completeness we list below
the corresponding expressions. We denote by $x_i$ the points of spacial mesh and by $u_i$ the approximation of the solution
at the point $x_i$. We also use the notation $f_i := f(u_i)$.
\
Based on the above we arrive at semi-discrete schemes for the functions $u_i=u_i(t)$.
For instance, using fourth order discretizations above we obtain
 \be
 \begin{split}
 {d u_i \over dt} = & - {1 \over h} \, \Big( {1\over
 12}f_{i-2}-{2\over 3}f_{i-1}+{2\over 3}f_{i+1}-{1\over 12}f_{i+2}
 \Big)
 \\
 & + {\eps \over h} \, \Big( -{1\over 12}u_{i-2}+{4\over
 3}u_{i-1}-{5\over 2}u_{i}+{4\over 3}u_{i+1}-{1\over 12}u_{i+2} \Big)
 \\
 & + {\alpha \, \eps^2 \over h} \, \Big( -{1\over 2}u_{i-2}+ u_{i-1}-
 u_{i+1}+{1\over 2}u_{i+2} \Big).
 \end{split}
 \label{EQ.11} \ee
Observe that this is a fully conservative scheme, in the sense that
$\sum_i u_i(t)$ is independent of $t$ (assuming, for instance,
periodic boundary conditions).
\
To actually implement the above algorithm, we use a Runge-Kutta scheme. Defining $U(t)= (u_i(t))$
(with ${i=\ldots, -1,0,1,\ldots}$),  
the semi discrete scheme takes the form
\be
{dU \over dt}(t) = R[U(t)],
\label{EQ.12}
\ee
where $R[U(t)]$ represents the spatial discretization. This system of ordinary differential equations
is solved numerically by employing an $s$-stage Runge-Kutta scheme defined
as follows
 \be
 \begin{split}
 & g^k = R\Big(U^n + \Delta t \, \sum_{j=1}^{k-1} a_{k,j} \,
 g^j\Big),
 \\
 & U^{n+1} : = U^n + \Delta t \, \sum_{k=1}^{s} b_{k} \, g^k.
 \end{split}
 \label{EQ.13} 
 \ee
For example, the non-zero coefficients of the fourth-order
Runge-Kutta scheme are given by 
\be
\aligned 
& a_{2,1}=1/2,a_{3,2}=1/2,a_{4,3}=1,
\\
& b_1=1/6,b_2=1/3,b_3=1/3,b_4=1/6.
\endaligned 
\label{EQ.14} 
\ee 
Coefficients of a sixth and an eighth-order
Runge-Kutta scheme are found in  \cite{TsitourasPapakostas} and \cite{Tsitouras}, respectively.
\ 


\subsection{Numerical experiments} 

\begin{figure} 
\begin{center}
\begin{tabular}{cc}
\begin{minipage}{6.5cm}
\begin{center}
\includegraphics*[height=6.5cm, width=6.5cm]{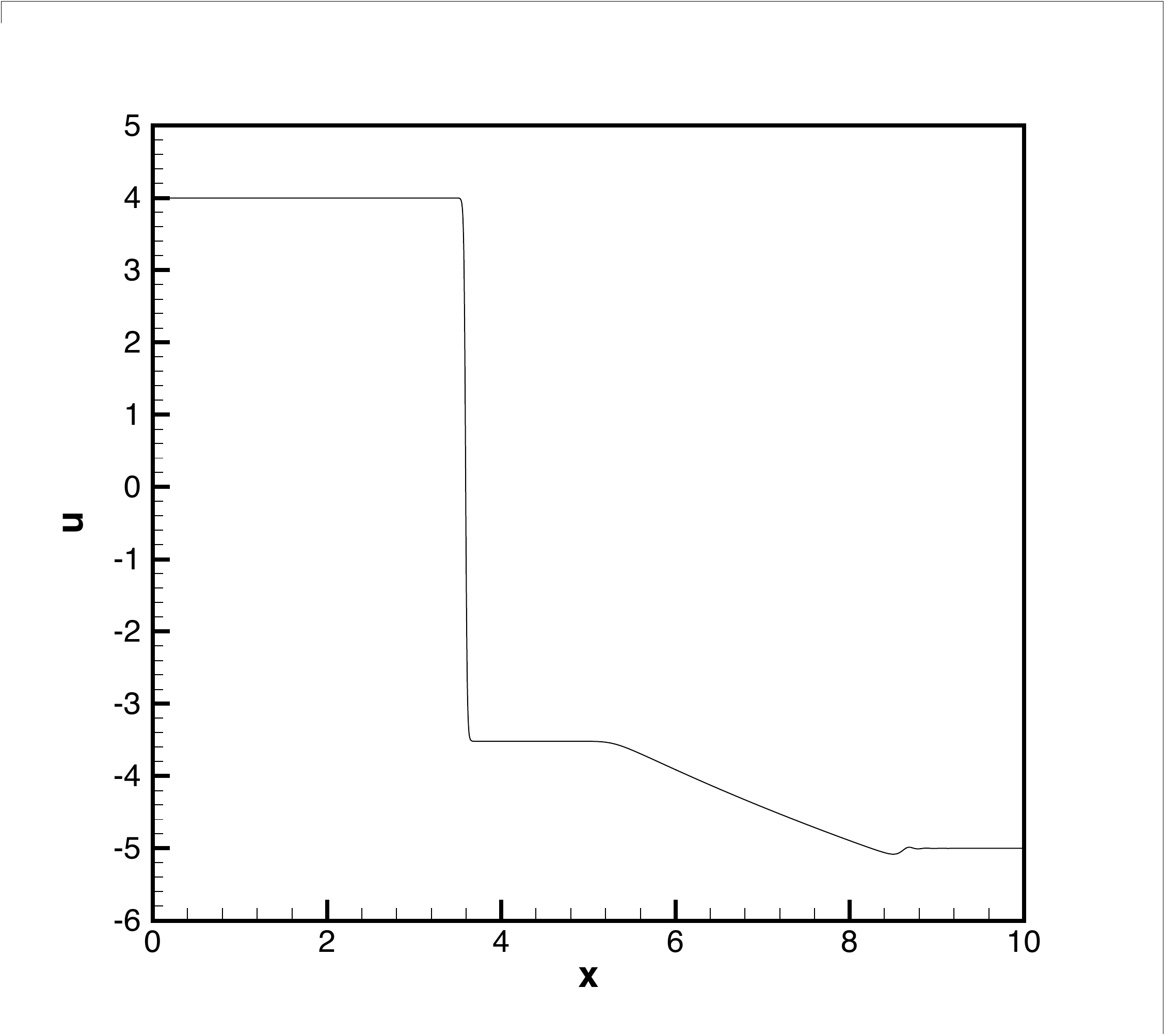}
\end{center}
\end{minipage}
&
\begin{minipage}{6.5cm}
\begin{center}
\includegraphics*[height=6.5cm, width=6.5cm]{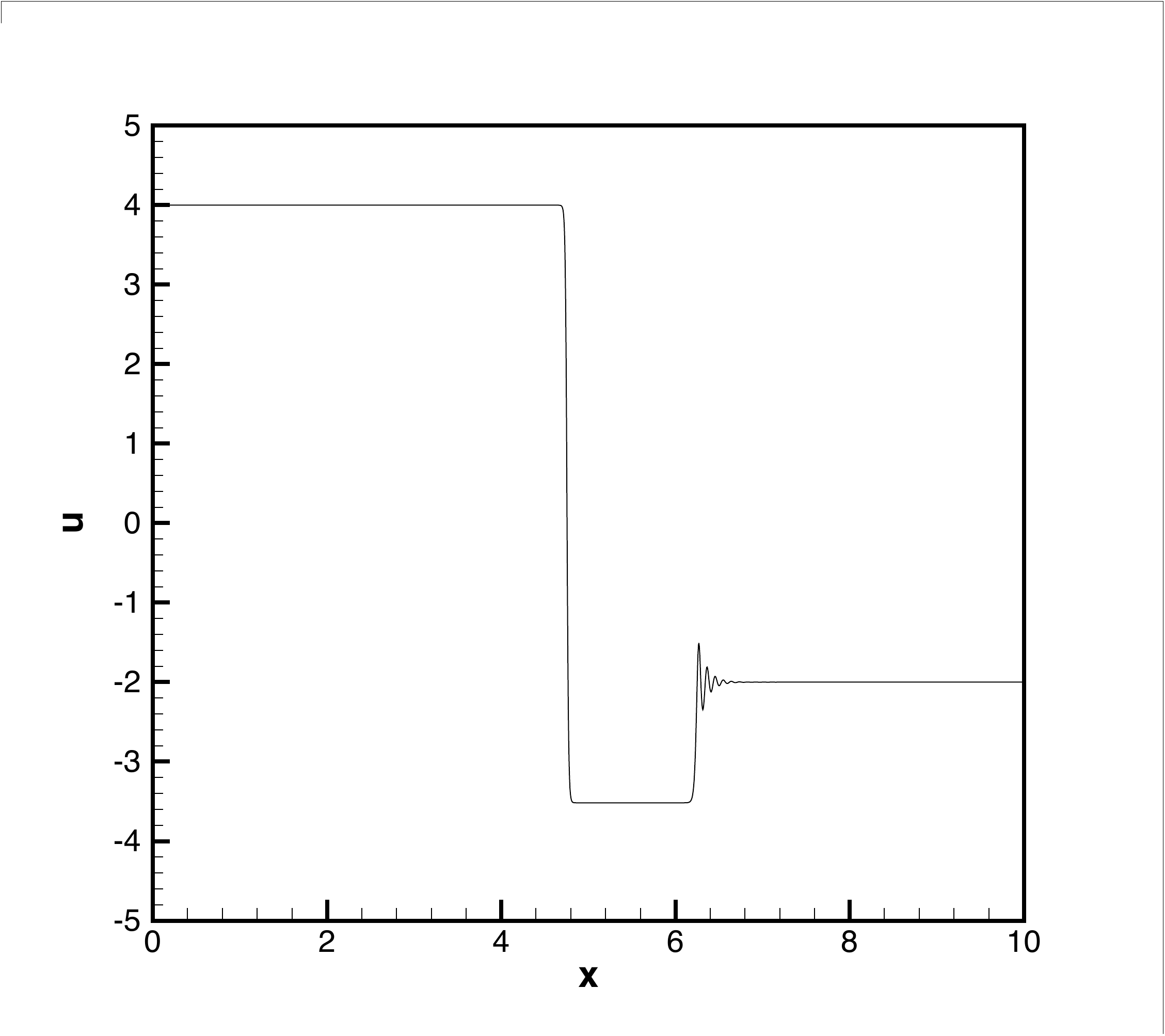}
\end{center}
\end{minipage}
\end{tabular}
\caption{Typical cases of non-classical waves for the cubic flux function.}
\label{typicalwaves}
\end{center}
\end{figure}

We now determine numerically the kinetic function associated
with the schemes described in the previous section. For each $q=4,6,8,10$ and for selected
values of the parameter $\alpha$ we compute the function $\psifaq$.
For the cubic flux function, two typical types of non-classical waves arise, which are 
under-compressive shock followed by a rarefaction wave (Figure \ref{typicalwaves}, left),
 and a double shock structure (Figure~\ref{typicalwaves}, right).
In order to plot a single kinetic function it is necessary to solve a large number of Riemann problems for various values
of left-hand state $u_-$ and to ensure that the right-hand state is picked up in an interval
where a non-classical shock does exist. In the numerical experiments it is often convenient to
select the right-hand state so that the Riemann solution contain two shocks. According to the
construction algorithm in \cite{LeFloch3} the middle state between the two shocks is therefore
the kinetic value $u_+ =\psifaq(u_-)$.

Several difficulties arise. First, when $u_-$ is close to the origin all the waves become very weak and
it is numerically difficult to identify the middle state. Second, in a range of parameter values, the two
waves may propagate with speeds that are very close and this again makes difficult the computation of the
middle state. Finally, solutions do contain mild oscillations, especially for small values of
$\alpha$ and this again introduces some numerical error. Due to these constraints we need to use
a rather fine mesh, with $h$ of the order of $1/1000$.

The following plots allow us to investigate numerically the convergence of the kinetic function $\psifaq$
toward the exact kinetic function, (\ref{EQ.7}), which is a piecewise affine function of the variable $u_-$.

In the following numerical tests, the CFL number was taken to be as large as possible in each run and
was identical for all schemes (4th order to 10th order).
It was observed that increasing the order of Runge-Kutta scheme (in the temporal integration) from four to six and eighth practically
does not change the results,
therefore, a tenth order Runge-Kutta scheme was not examined. On the other hand,
the order of spatial discretization was found to be very important and effectively made an important difference.

\begin{figure}
\begin{center}
  \includegraphics[width=8cm,height=8cm]{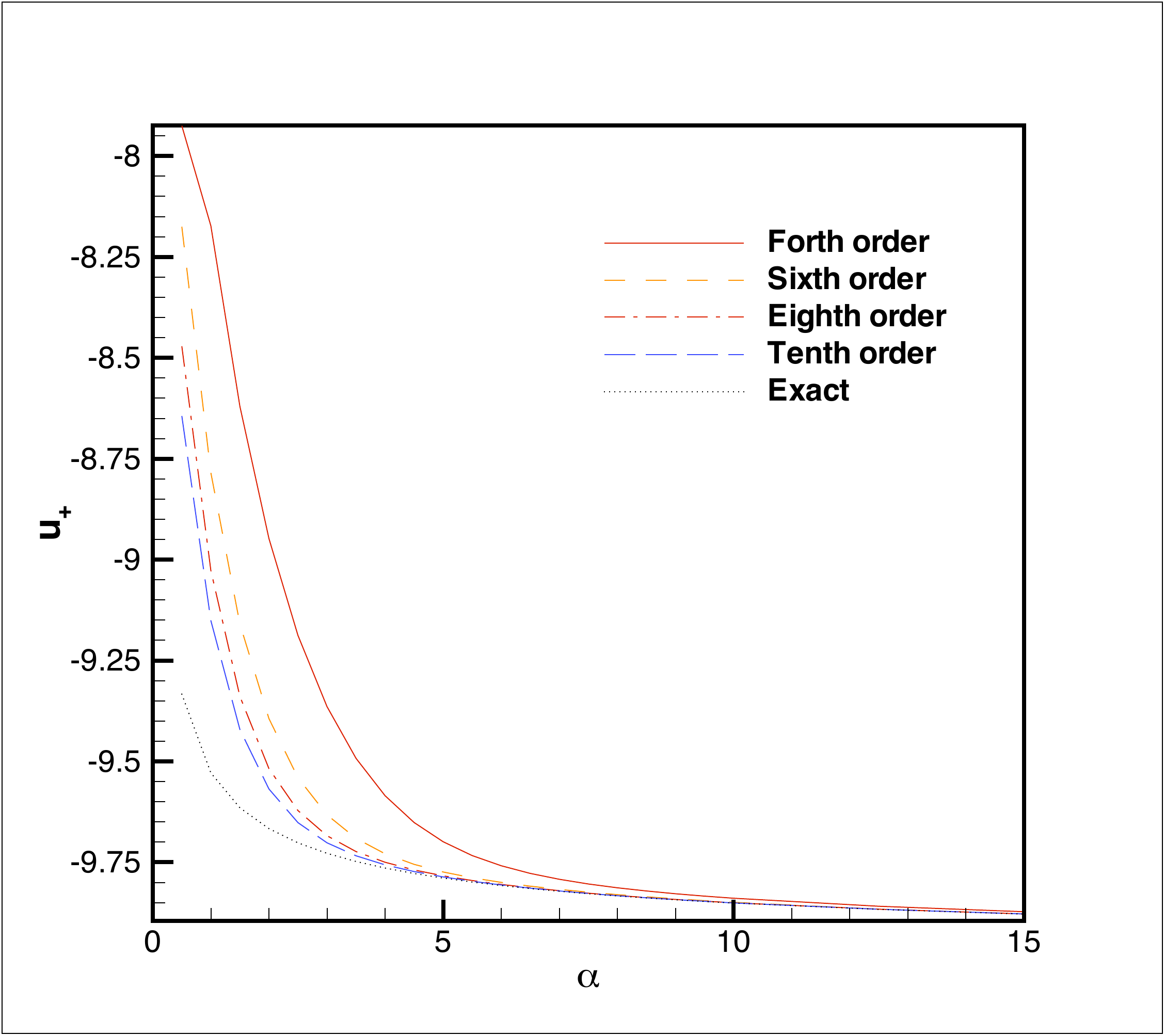}\\
  \caption{$u_+$ versus $\alpha$  for different schemes.}
  \label{alpha}
\end{center}
\end{figure}

First of all, Figure \ref{alpha} shows the right-hand state $u_+$ versus the parameter $\alpha$  for different schemes.
Here, we have used $h=0.001$, $u_-=10$, and $\eps=5h$. 
As $\alpha$ increases, the results of different order numerical schemes become closer and closer to the exact solution. 
For large values of $\alpha$, e.g. $\alpha>10$, even the fourth-order scheme gives satisfactory results.  
Since the solution of (\ref{EQ.1}) is, in fact, the limit of diffusive-dispersive solutions of
 (\ref{EQ.2}), with $\alpha$ fixed and $\eps \to 0$, we conclude that (provided $\alpha$ is sufficiently large)
 even the fourth order scheme converges to the exact solution. Figure~\ref{alpha} shows also that, for a fixed value of 
 $\alpha$ (sufficiently large), as the order of accuracy increases, the numerical solution converges to the exact one.
 These results strongly support Conjecture~\ref{conjecture}.

\begin{figure} 
\begin{center}
\begin{tabular}{ccc}

  & $u_+$ versus $u_-$ &  $\phi(s)/s^2$ versus $s$%
\\
$\alpha$=1 &
\begin{minipage}{6.0cm}
\begin{center}
\includegraphics*[height=5.5cm, width=5.5cm]{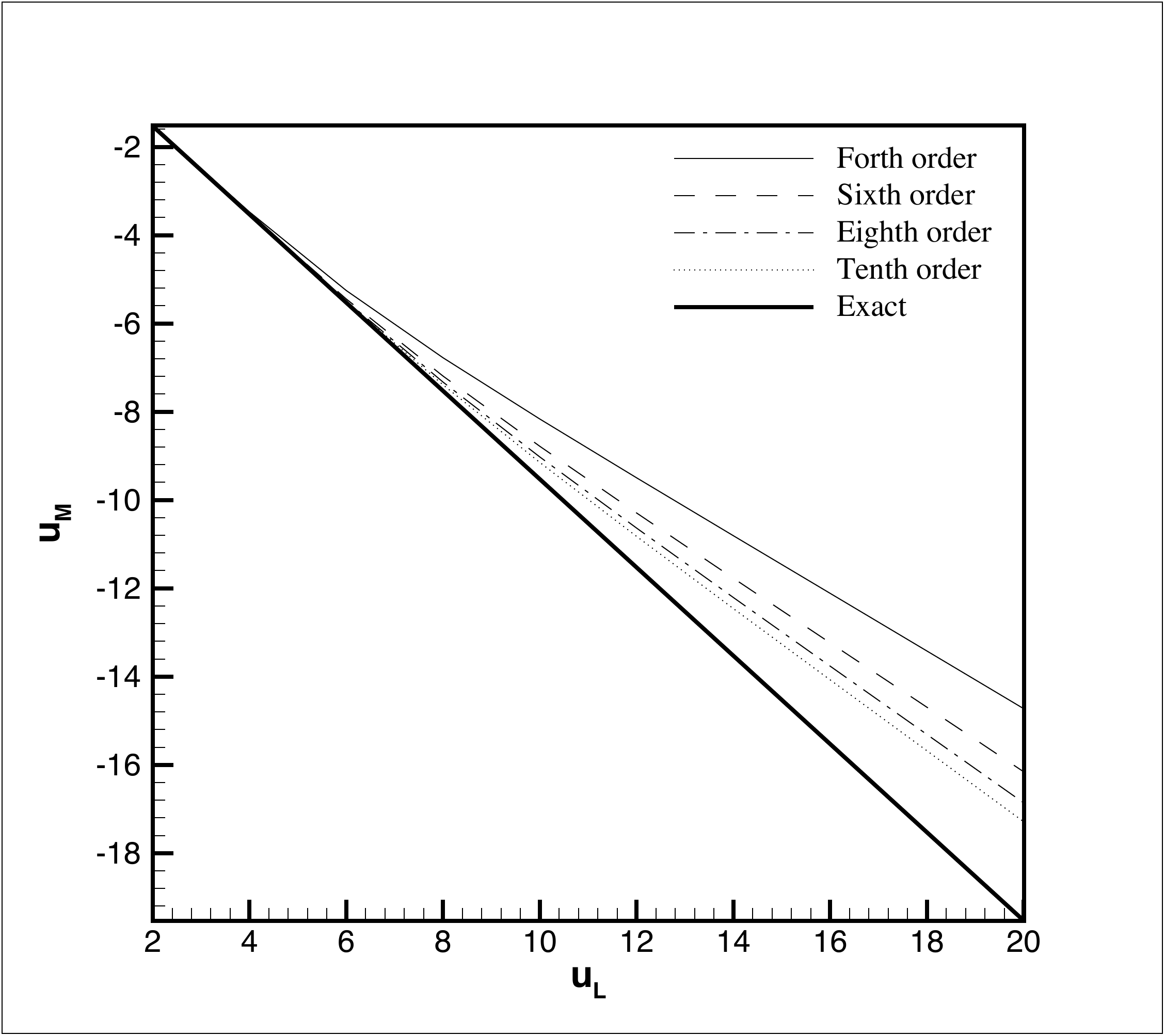}
\end{center}
\end{minipage}
&
\begin{minipage}{6.0cm}
\begin{center}
\includegraphics*[height=5.5cm, width=5.5cm]{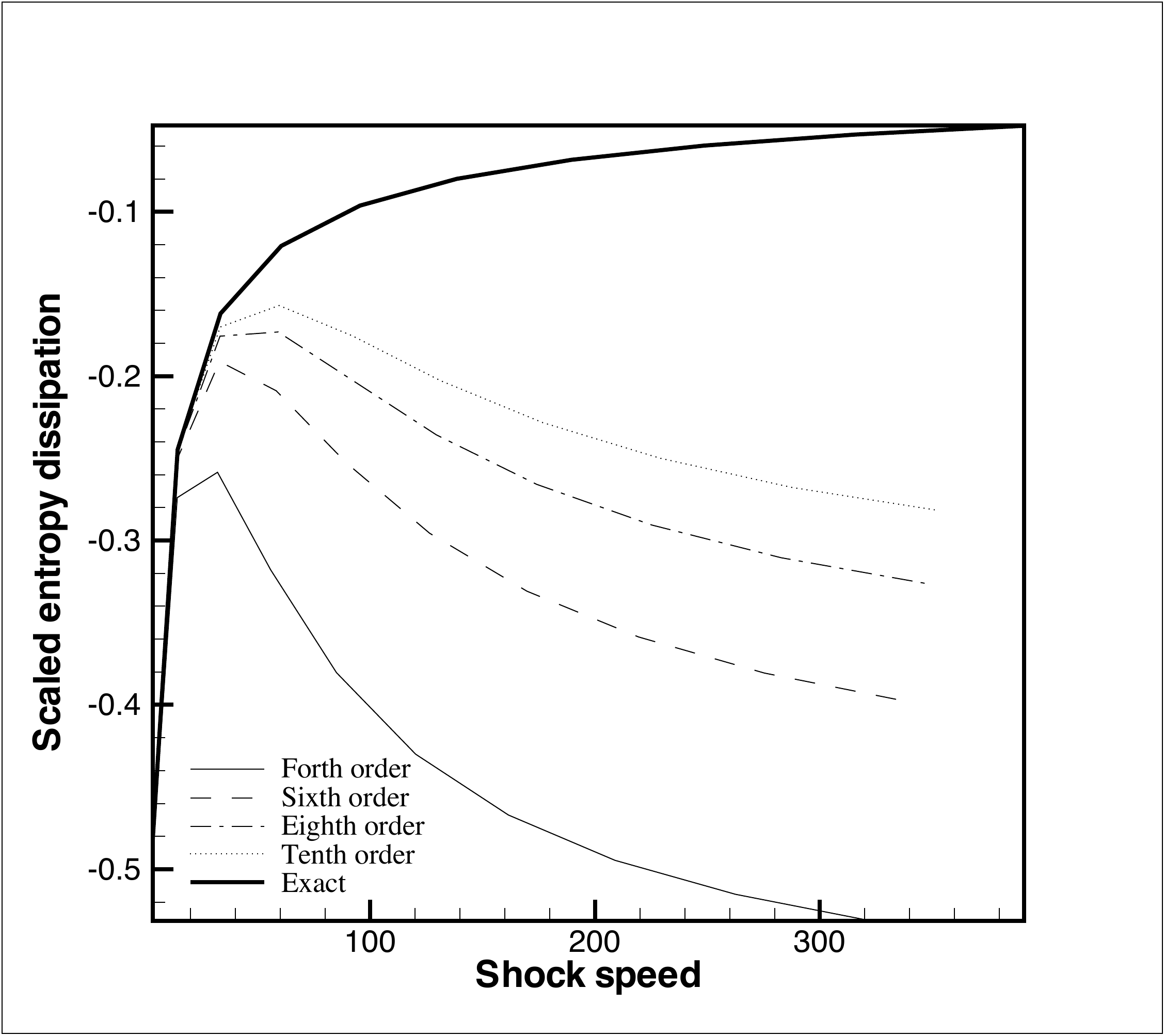}
\end{center}
\end{minipage}

\\
\\

$\alpha$=4 &
\begin{minipage}{6.0cm}
\begin{center}
\includegraphics*[height=5.5cm, width=5.5cm]{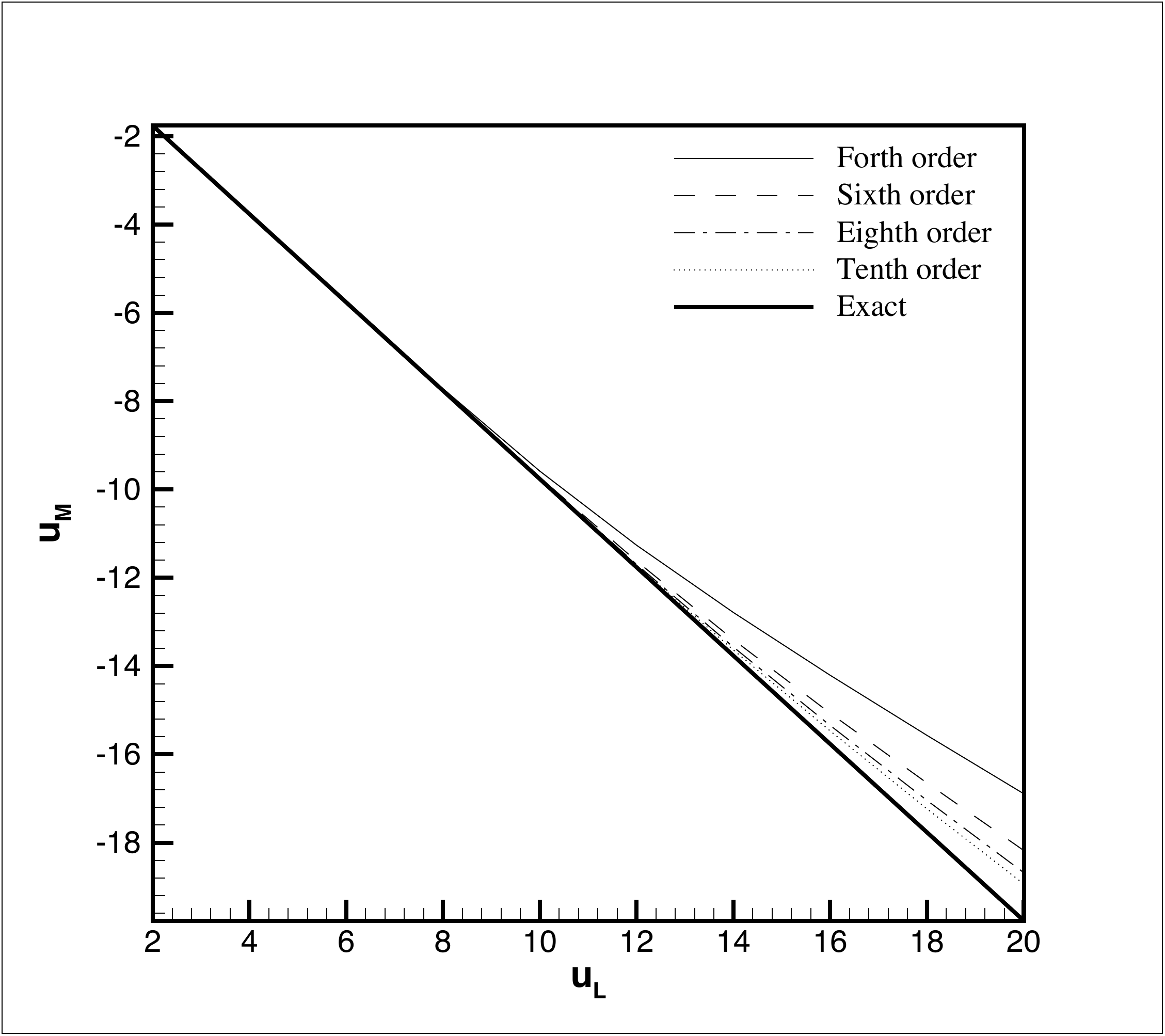}
\end{center}
\end{minipage}
&
\begin{minipage}{6.0cm}
\begin{center}
\includegraphics*[height=5.5cm, width=5.5cm]{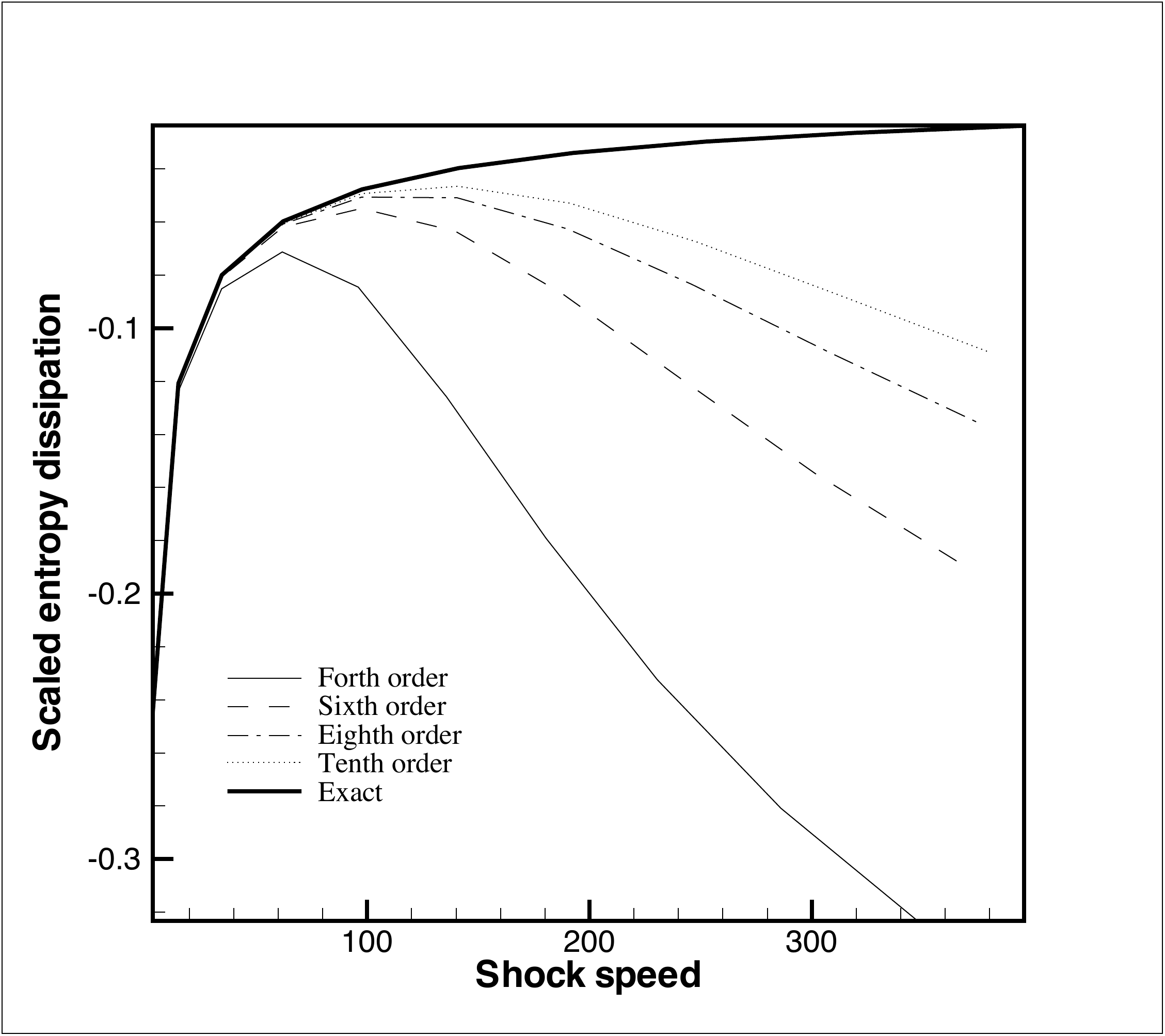}
\end{center}
\end{minipage}
\\
\\

$\alpha$=6 &

\begin{minipage}{6.0cm}
\begin{center}
\includegraphics*[height=5.5cm, width=5.5cm]{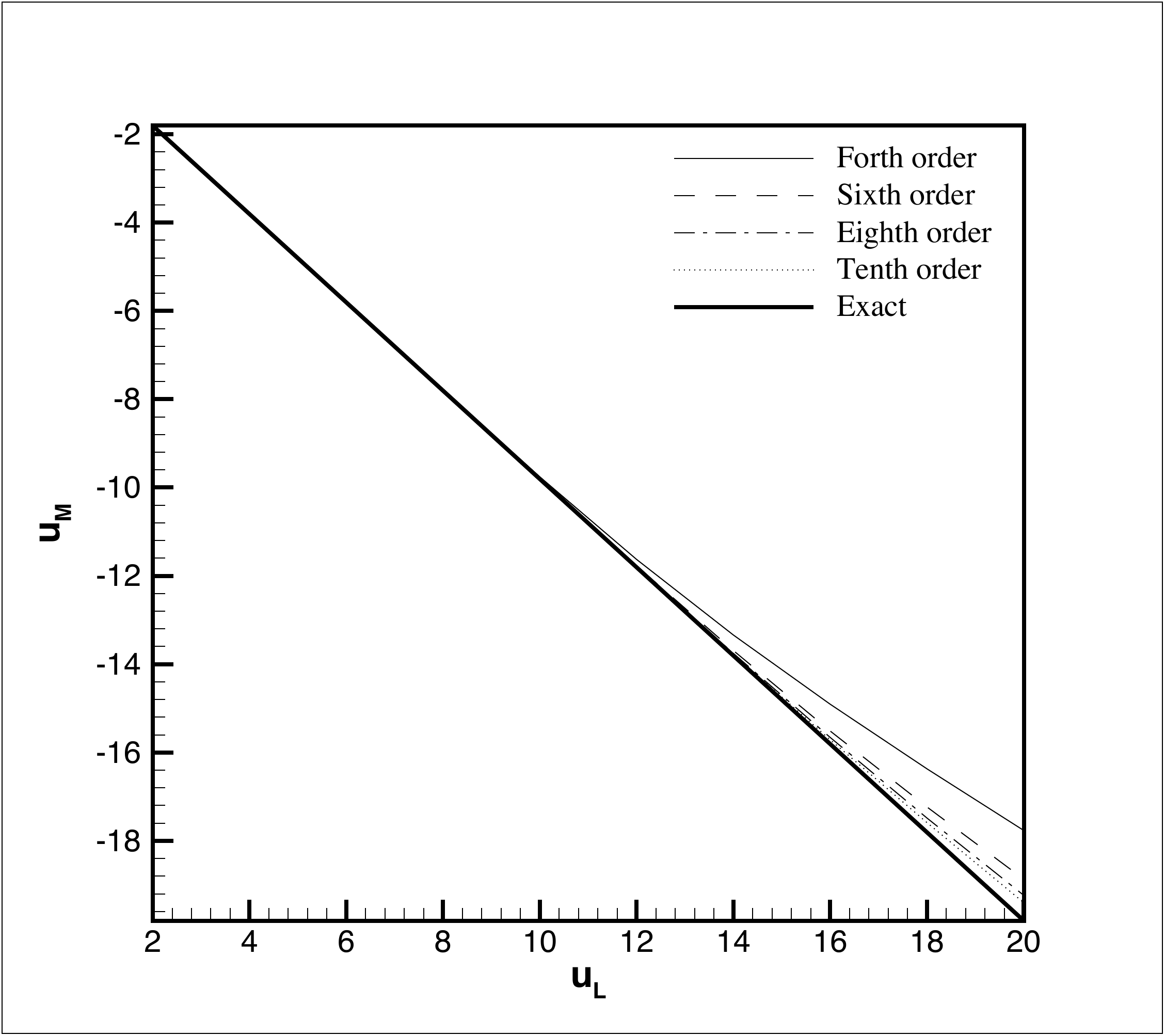}
\end{center}
\end{minipage}
&
\begin{minipage}{6.0cm}
\begin{center}
\includegraphics*[height=5.5cm, width=5.5cm]{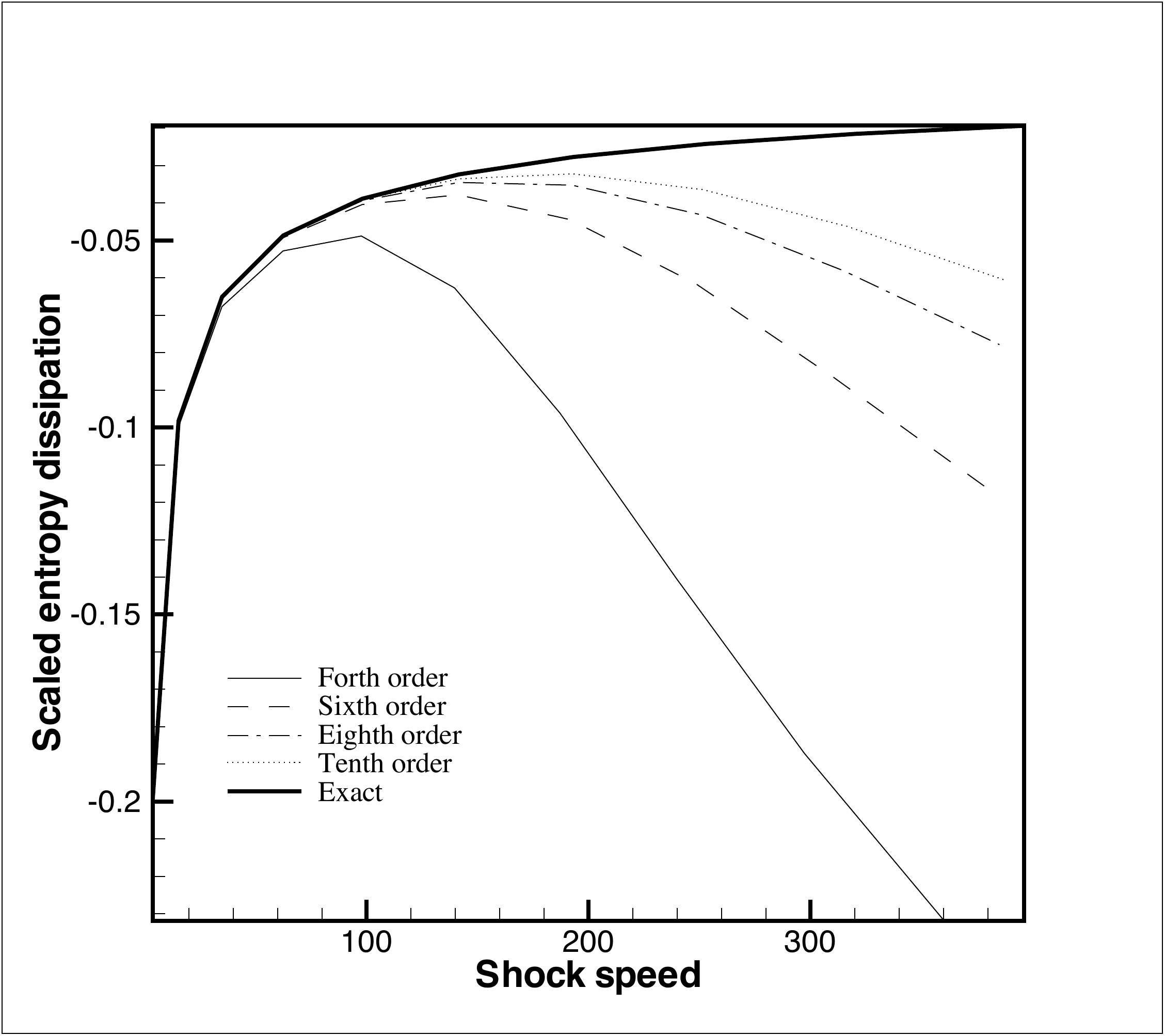}
\end{center}
\end{minipage}
%


\end{tabular}

\caption{$u_+$ versus $u_-$  (left) and scaled entropy dissipation
$\phi(s)/s^2$ versus shock speed $s$ (right) for different
values of $\alpha$. }

\label{cubicfig}

\end{center}
\end{figure}

Next, Figure~\ref{cubicfig} (left column) shows $u_+$ versus $u_-$ for $\alpha$=1, 4 and 6 respectively, with $h=0.005$, and $\epsilon=5h$. 
A grid of $2000$ nodes was used in all runs. Again, by increasing the order of accuracy,
the numerical solution converges to the exact one, for sufficiently small left-hand state $u_-$.
However, the suitable value of $\alpha$ depends on  $u_-$ and increases as $\alpha$ does. It should be mentioned that
very high values of $\alpha$ (which are in need for large $u_-$), are not satisfactory since  they lead to high oscillatory
results (recall that $\alpha$ is dispersion-diffusion ratio). It was also observed that for small values of $\alpha$, high frequency oscillations occur
before the non-classical shock, while for large values of $\alpha$, low frequency oscillations take place after the rarefaction (not shown).

Finally, Figure \ref{cubicfig} (right column) shows the scaled entropy dissipation $\phi(s)/s^2$ versus the
shock speed $s$ for $\alpha=1,4,6$, respectively. We use the quadratic entropy $U(u)=u^2/2$ in computing the entropy dissipation. In terms of $u_-$ and $u_+$, this is
\be
\phi(u_-,u_+)=(u_+-u_-)^2(u_+^2-u_-^2),
\label{EQ.15}
\ee
and from the Runkine-Hugoniot relation, the shock speed is
\be
 s=u_-^2+u_-u_++u_+^2.
\label{EQ.16}
\ee
The same feature is observed as before and by increasing the order of accuracy; the numerical entropy dissipation converges to the exact one (provided $\alpha$ is sufficiently large), which is in accordance with 
Conjecture~\ref{conjecture}.

%

\

\subsection{Spatial discretization parameter versus diffusive-dispersive parameters}

To conclude this section let us  
set $h = c \, \eps$ and plot the corresponding numerical results
with the sixth-order scheme and some fixed values of $h$ and $\eps$. It was observed that
when $c$ is either too small or too large, the numerical results deteriorate and
an approximate kinetic function can no longer be associated with the numerical approximations.

\begin{figure}
\begin{center}
  \includegraphics[width=8cm,height=8cm]{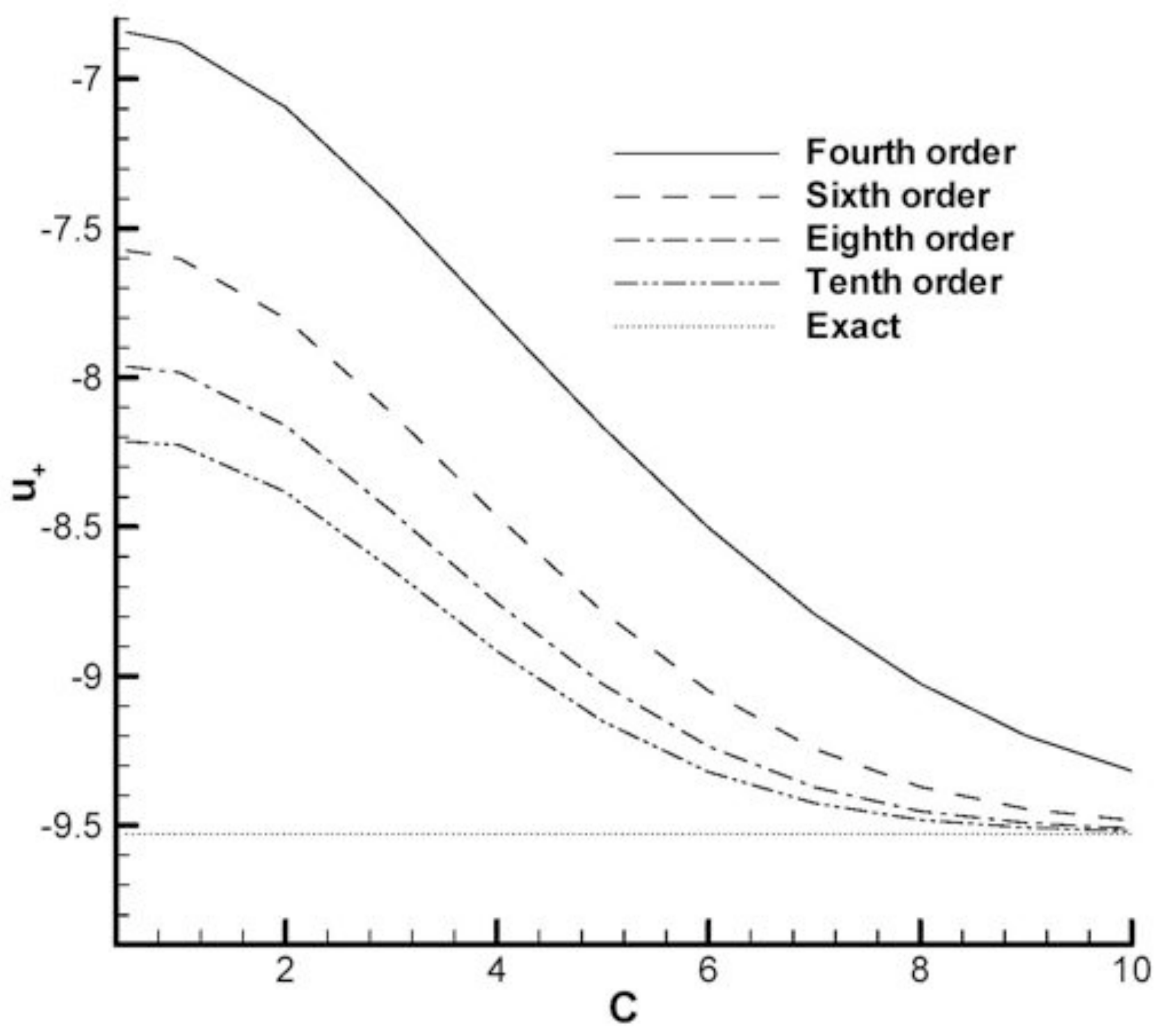}
  \\
  \caption{$u_+$ versus  $c$ for different schemes.}
 \label{c_coeff}
\end{center}
\end{figure}

\begin{figure} 
\begin{center}
\begin{tabular}{cc}
\begin{minipage}{6.0cm}
\begin{center}
\includegraphics*[height=6.5cm, width=6.5cm]{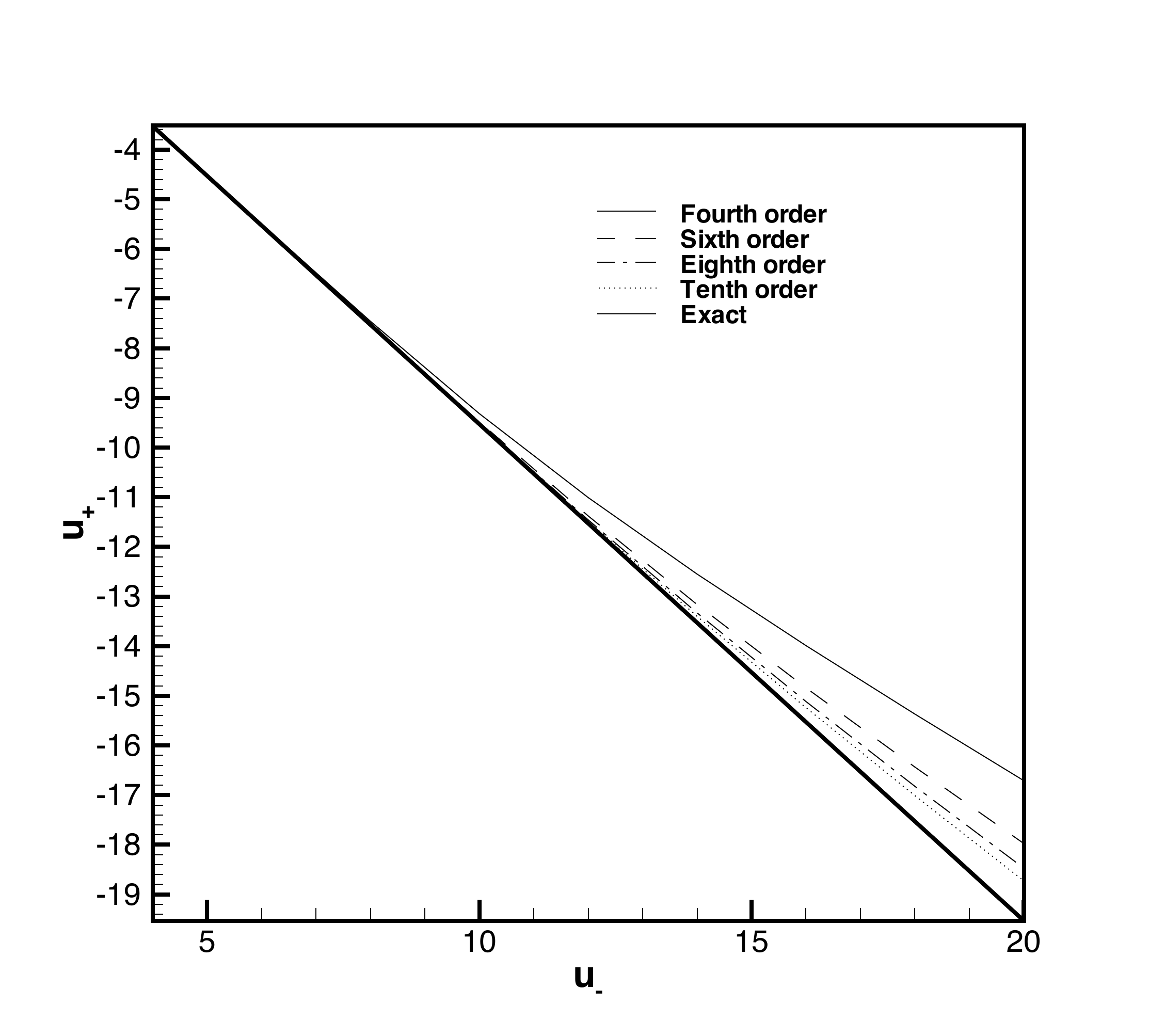}
\end{center}
\end{minipage}
&
\begin{minipage}{6.0cm}
\begin{center}
\includegraphics*[height=6.5cm, width=6.5cm]{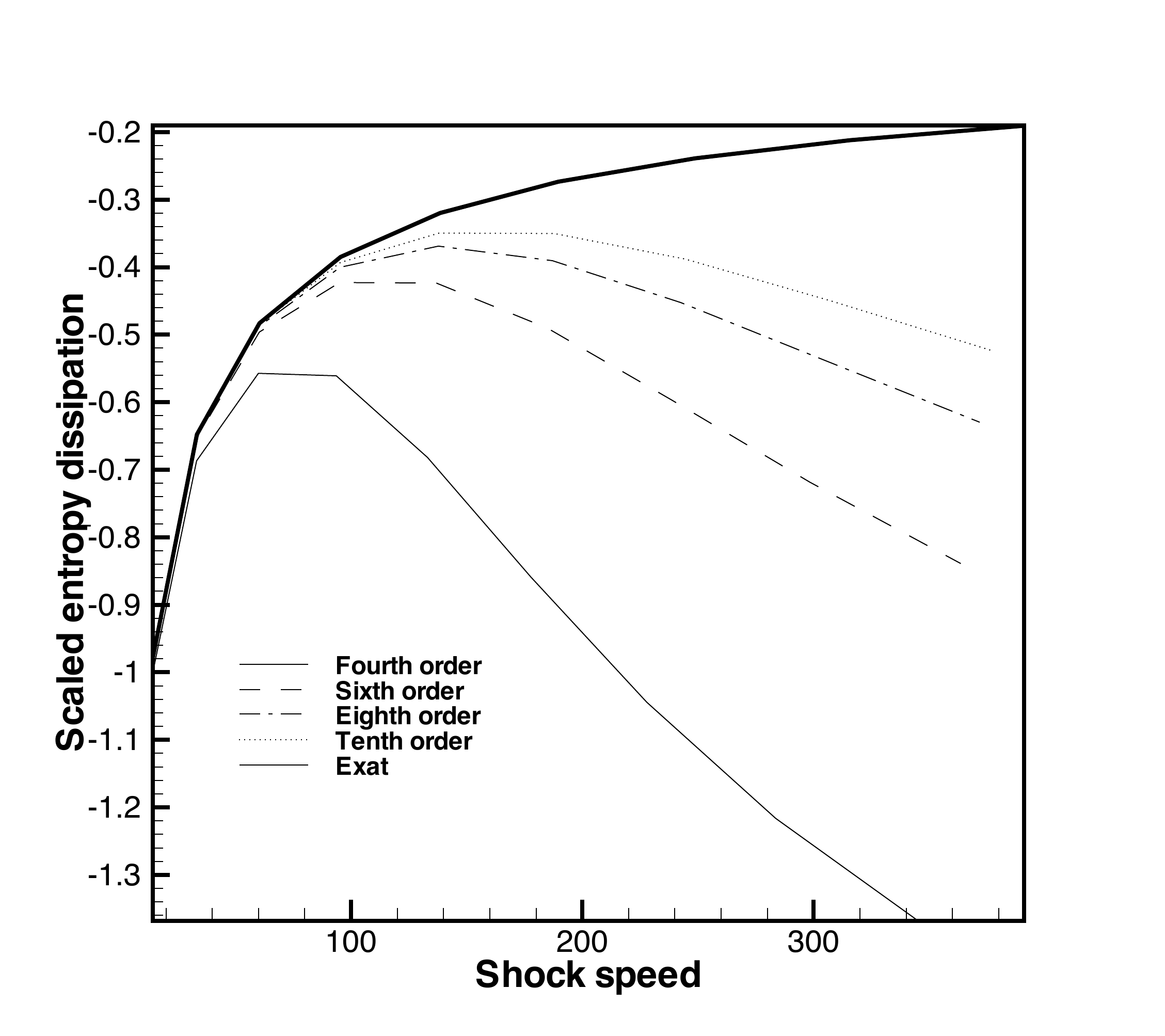}
\end{center}
\end{minipage}
\end{tabular}
\caption{$u_+$ versus $u_-$  (left) and scaled entropy dissipation $\phi(s)/s^2$ versus shock speed $s$ 
(right) for $\alpha=1$ and $c=10$.}
\label{c=10}
\end{center}
\end{figure}

As shown in Figure \ref{c_coeff} for typical values $u_-=10$,
$\alpha=1$ and CFL=0.5, the specific choice of $c$ does not change the general feature of
the results. However, as we increase the value of $c$, all schemes
converge to the exact solution and the tenth order scheme converges
faster. Moreover, it is concluded from Figure \ref{c_coeff} that the sixth order
scheme with $c=10$ is a reasonable choice in terms of accuracy and computational efficiency. For this value of $c$ (that is, $\epsilon=10h$) with $h=0.005$,
Figure \ref{c=10} shows the kinetic function $u_+$ versus $u_-$  (left) and the
corresponding scaled entropy dissipation $\phi(s)/s^2$ versus the
shock speed $s$ for $\alpha=1,4,6$ (right) which are in accordance with the previous
results about the convergence of the numerical results to the exact one.

\

\

\

\


\

\section{Fourth-order models}  
\label{KI-0}

We now return to the examples of Section~\ref{MO-0} and discuss the thin film and Camassa-Holm models. 
We will show that solutions of the Riemann problem may admit under-compressive
shocks, and we will determine the corresponding kinetic function in a significant range of data.  
Our results in the previous section indicate that we can continue this investigation
with six-order schemes, which we will do in (most of) the forthcoming tests.

\subsection{Thin liquid films}

Consider the model of thin liquid films together with the physical regularization terms. Our goal is to 
(numerically) compute a kinetic function associated with this model.
It should be reminded that the physically relevant range for the solutions is the interval $(0,1)$. 
This model has a flux with a single inflection point, as is the case with the cubic conservation laws, 
and the interesting aspect is the particular form of the regularization. In view of experiments performed in \cite{BertozziShearer},
it may be anticipated that a unique Riemann solver can not be associated to a given set of initial data, and
that the specific initialization adopted in implementing a given scheme may play a role. 
We will first show that a kinetic relation can be associated to a family of initial data and a given approximation of these initial data. 

Since the flux has a convex-concave shape rather than a concave-convex form (contrary to the cubic flux function), the nonclassical shock is fast undercompressive and it is more convenient
to draw the kinetic functions from right to left, that is, 
$$
u_- = \varphi (u_+).
$$
The calculations here are more delicate than those performed with the third-order regularization,
since we are now dealing with a fourth-order regularization. Therefore, for the discretization of the terms 
arising in \eqref{AP.1}, a scheme must be fifth order in accuracy at least. 
We thus exclude the fourth-order scheme used earlier and we concentrate attention on the sixth-order scheme. 

It should be mentioned that the (fourth-order) regularization is
singular at $u=0$, so tests involving values close to $u=0$ are
challenging and the kinetic function should have some degenerate behavior as some
values in the solution approache $u=0$.

The thin film model (\ref{AP.1}) may be written in the conservative form
\be
u_t+g(u)_x =0,
\label{TF.1}
\ee
with
\be
g(u):=u^2-u^3-u^3 \, ( \delta  \, u_x  - \, u_{xxx} ),
\label{TF.2}
\ee
where we have set $\eps = \delta$ and $\alpha \eps^2 =1$. 

The numerical solution of (\ref{TF.1}) is obtained by calculating the terms $u_x$
and $u_{xxx}$ in (\ref{TF.2}) using the sixth order scheme (Appendix A), computing
the nodal values of $g$, and finally calculating $g_x$
again using the sixth order scheme (Appendix A).
The numerical experiments are performed here with
 $ \delta=0.1h $ and $h=1$ using a computational grid
 of $1000$ nodes and the initial condition used in \cite{BertozziShearer}
\be
u(x)=(\tanh(-x+100)+1){(u_L-u_R)\over 2}+u_R.
\label{Ini.1}
\ee

\begin{figure}
\begin{center}
  \includegraphics[width=8cm,height=8cm]{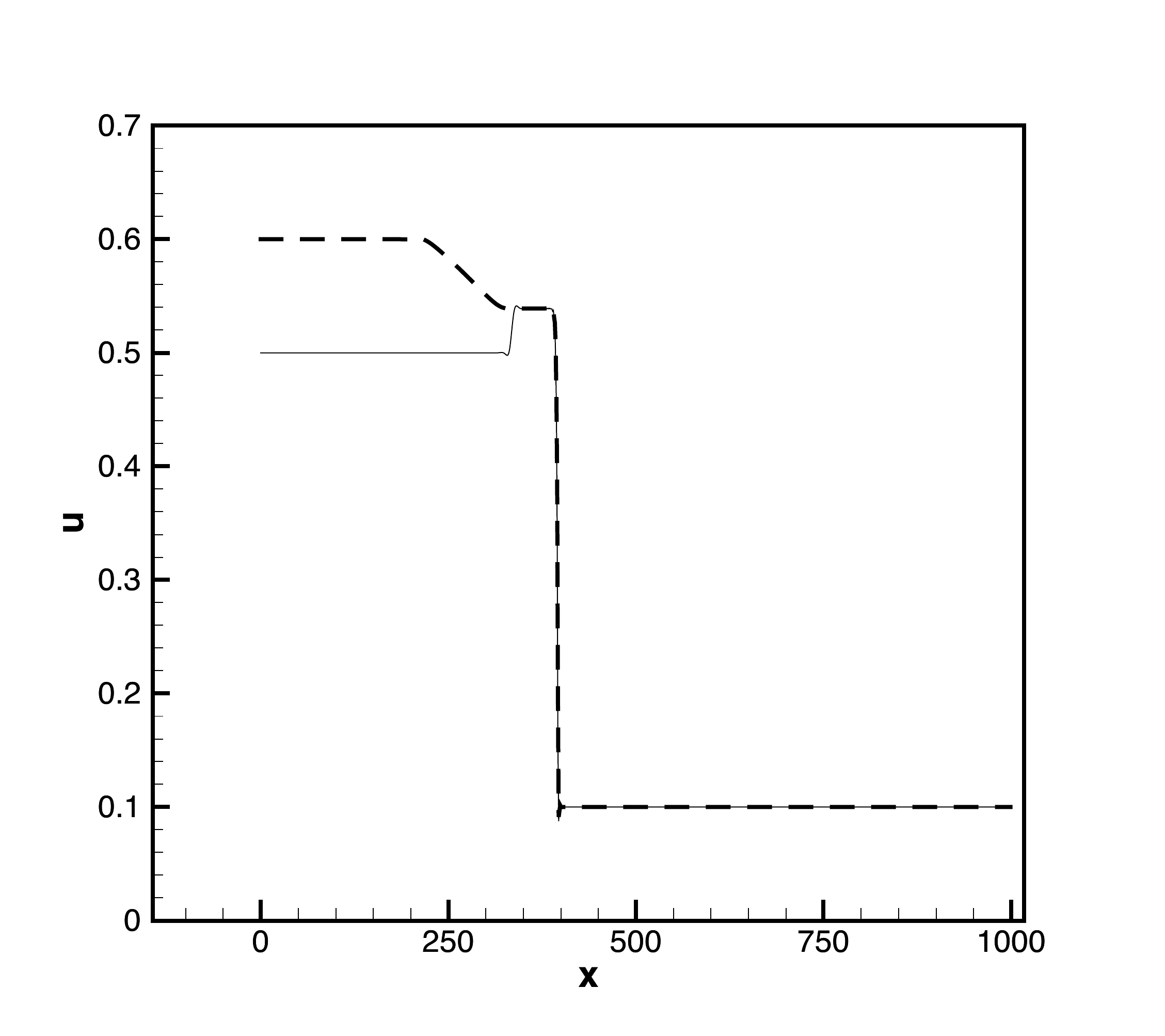}
  \caption{Wave structure for $u_R=0.1$ and $u_L=0.5$ and $0.6$.}
 \label{u_R=0.1}
\end{center}
\end{figure}

The numerical solutions of (\ref{TF.1}) for $u_R=0.1$ and $u_L=0.5$ and
 $0.6$ with $\Delta t=0.6481 $ at time $t=1037$ are shown
 in Figure \ref{u_R=0.1}. For the case $u_L=0.5$ a 
 double shock structure is observed while for $u_L=0.6$ a rarefaction wave and an under-compressive shock are generated. The speed of the
 non-classical shock is similar for the two cases, as expected.

\begin{figure}
\begin{center}
  \includegraphics[width=8cm,height=8cm]{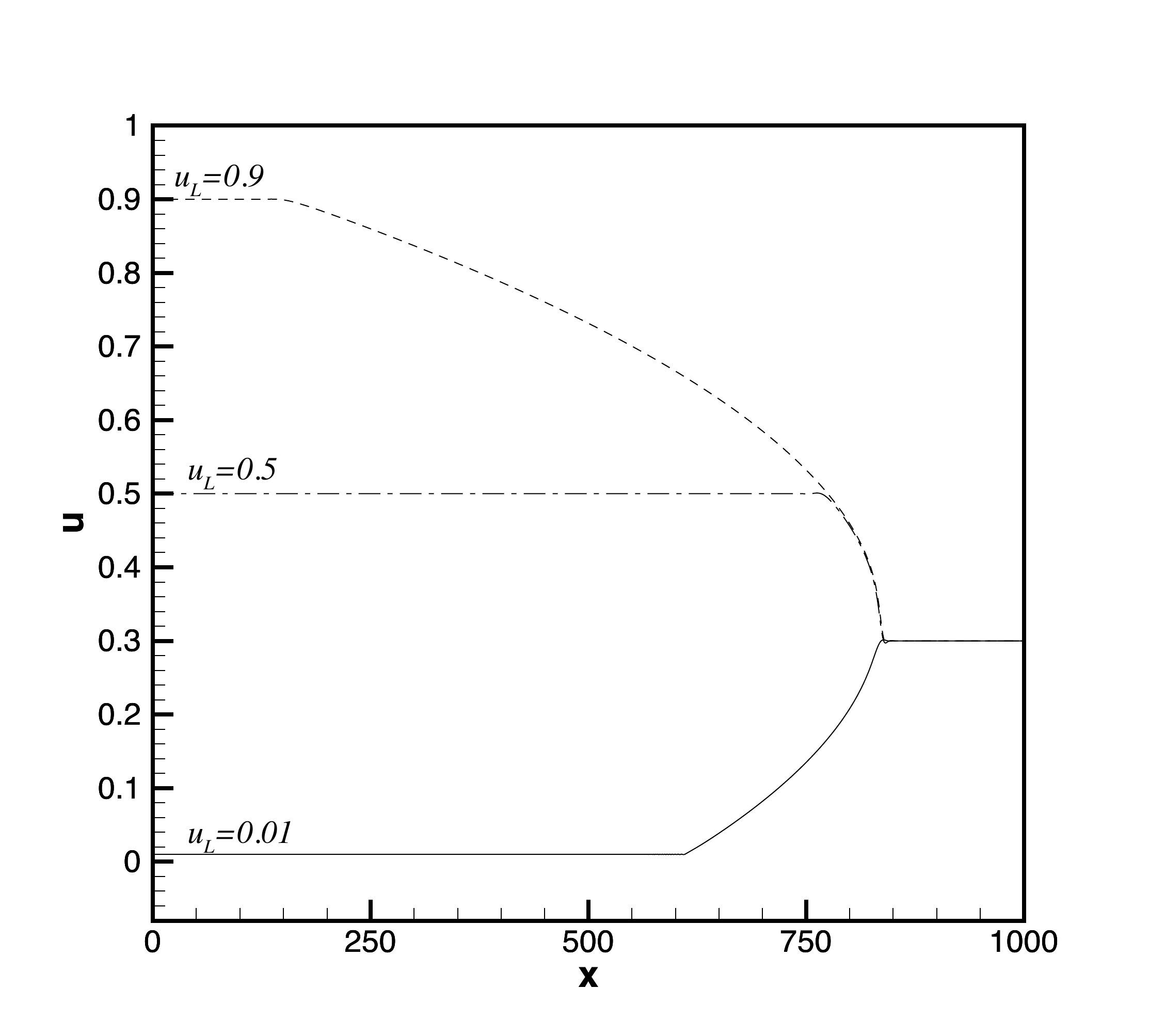}
  \caption{Wave structure for $u_R=0.3$ and various values of $u_L$.}
 \label{u_R=0.3}
\end{center}
\end{figure}

\begin{figure}
\begin{center}
  \includegraphics[width=8cm,height=8cm]{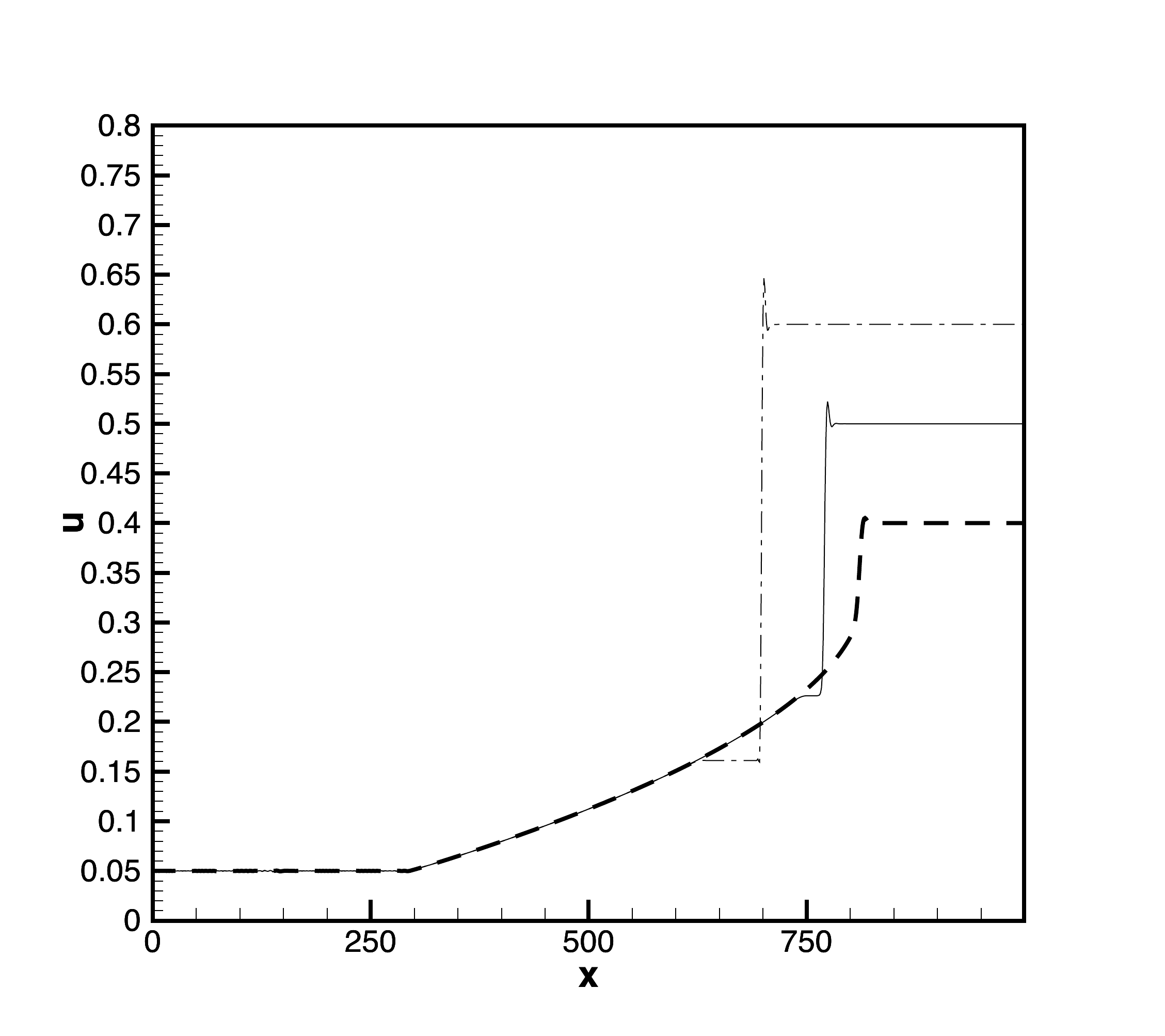}
  \caption{Wave structure for $u_L=0.05$ and $u_R=0.4$, 0.5 and $0.6$.}
 \label{u_R=0.4-0.5}
\end{center}
\end{figure}

In agreement with what was observed (and proven rigorously in \cite{BedjaouiLeFloch}) 
for the diffusive-dispersive regularizations,  when $0.2<u_R<0.4$, i.e.~around the inflection
point $u=1/3$,  only {\sl classical} waves are observed as shown e.g. in Figure \ref{u_R=0.3}
for $u_R=0.3$ and different values of $u_L$  (with $\Delta t=0.388$ at time $t=698$).
For $0.4<u_R$ {\sl non-classical} waves arise even when $u_L$ is small 
and become clearly visible as $u_R$ increases as shown in Figure~\ref{u_R=0.4-0.5}.

\begin{figure}
\begin{center}
  \includegraphics[width=8cm,height=8cm]{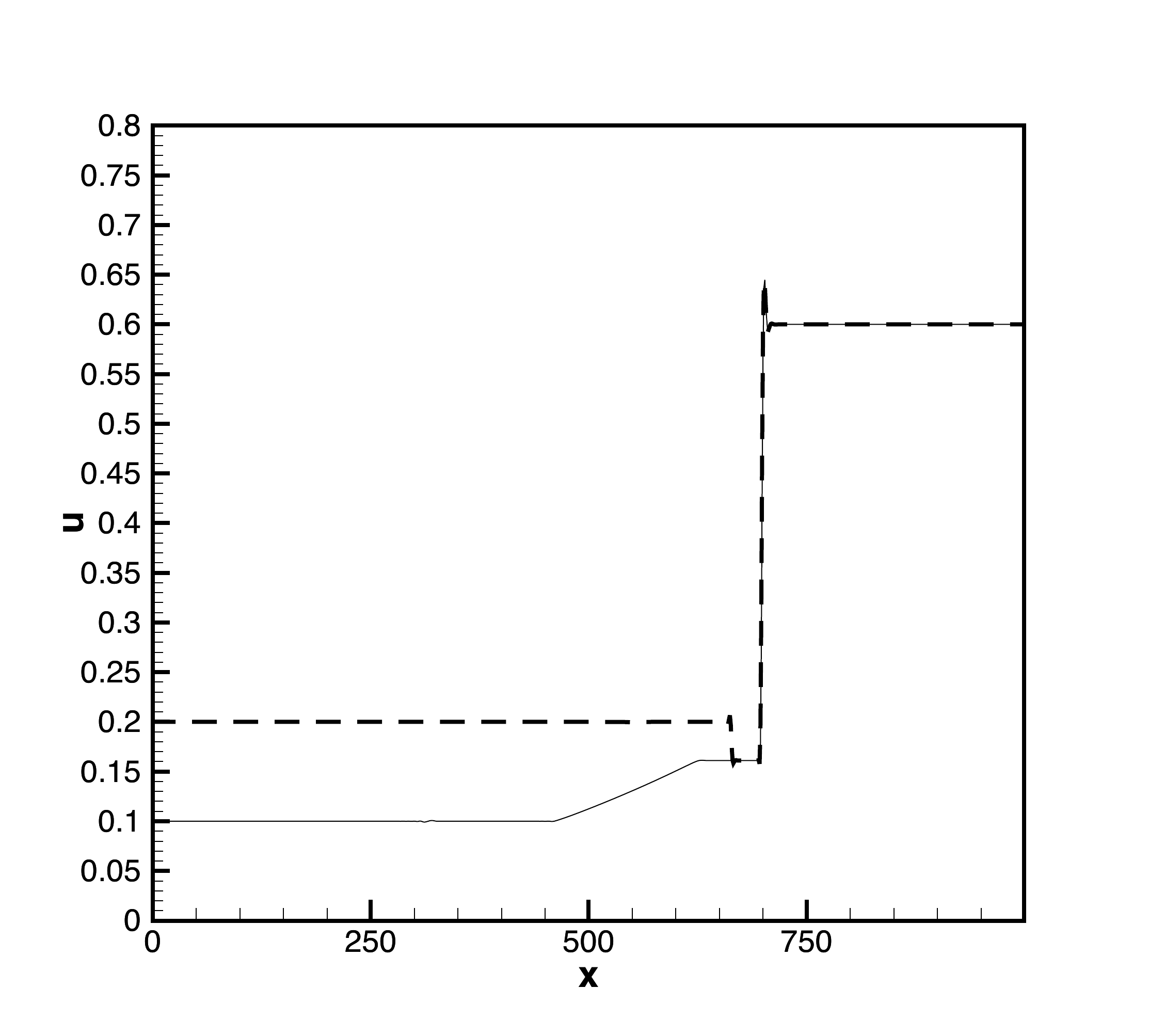}
  \caption{Wave structure for $u_R=0.6$ and $u_L=0.1$ and $0.2$.}
 \label{u_R=0.6}
\end{center}
\end{figure}

\begin{figure}
\begin{center}
  \includegraphics[width=8cm,height=8cm]{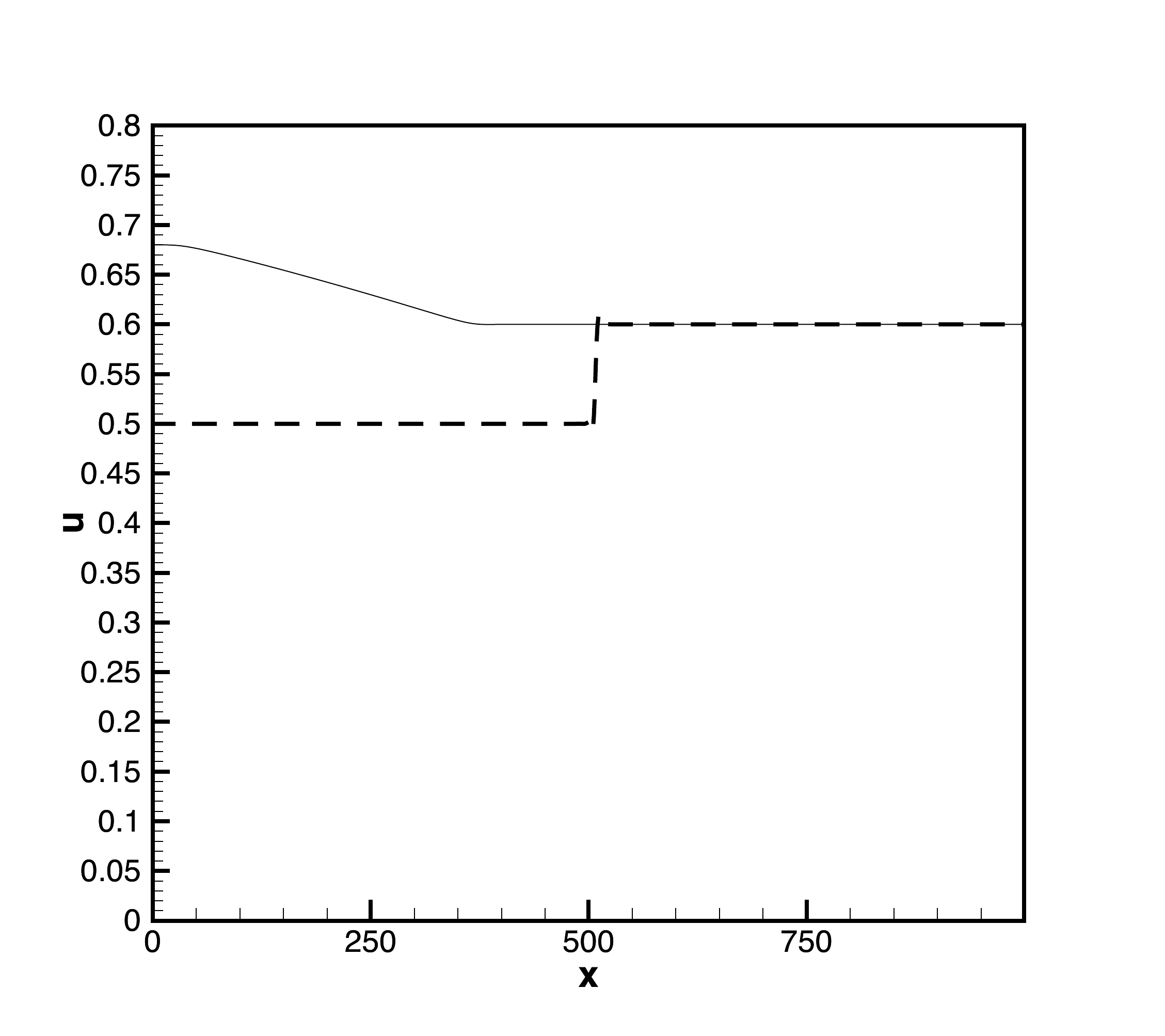}
  \caption{Wave structure for $u_R=0.6$ and $u_L=0.5$ and $0.68$.}
 \label{ur=0.6-II}
\end{center}
\end{figure}

Figure \ref{u_R=0.6}  shows the numerical solutions
 obtained for $\Delta t=0.37 $, at time $t=2148$
for $u_R=0.6$ and two cases $u_L=0.1$ and $u_L=0.2$. For the
case $u_L=0.2$ a double shock structure
is observed while for $u_L=0.1$ a rarefaction wave and
an under-compressive shock are generated. It should be mentioned that when $u_L$ is selected large enough, only
classical  waves are observed as shown for $u_L=0.5$
and $0.68$ (with the same data) in Figure~\ref{ur=0.6-II}.

\begin{figure} 
\begin{center}
\begin{tabular}{cc}
\begin{minipage}{6.0cm}
\begin{center}
\includegraphics*[height=6.5cm, width=6.5cm]{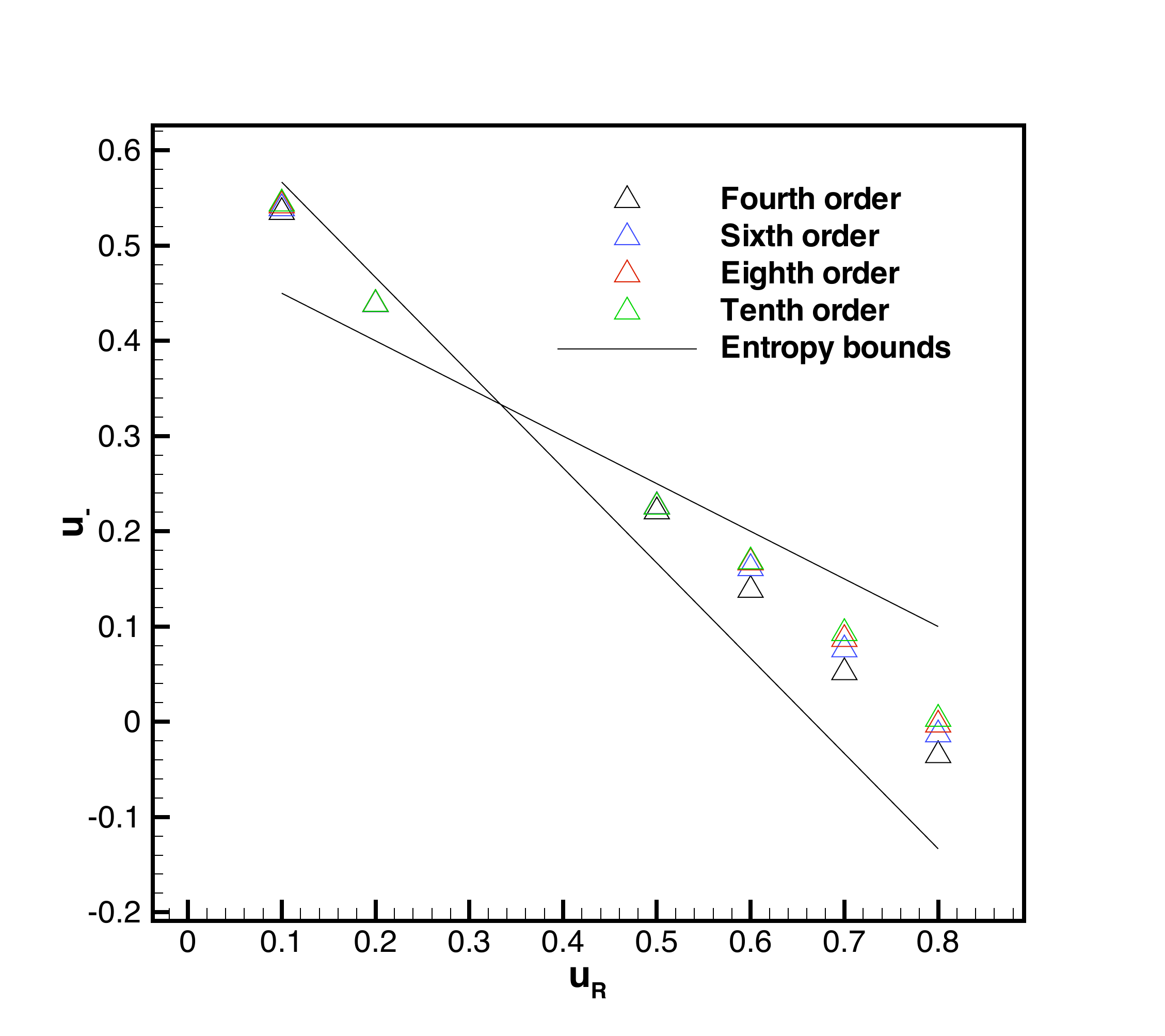}
\end{center}
\end{minipage}
&
\begin{minipage}{6.0cm}
\begin{center}
\includegraphics*[height=6.5cm, width=6.5cm]{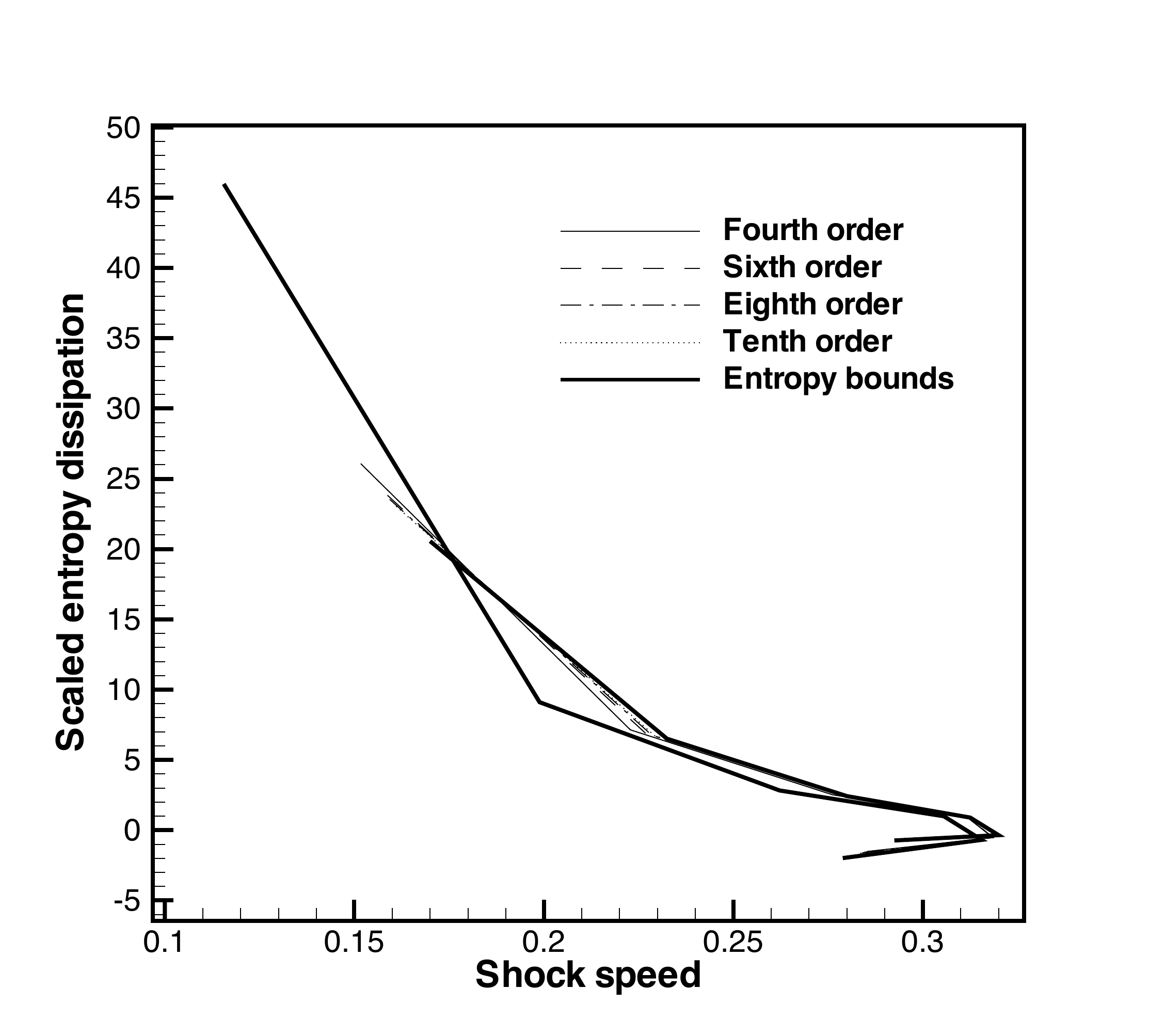}
\end{center}
\end{minipage}
\end{tabular}
\caption{Kinetic function (left) and scaled entropy dissipation (right), for $\eta=0.1$  based on the initial conditions (\ref{Ini.1}).}
\label{FigkineticThin}
\end{center}
\end{figure}

Varying the left-hand and right-hand initial states and keeping fixed all other parameters,
for $\delta=\eta h$ with $\eta=$ 0.1 we obtain the numerical kinetic function ($u_-$ as a function of $u_R$) 
shown in Figure \ref{FigkineticThin} (left). 
Interestingly enough, this function is decreasing, as required in the general theory of nonclassical shocks. It is also
observed that the numerical kinetic functions converge as the order of the accuracy of the numerical scheme is increased.
Note that, for sufficiently large right-hand values of $u_R$, negative left-hand values $u_-$ are obtained which are physically irrelevant. The scaled entropy dissipation $\phi(s)/s^2$ versus the shock speed $s$ is also shown in Figure \ref{FigkineticThin} (right). 

\begin{figure}
\begin{center}
  \includegraphics[width=8cm,height=8cm]{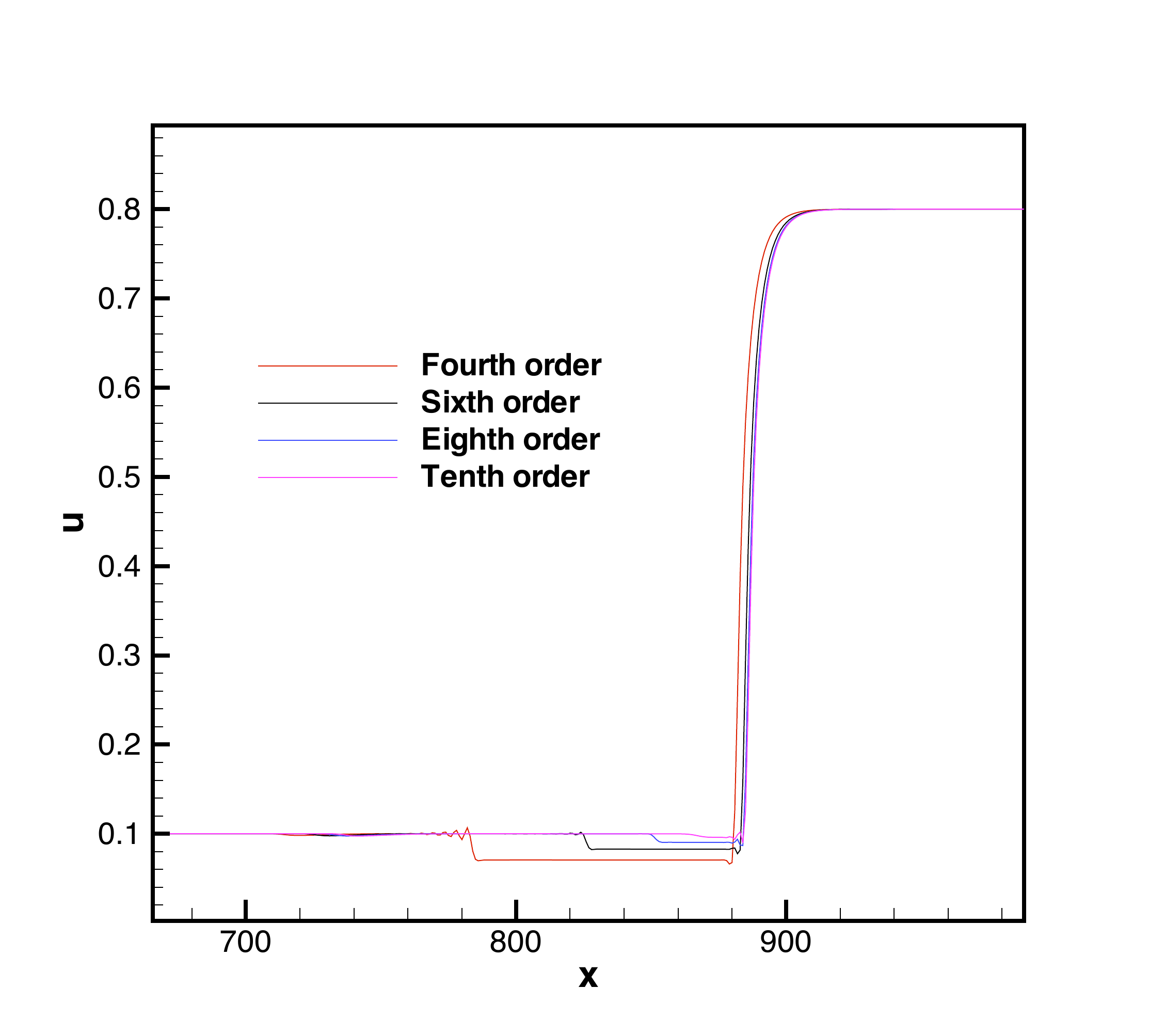}
  \caption{ A sample of the solution with various schemes for $u_r=0.8$ and $\eta=4.7$. }
 \label{sample}
\end{center}
\end{figure}

\begin{figure}
\begin{center}
  \includegraphics[width=8cm,height=8cm]{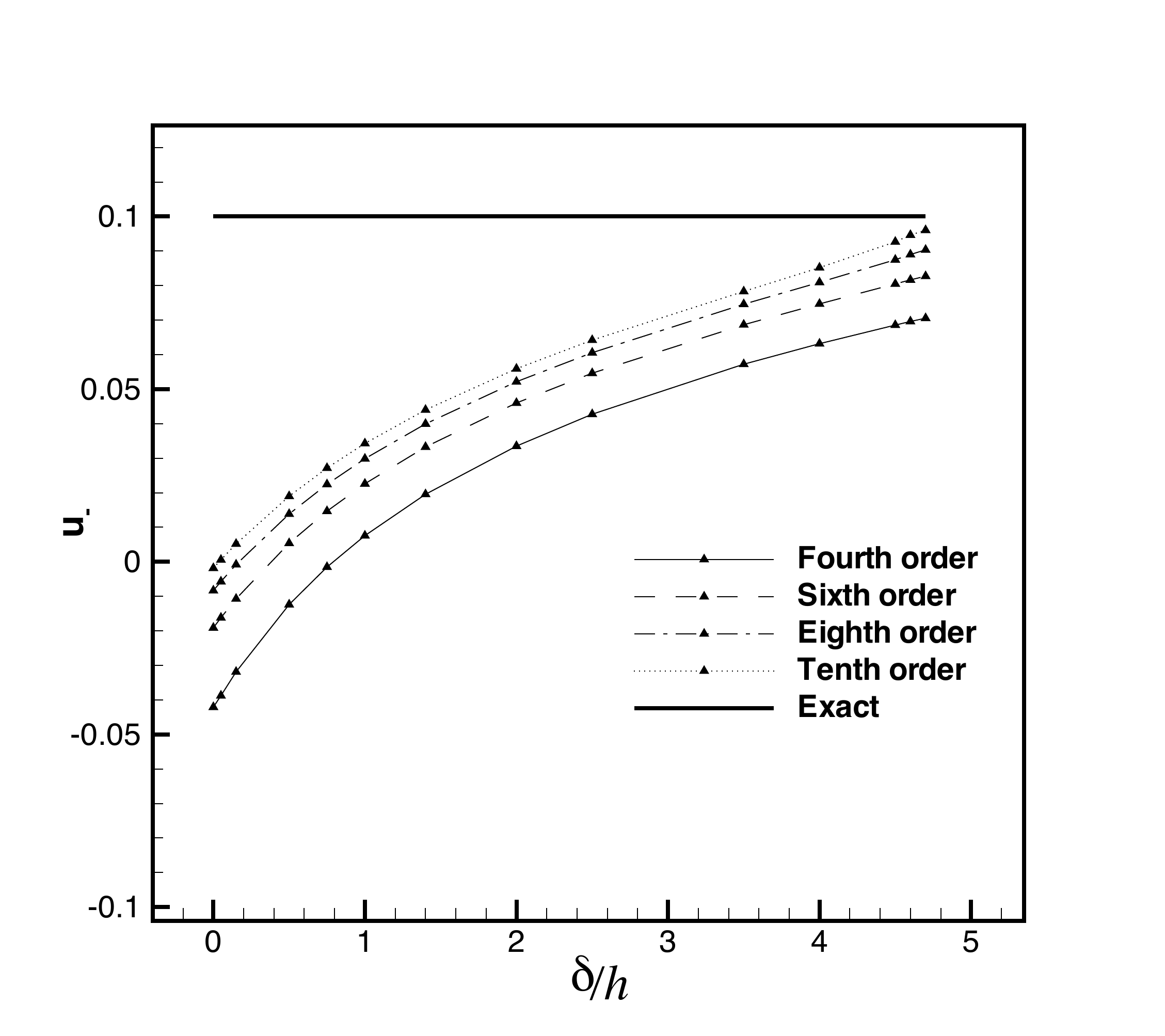}
  \caption{$u_-$ versus $\eta=\delta /h$ for $u_r=0.8$.}
 \label{convergence}
\end{center}
\end{figure}

\begin{figure}
\begin{center}
  \includegraphics[width=8cm,height=8cm]{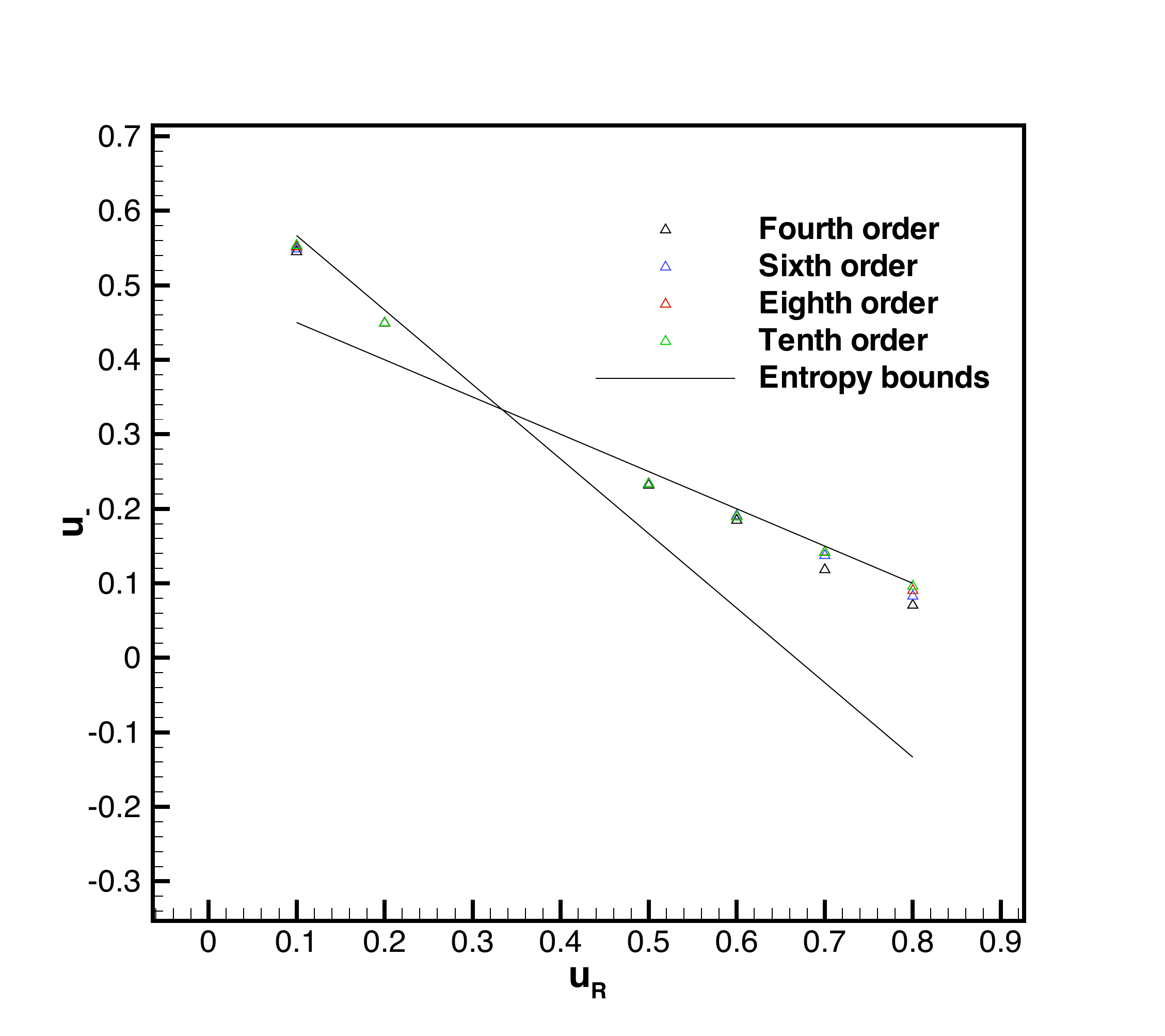}
  \caption{Kinetic function for thin films, based on the initial conditions (\ref{Ini.1}) with optimum $\eta$.}
 \label{KineticVariable}
\end{center}
\end{figure}

Finally, we study the effect of  $ \eta=\delta /h $ for $u_r=0.8$. Figure \ref{sample} presents a sample of the solution with various schemes for $\eta=4.7$. In Figure \ref{convergence}, it is shown that by increasing $\eta$, the
numerical results converge to the limiting value with the tenth order scheme.  The optimum value 
of $\eta$ depends on $u_r$ and it is increased as $u_r$ does (not shown). In Figure \ref{KineticVariable} we have plotted the
kinetic functions for various schemes where a variable (optimum) $\eta$ (depending on $u_r$) is employed. 
The convergence 
predicted in Conjecture~\ref{conjecture} is clearly confirmed here.

\begin{figure} 
\begin{center}
\begin{tabular}{cc}
\begin{minipage}{7.0cm}
\begin{center}
\includegraphics*[height=7.0cm, width=7.0cm]{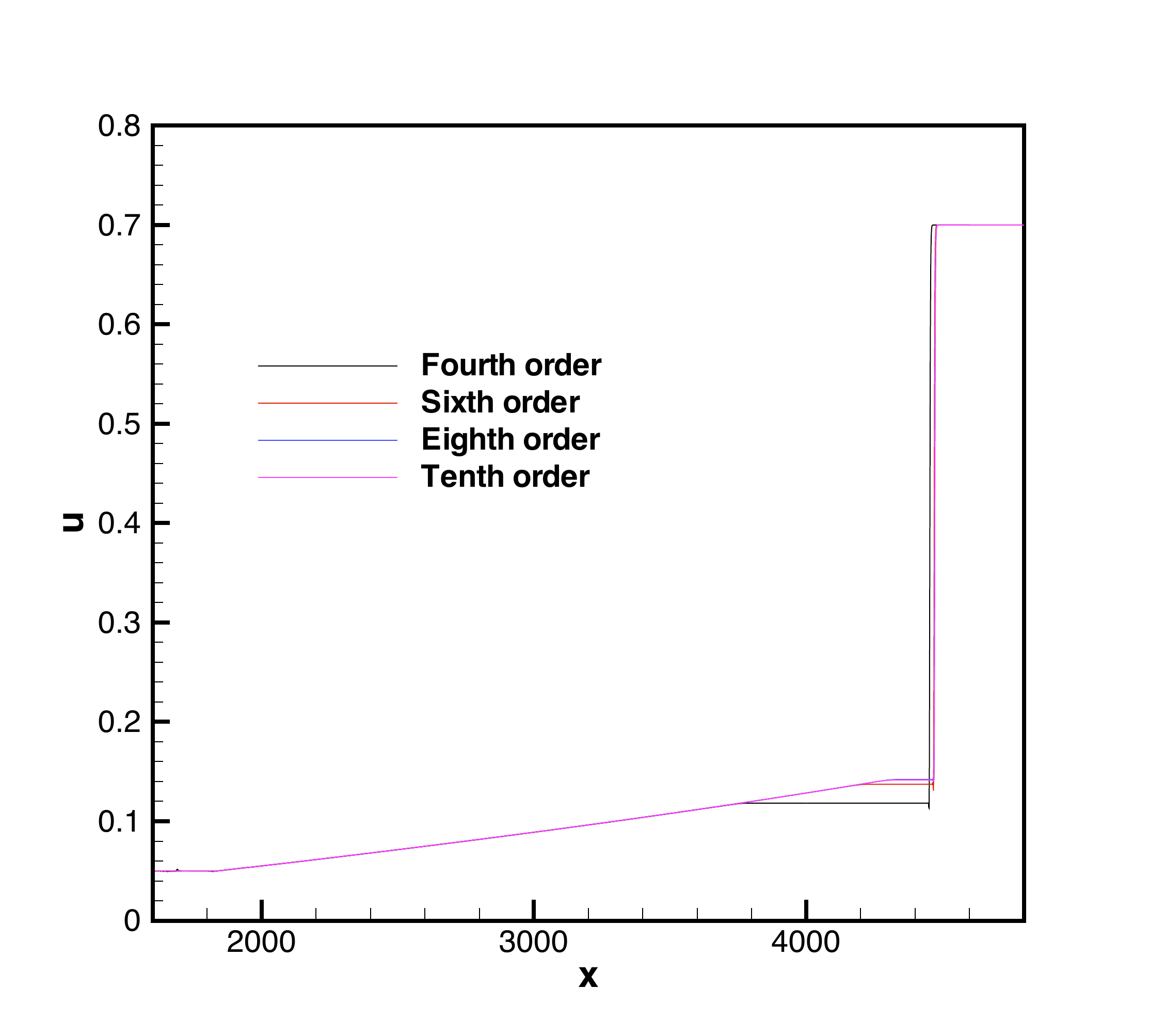}
\end{center}
\end{minipage}
&
\begin{minipage}{7.0cm}
\begin{center}
\includegraphics*[height=7.0cm, width=7.0cm]{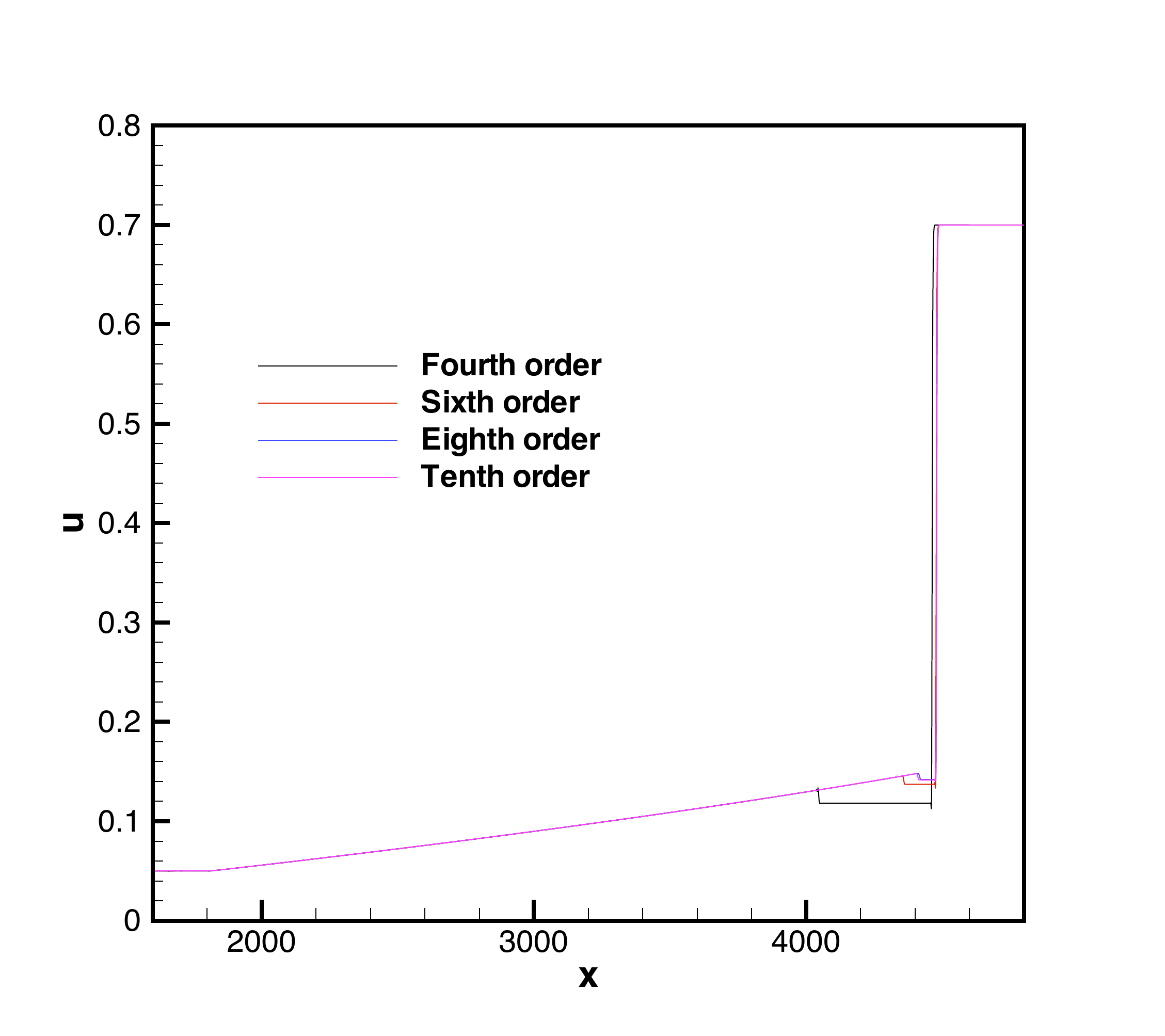}
\end{center}
\end{minipage}
\end{tabular}
\caption{Numerical results obtained using the initial conditions (\ref{Ini.1}) (left) and (\ref{Ini.2})(right).}
\label{figIni}
\end{center}
\end{figure}

In order to investigate the effect of the initial conditions, we consider the following initial condition
 \be
u(x)=
 \begin{cases}
(0.35-u_L)\tanh(x-80 )/2+(0.35+u_L)/2,     &         x \leq 105,
    \\
(u_R-0.35)\tanh(x-130)/2+(0.35+u_R)/2,    &         x >  105.
 \end{cases}
\label{Ini.2}
\ee
The numerical solutions of (\ref{TF.1}) for $u_R=0.7$ and $u_L=0.05$ 
 with $\Delta t=0.27 $ and $\Delta x=1 $ at time $t=18792$ are shown
 in Figure~\ref{figIni}. As shown in Figure~\ref{figIni}, the two choices of numerical solutions lead to different
 results. However, such a dependence on the initial data is only limited to few cases.
%


\subsection{Generalized Camassa-Holm model}

It is expected that the kinetic function of the Camassa-Holm model is qualitatively similar to that of the
diffusive-dispersive model. However, some differences, at least quantitatively are expected. 
In the following, we draw a kinetic function for the Camassa-Holm model 
and compare with the diffusive-dispersive model with the same coefficients.
First, we investigate the behavior of the function for large values of $u$.

\begin{figure} 
\begin{center}
\begin{tabular}{cc}
\begin{minipage}{6.0cm}
\begin{center}
  \includegraphics[width=6.5cm,height=6.5cm]{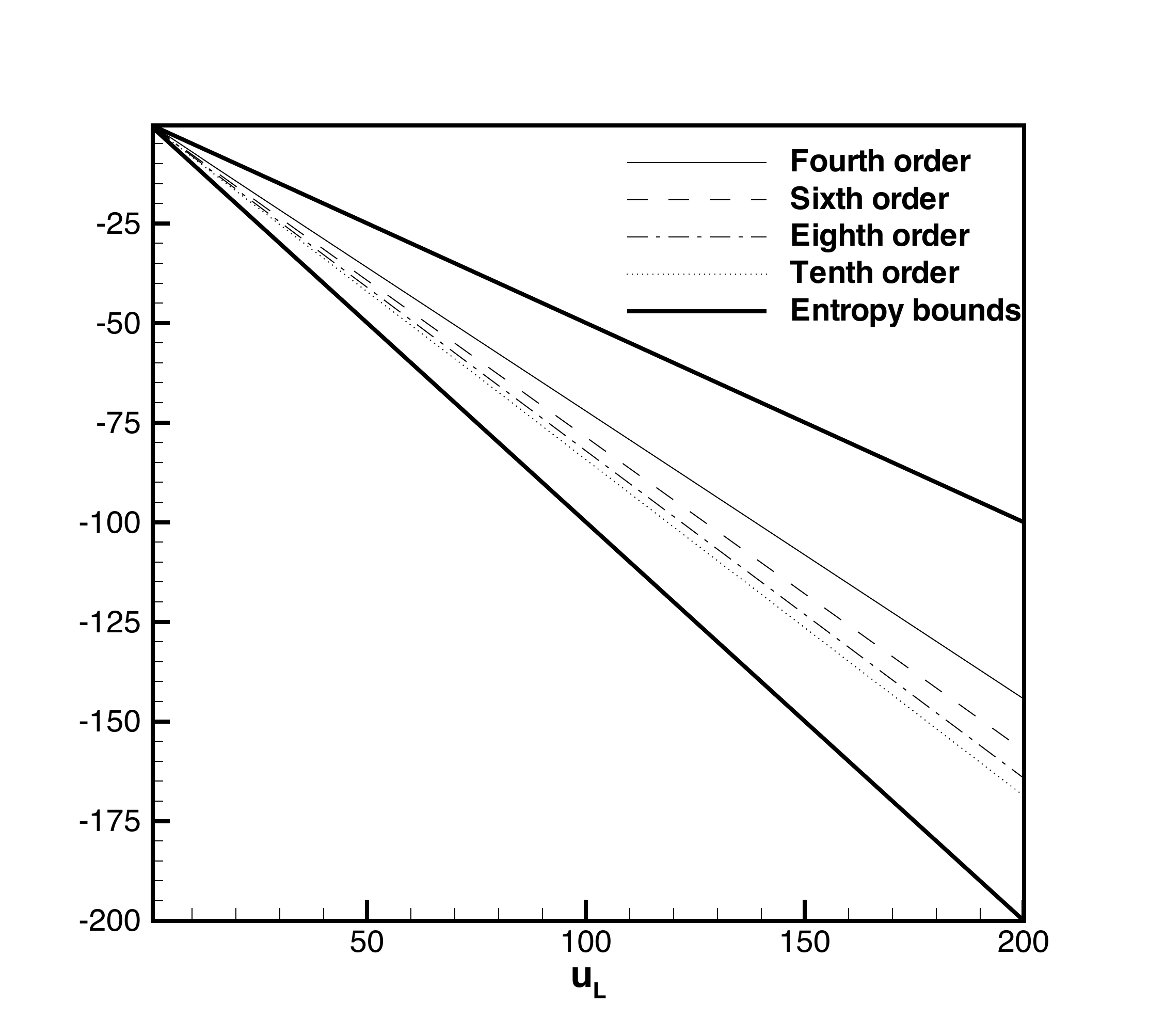}
\end{center}
\end{minipage}
&
\begin{minipage}{6.0cm}
\begin{center}
\includegraphics*[height=6.5cm, width=6.5cm]{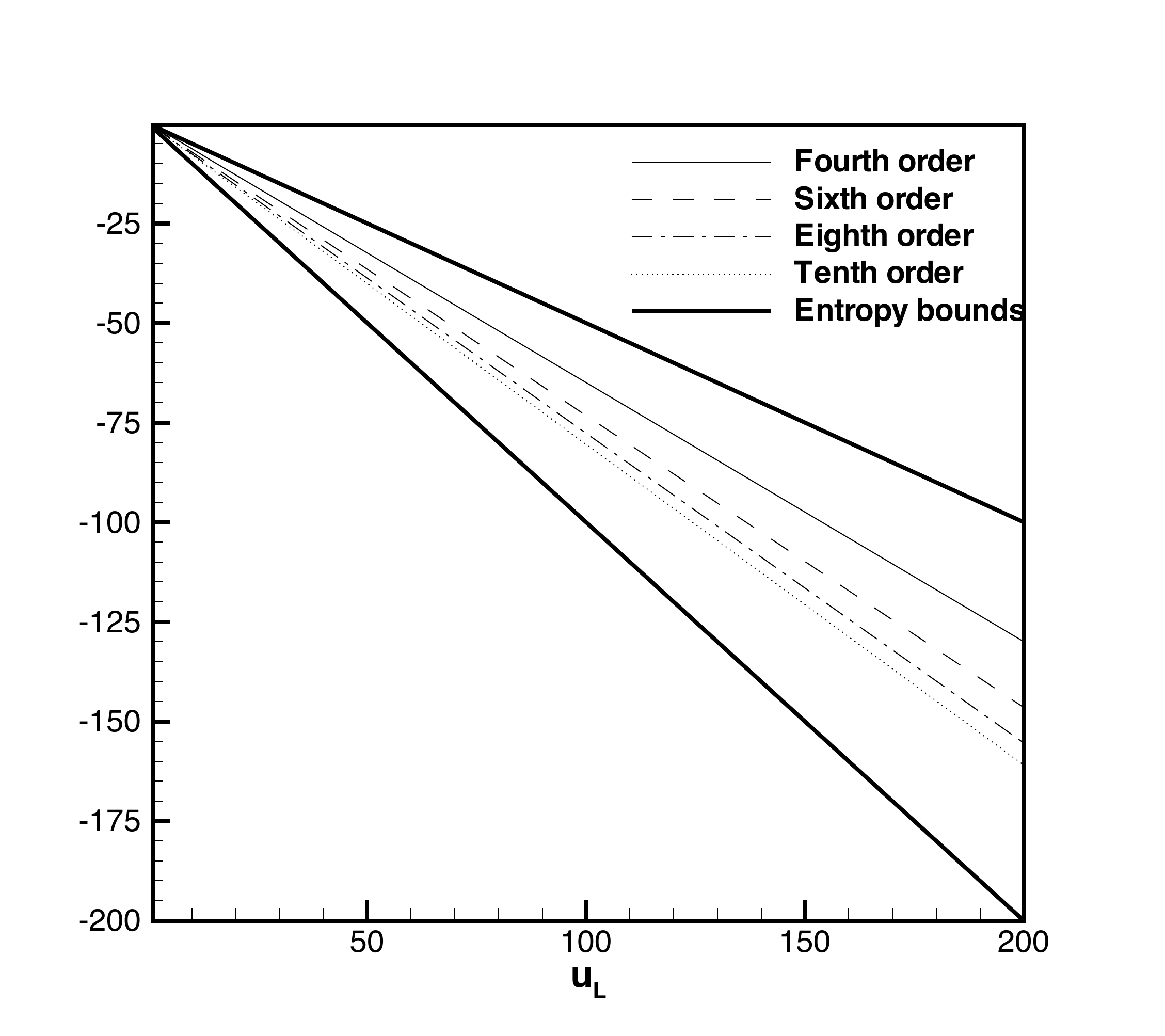}
\end{center}
\end{minipage}
\end{tabular}
\caption{The kinetic function (large values of $u$) for the linear diffusion-dispersion case (left) and Camassa-Holm model (right).}
\label{Camassa1}
\end{center}
\end{figure}

Figure \ref{Camassa1} shows the results obtained with $h=1$,\ $\eps=0.005h$,\ $\alpha=1$ and CFL=0.5 using a computational grid of $1000$ nodes. 
The kinetic function for the Camassa-Holm model (Figure~\ref{Camassa1}, right) is larger than the linear diffusion-dispersion case (Figure~\ref{Camassa1}, left).

The comparison between the Camassa-Holm model and the linear diffusion-dispersion model
allows to see that 
the kinetic function of the two models is similar for {\sl small} values of  
$u$ but {\sl different} for large values of $u$. Typically, 
for large values of $u$, the Camassa-Holm model leads to larger kinetic values.
This behavior is due to the nonlinearity of the regularization in  the right-hand side of \eqref{AP.camassa}. 
It would be interesting (and challenging) 
to apply the techniques in \cite{BedjaouiLeFloch} and 
prove the existence and monotonicity of the kinetic function associated with the Camassa-Holm model.

\begin{figure}
\begin{center}
  \includegraphics[width=6.5cm,height=6.5cm]{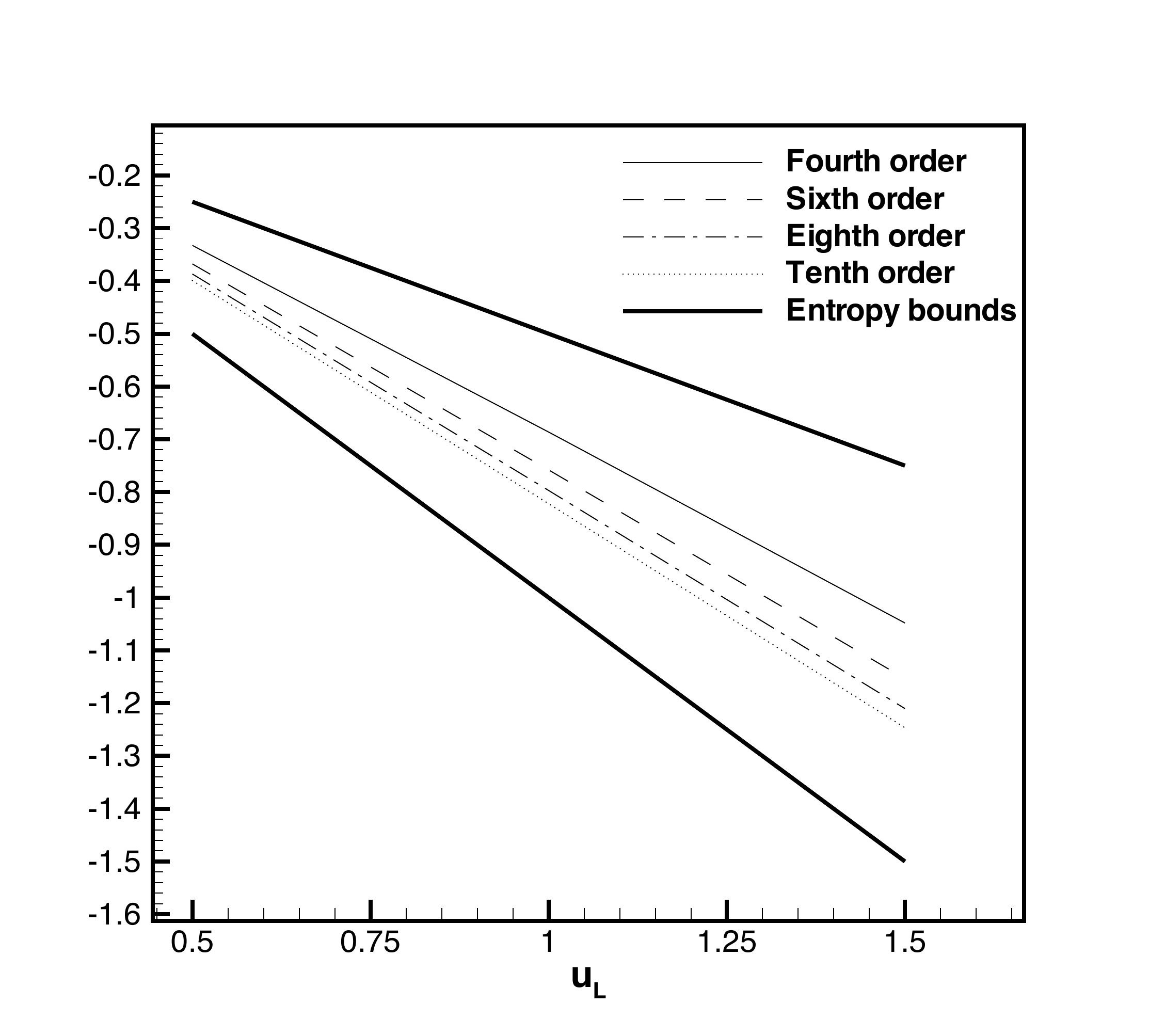}
  \caption{The kinetic function (small values of $u$) for the linear diffusion-dispersion case and Camassa-Holm model.}
 \label{Camassa2}
\end{center}
\end{figure}

We then consider small values of $u$. Figure \ref{Camassa2} shows the results obtained for both models
with $h=1$,\ $\eps=0.1h$,\ $\alpha=4$ and CFL=0.5 using a computational grid
 of $1000$ nodes. Contrary to the case of large values of $u$, the behavior of the function for small values of $u$ is the same for both models.

 \


\

\section{Kinetic functions associated with van der Waal fluids}
\label{MU-0}

Our next objective is to investigate the properties of non-classical
solutions to hyperbolic conservation laws whose flux-function admits two inflection points.
We consider the phase transition model \eqref{VD.1} presented in Section~\ref{MO-0}. 
As we will see, the dynamics of non-classical shock waves for general flux
is much more intricate than the now well-understood case of a single inflection point.
In particular, we demonstrate numerically that the Riemann problem
may well admit {\sl arbitrarily large number} of solutions.
This happens only when the flux admits two inflection points at least and provided the dispersion
parameter is sufficiently large. We then determine the kinetic function and observe its lack of monotonicity.

We consider an equation of state having the same shape as that of van der Waals fluids,
described by the (normalized) equation
\be
p(\tau) :=  { 1 \over (3\tau-1)^{1+1/\zeta} } - {3 \over \tau^2}, \quad u > 1/3,
\label{MU.1}
\ee
for some positive constant $\zeta = 1/(\gamma - 1)$ where $\gamma \in (1,2)$.
 
To exhibit the dynamics of non-classical shocks for general flux, we solve the Riemann problem 
for left- and right-data:
\begin{enumerate}
\item First, we determine the wave structure of the solution for each
fixed left-hand state, that is, we identify the waves (classical/non-classical shock or rarefaction) within the
Riemann solution. In turn, by varying the left-hand state, we identify regions in the plane
in which the structure of the Riemann solution remain unchanged as we change the Riemann data.
This provides us with a representation of the Riemann solver in the plane.
\item Second, within the range of $u_l$ where non-classical shocks leaving from $u_l$
are available, we determine the corresponding kinetic functions that are
needed to characterize the dynamics of non-classical shocks. It is expected that
more than one kinetic function will be needed in some range of $u_l$, but
a single kinetic function maybe sufficient in certain intervals.
\end{enumerate}

\begin{figure}
\begin{center}
  \includegraphics[width=8cm,height=8cm]{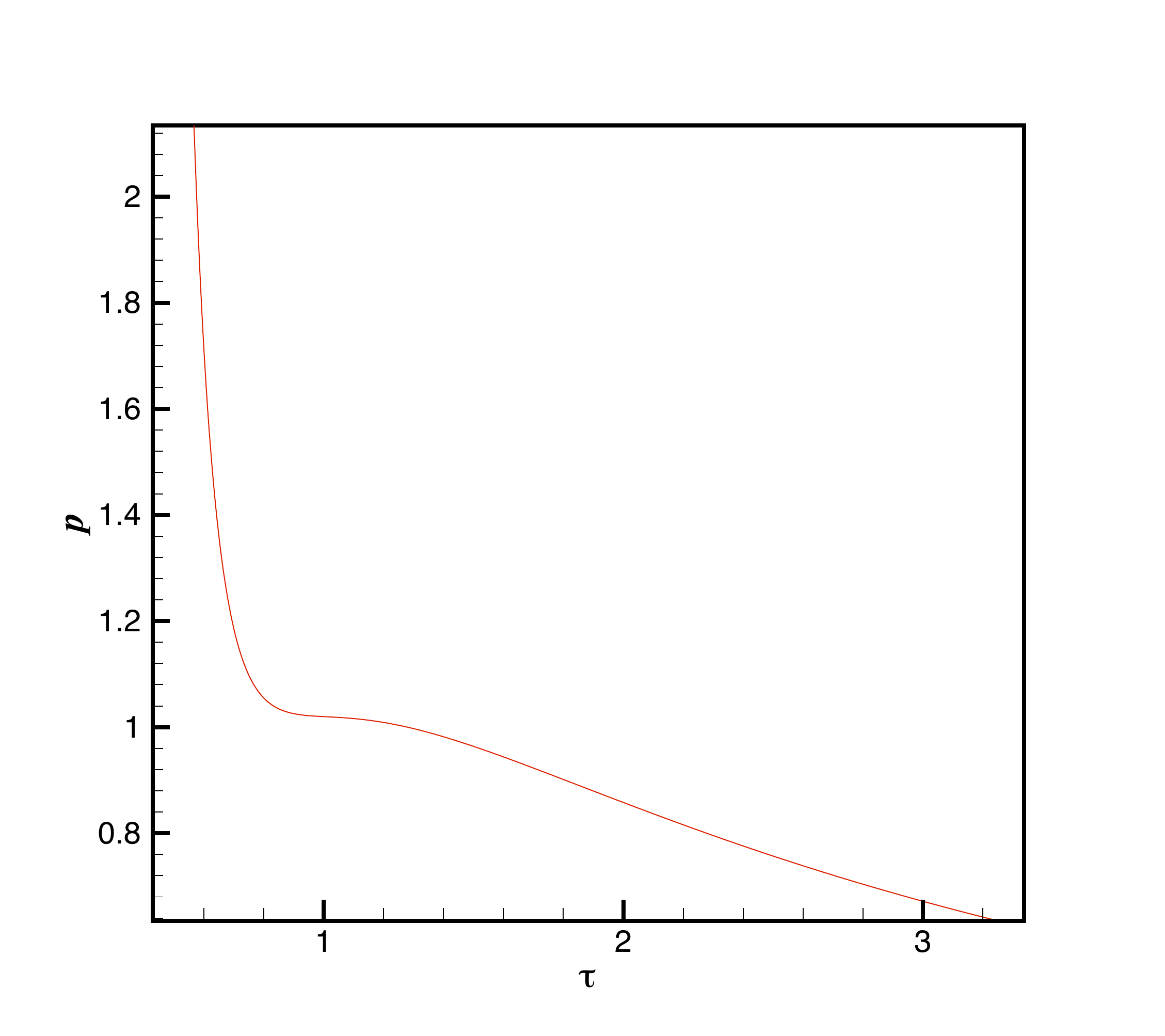}
  \caption{Equation of state considered in (\ref{MU.2}).}
 \label{flux}
\end{center}
\end{figure}

To begin with, we consider a simplified case
\be
\begin{split}
& \del_t \tau - \del_x u = 0, \\
& \del_t u + \del_x p(\tau) = \eps \, \del_{xx}  u,
\end{split}
\label{VD.smpl} 
\ee 
with
\be
 p(\tau) :=  { RT \over (\tau- \frac13)} - \displaystyle{3 \over \tau^2},
\label{MU.2}
\ee
and $R=\frac83$ and $T=1.005$. This flux function has two inflection points at $\tau=1.00996$ and 1.8515 and
 it is shown in Figure  \ref{flux}. 
 
 In order to investigate possible wave cases, various test cases are performed by keeping three of
 four variables ($u_L$, $\tau_L$, $\tau_R$ and $u_R$) constant and changing the fourth one. Some typical wave structures for $u_L=0.35$, $\tau_L=0.8$, $\tau_R=2$ and $u_R=0.5, 1.5$ and 2 are shown in Figures~\ref{vanderwaals0.5}, \ref{vanderwaals1.5} and \ref{vanderwaals2}. 

Next, we select $u_R=1$, $\tau_L=0.8$, $\tau_R=2$ and $u_L$ variable from 0 to 2. Note that the two inflection points are between the selected  $\tau_L$ and $\tau_R$.
We study the kinetic function for this case.  Numerical results are obtained using a computational grid of $2000$ nodes and $h=0.0005$ with $\eps=0.00003$. The wave structures for this case are given in Figures~\ref{tau} and \ref{u} (corresponding to the case $u_R=1$).

The kinetic function obtained for $u_R=1$ is shown in Figure~\ref{Kinetiv_v--ur=1}, which shows the relation between right and left states of the nonclassical shock for $\tau$. Note that the kinetic functions are obtained only for $u_L<1.2$, where clear nonclassical shocks are generated. It is observed that for large $u_-$, all schemes give close results, while for small $u_-$,  the results largely depend on order of accuracy of the numerical scheme. Again, as the order of accuracy increases, the numerical solutions appear to converge, which supports Conjecture  \ref{conjecture}.
The kinetic function in this case is monotonic and strictly decreasing.

In order to check whether the kinetic function is single-valued, we now consider $u_R=1.5$. 
The kinetic function obtained in this case 
 for $\tau$ is shown in Figures \ref{Kinetiv_v--ur=1.5}. As it it observed, the results are similar to those of kinetic function in Figure \ref{Kinetiv_v--ur=1}. Similarly, the cases with $u_R=0$, $u_R=3$ and $u_R=5$ were also tested, and they led to the same kinetic function.
This shows that the kinetic function is single-valued here, at least for the cases considered above.

Next, we add the capillarity effects in (\ref{VD.smpl})
\be
\begin{split}
& \del_t \tau - \del_x u = 0, \\
& \del_t u + \del_x p(\tau) = \eps \, \del_x ( \del_x u ) - \alpha \eps^2 \del_{xxx} \tau.
\end{split}
\label{VD.000} 
\ee 
We repeat the experiments with the same data (i.e. $u_R=1$, $\tau_L=0.8$, $\tau_R=2$ and varying $u_L$), but now including capillarity effect with $\alpha=1$. The capillarity coefficient is $\alpha \eps^2$ as in the previous cases.
The kinetic functions for various schemes are shown in Figure~\ref{disp}, (which shows the relation between right and left states of the nonclassical shock for $\tau$), and they are compared with the case of no capillarity effect. As it is observed, for a given $\tau_-$, the capillarity effect 
leads to smaller values of $\tau_+$.

\subsection{A piecewise linear pressure function}
In order to highlight the effect of a pressure function with two inflection points, here we use a 
piecewise linear flux function which was already studied in \cite{BCCL}, given by 
 \be
 p(\tau)=
 \begin{cases}
  -7 \tau + 10,         &          \tau \leq 1,    \\
  4 \tau - 1,           &        1< \tau \leq 2,     \\
-\frac{5}{2} \tau +12,      &        2< \tau \leq 4,     \\
-\frac{1}{5} \tau +\frac{4}{15},      &        4< \tau,     \\
 \end{cases}
 \label{tauLin} \ee
and shown in Figure  \ref{fluxLin}. 
In order to explore various regimes occurring by this pressure function, we set $\tau_L=0.9$ and $\tau_R=4$. 
The two inflection points are again between the selected  $\tau_L$ and $\tau_R$.\\
\

Note that, since the pressure function is piecewise linear, $\tau$ in the system \ref{VD.000} 
does not depend on the specific choice of $u_L$ and $u_R$, but depends only on their difference $u_L-u_R$. 
This is easily checked by a change of variable $u \mapsto u+c$, where $c$ is any constant. 
Hence, in the numerical experiments without loss of generality, we can fix $u_R=1$, and change $u_L$ only. 

Numerical results are obtained using a computational grid of $2000$ nodes 
and $h=0.0005$ with $\eps=0.001$ and $\alpha=0$ (no capillarity effects).
 
\
 
Three regimes are identified, as explained now: 
\begin{itemize}

	\item \emph{Regime A: $1.4 \leq u_L$ and a stationary shock.}\\
In this regime, a stationary shock wave is generated in the center. A typical solution
is plotted in Figure~\ref{RegimeA} for $u_L=1.5$. The kinetic function in this regime based on
the fourth order scheme is shown in Figure~\ref{RegimeAKin}, and it appears to be a linear function.	
\\

	\item \emph{Regime B: $ -1.2 <u_L \leq 1.4$ and a non-stationary shock.}\\ 
In this regime, a (left-going) nonclassical shock wave is generated. A typical picture of this regime
 is shown in Figure~\ref{RegimeB} for $u_L=1$. As $u_L$ decreases, the shock speed is increased, and 
 the left-hand and right-hand values of the nonclassical shock reach some limiting values. 
The kinetic function in this regime is \emph{non-monotone} and \emph{not single-valued.}
 \\	

	\item \emph{Regime C: $ u_L \leq -1.2$ and a non-stationary shock with fixed left- and right-hand values.}\\
In this regime, a (left-going) nonclassical shock wave is generated, but the values $\tau_-$ and $\tau_+$ are fixed
to the limiting values of Regime B awhile the specific volume $\tau$ is increased on the right-hand side of
the nonclassical shock and generates a right-moving wave. A typical solution for this regime
is plotted in Figure~\ref{RegimeC} for $u_L=-2$.
No kinetic function arises in this regime since $\tau_-$ and $\tau_+$ do not change.

\

Finally, kinetic functions for both Regimes A and B are given in Figure~\ref{Both}, using
 fourth- and eighth-order schemes. Note that the kinetic functions are almost identical in the Regime A (stationary shock) 
for both schemes, but they are largely distinct in the Regime B.

\begin{figure} 
\begin{center}
\begin{tabular}{cc}
\begin{minipage}{6.0cm}
\begin{center}
\includegraphics*[height=7cm, width=7cm]{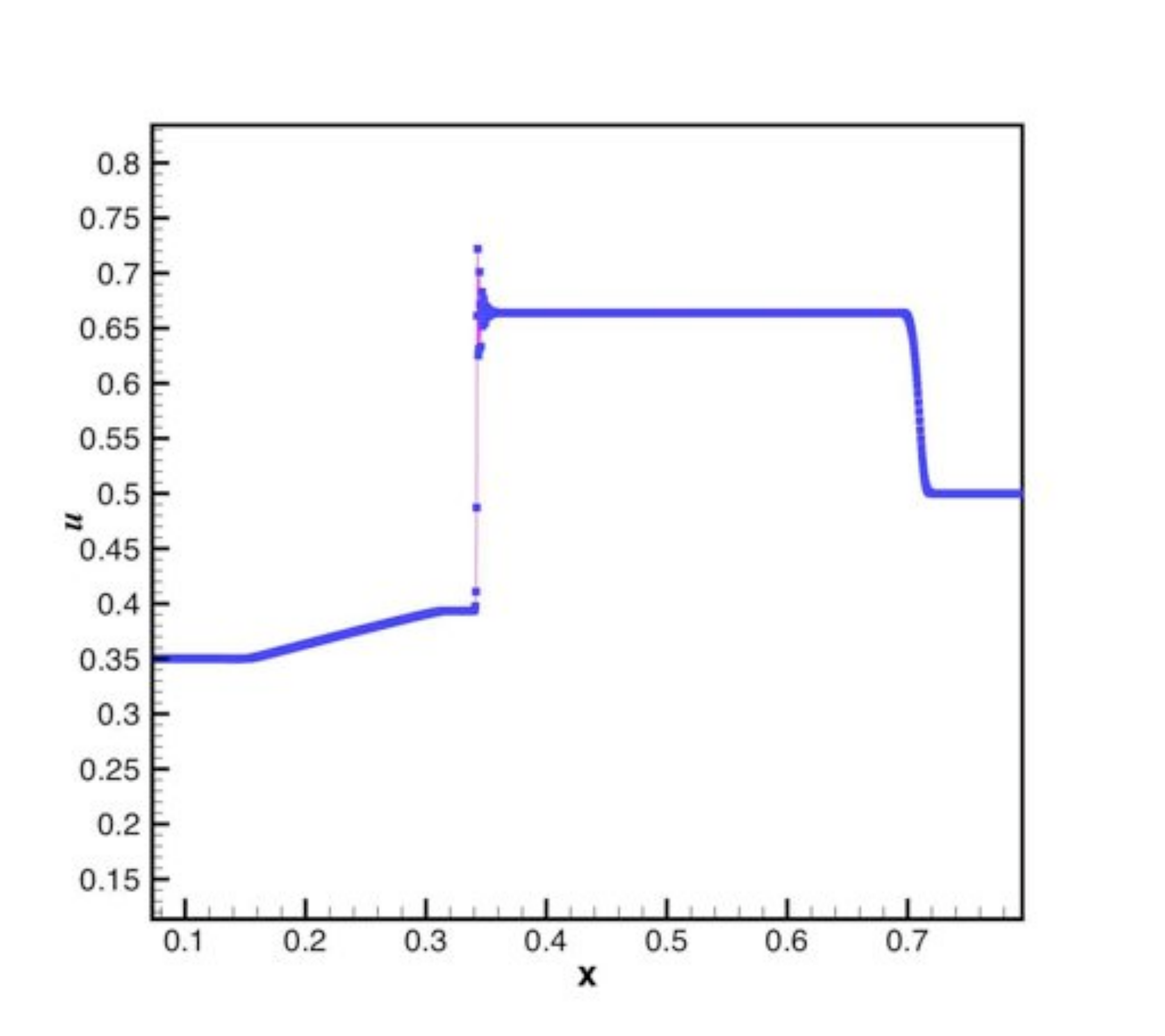}
\end{center}
\end{minipage}
&
\begin{minipage}{6.0cm}
\begin{center}
\includegraphics*[height=7cm, width=7cm]{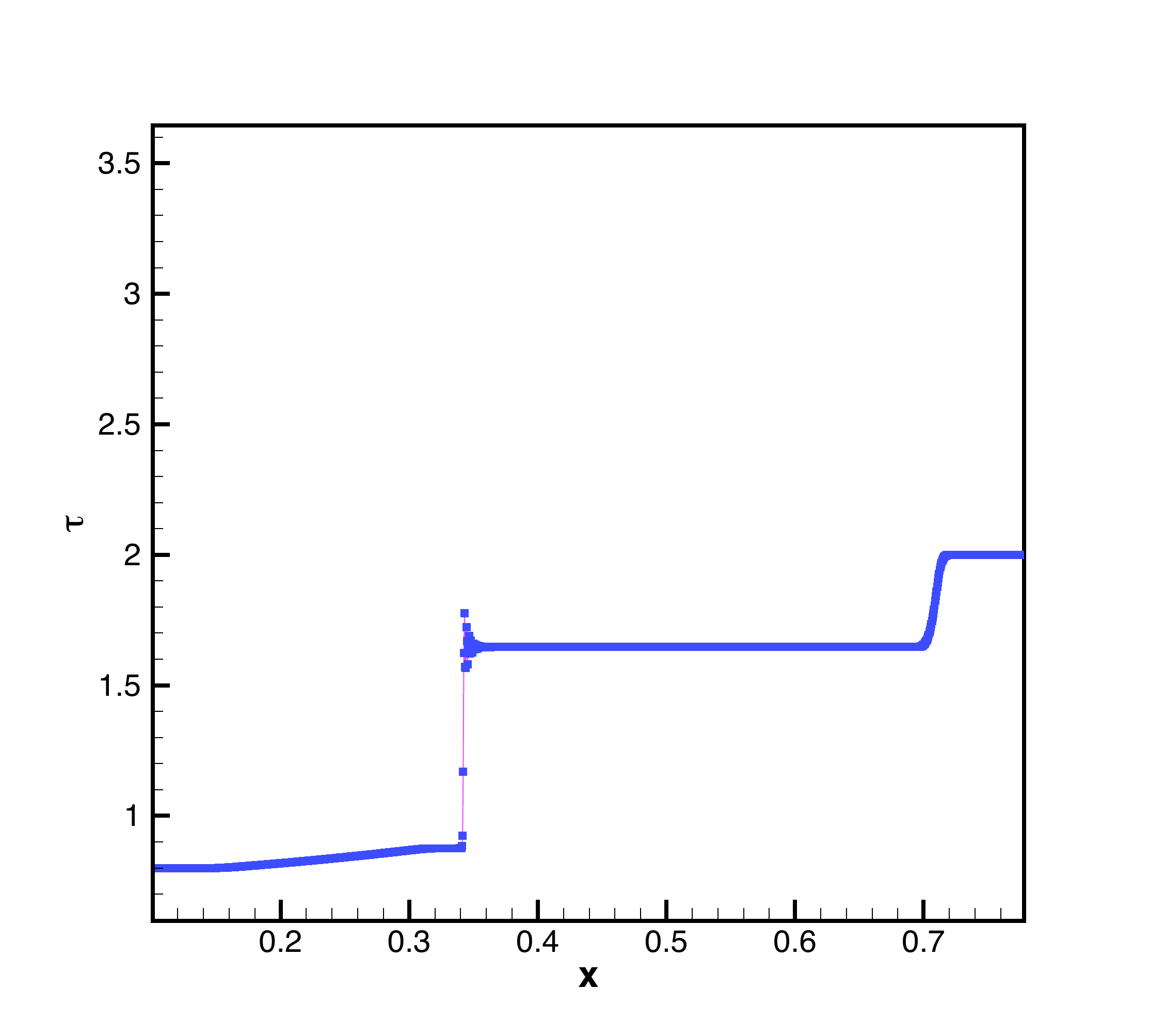}
\end{center}
\end{minipage}
\end{tabular}
\caption{Typical wave structure for van der Waal fluids; left $u$ and right $\tau$, for $u_R=0.5$.  }
\label{vanderwaals0.5}
\end{center}
\end{figure}
\begin{figure} 
\begin{center}
\begin{tabular}{cc}
\begin{minipage}{6.0cm}
\begin{center}
\includegraphics*[height=7cm, width=7cm]{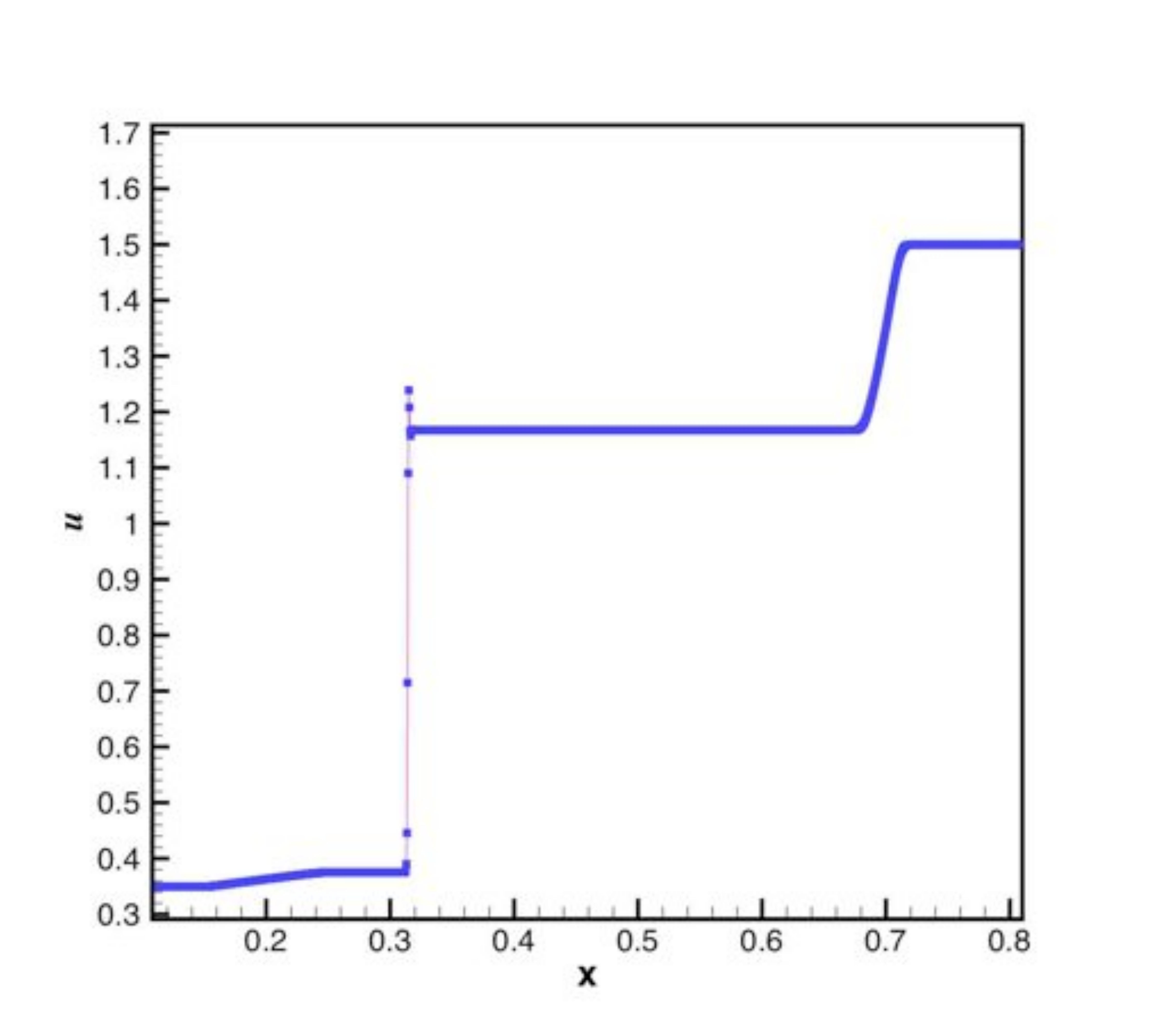}
\end{center}
\end{minipage}
&
\begin{minipage}{6.0cm}
\begin{center}
\includegraphics*[height=7cm, width=7cm]{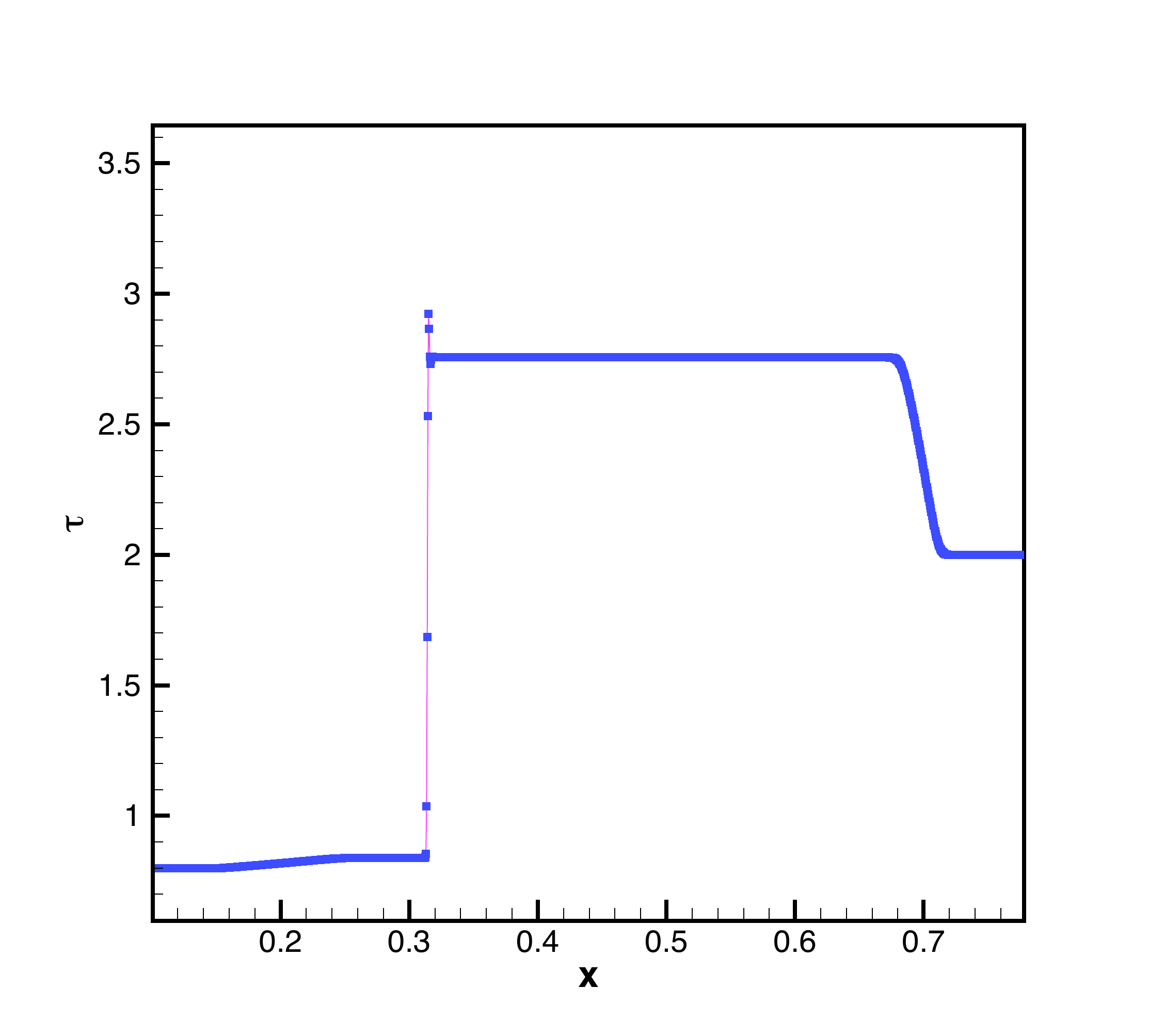}
\end{center}
\end{minipage}
\end{tabular}
\caption{Typical wave structure for van der Waal fluids; left $u$ and right $\tau$, for $u_R=1.5$.  }
\label{vanderwaals1.5}
\end{center}
\end{figure}
\begin{figure} 
\begin{center}
\begin{tabular}{cc}
\begin{minipage}{6.0cm}
\begin{center}
\includegraphics*[height=7cm, width=7cm]{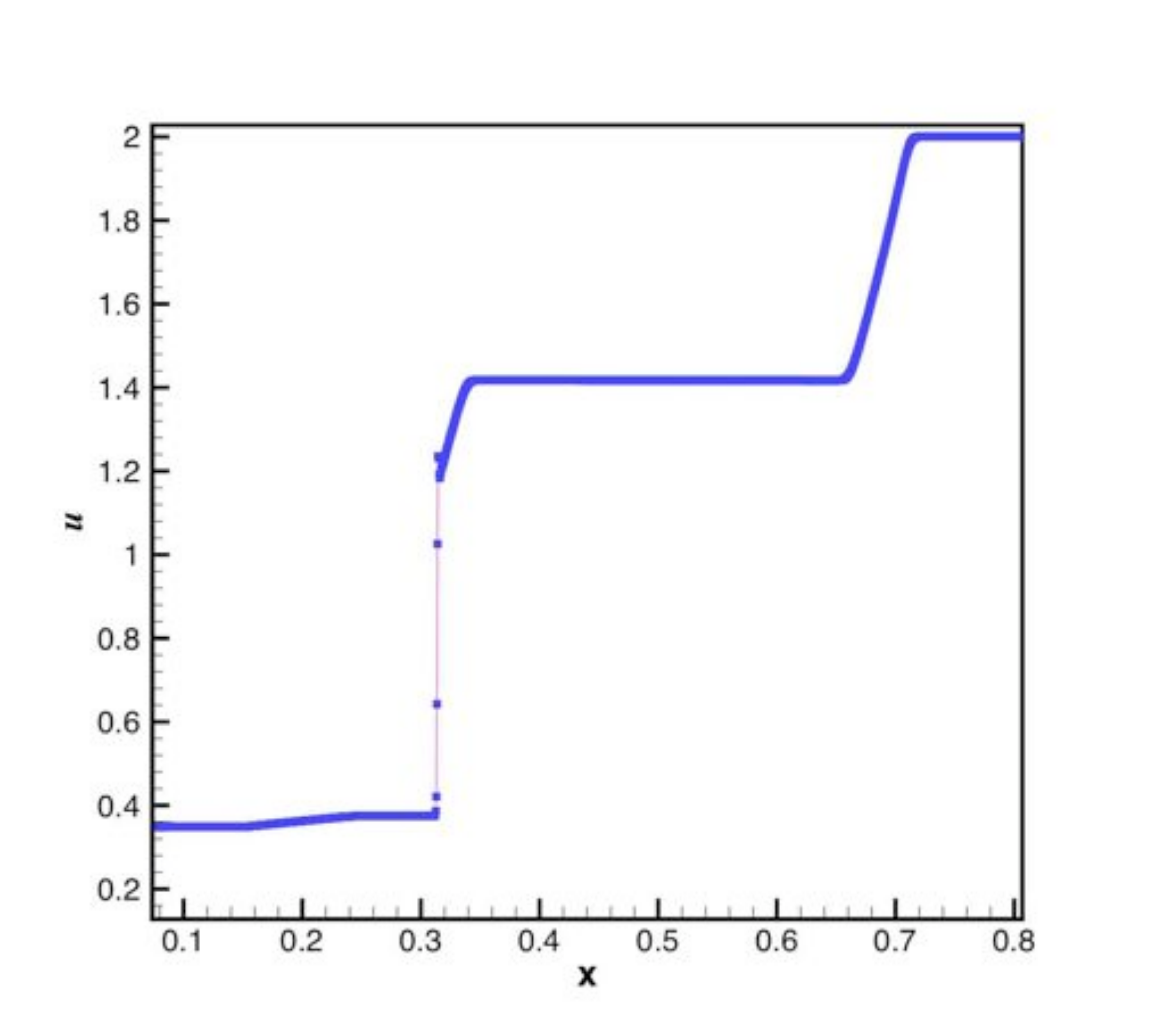}
\end{center}
\end{minipage}
&
\begin{minipage}{6.0cm}
\begin{center}
\includegraphics*[height=7cm, width=7cm]{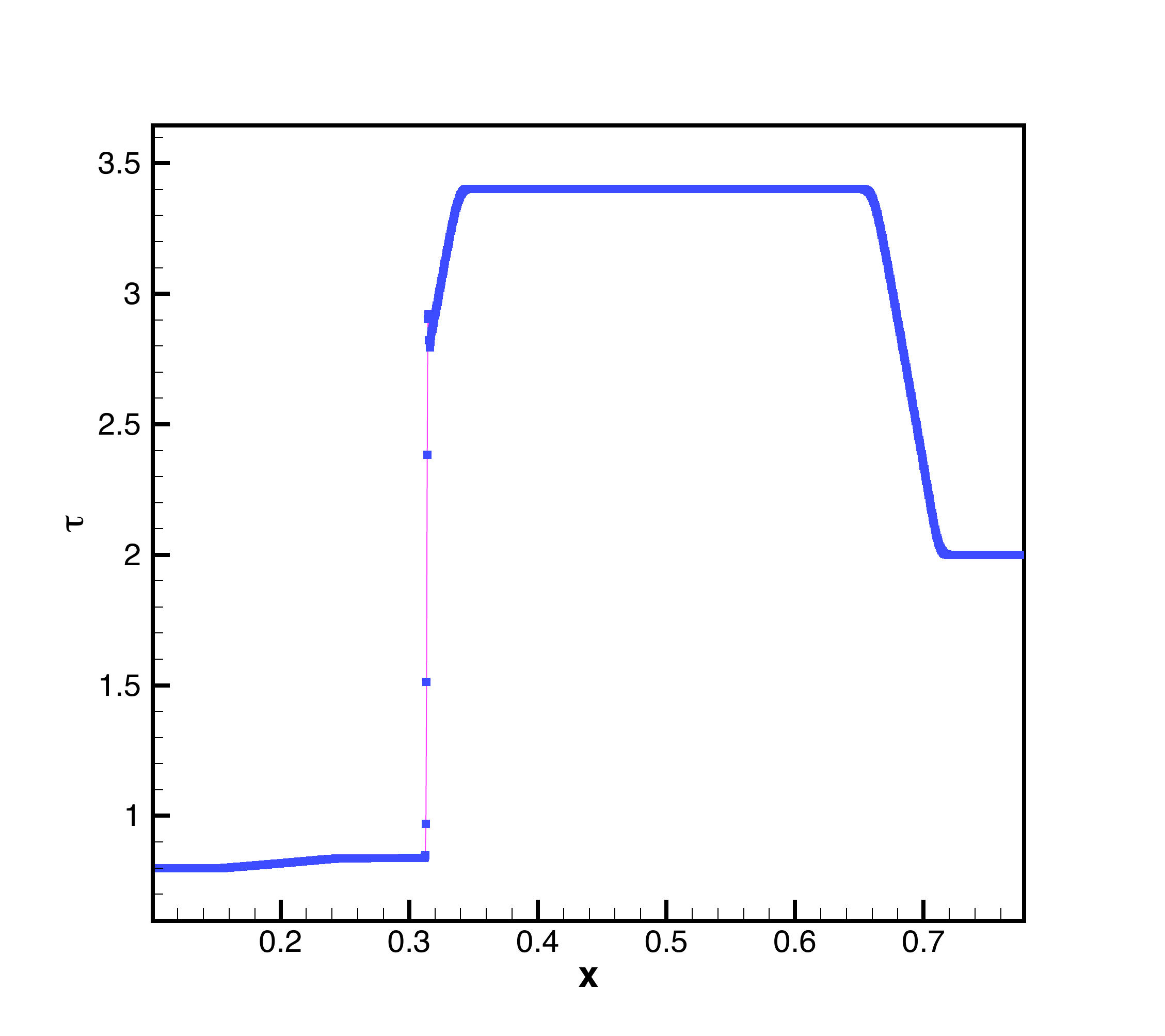}
\end{center}
\end{minipage}
\end{tabular}
\caption{Typical wave structure for van der Waal fluids; left $u$ and right $\tau$, for $u_R=2$.  }
\label{vanderwaals2}
\end{center}
\end{figure}
\begin{figure}
\begin{center}
  \includegraphics[width=8cm,height=8cm]{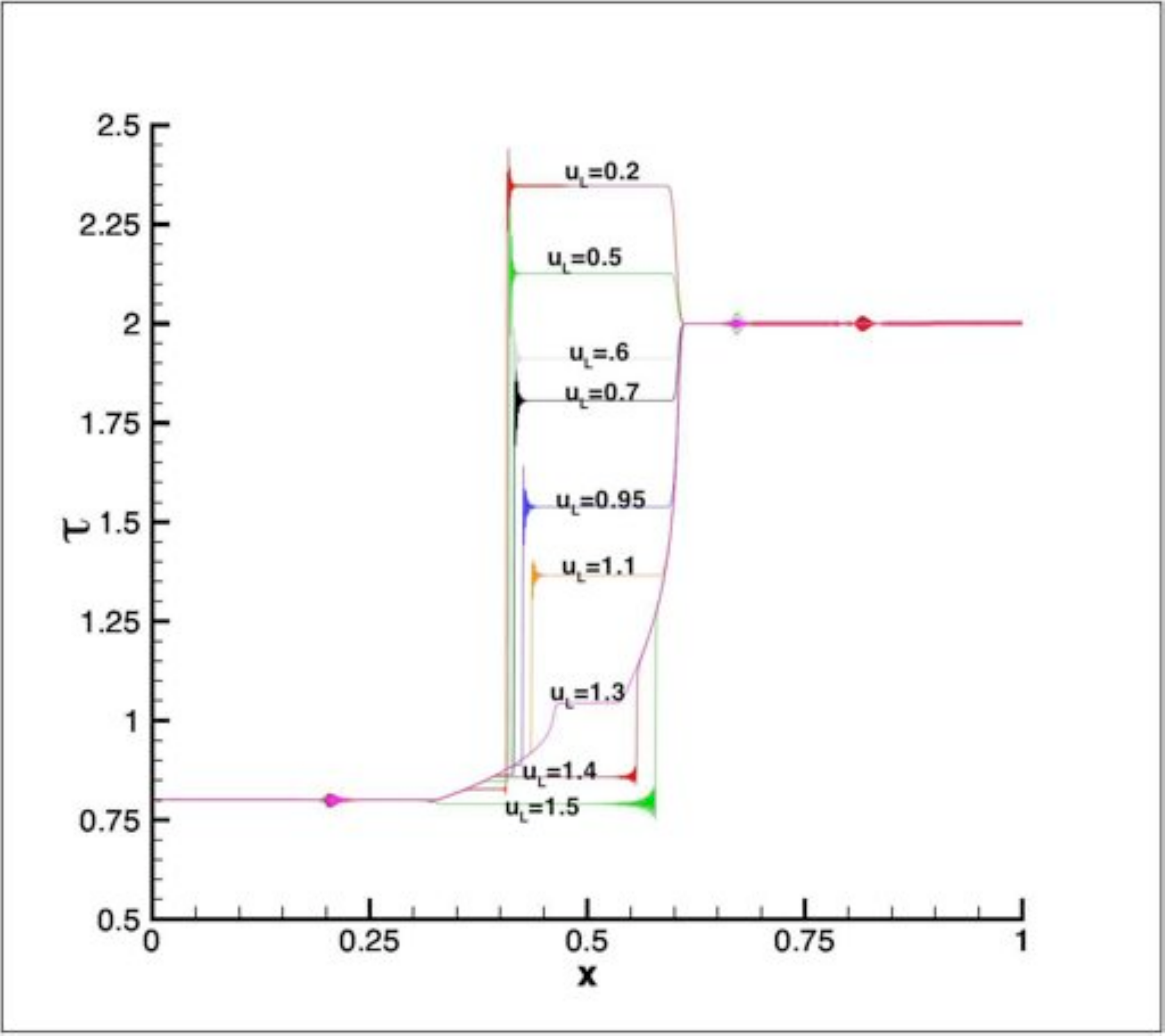}
  \caption{Wave structure ($\tau$), corresponding to the case $\tau_L=0.8$, $\tau_R=2$, $u_R=1$ and variable $u_L$.}
 \label{tau}
\end{center}
\end{figure} 
\begin{figure}
\begin{center}
  \includegraphics[width=8cm,height=8cm]{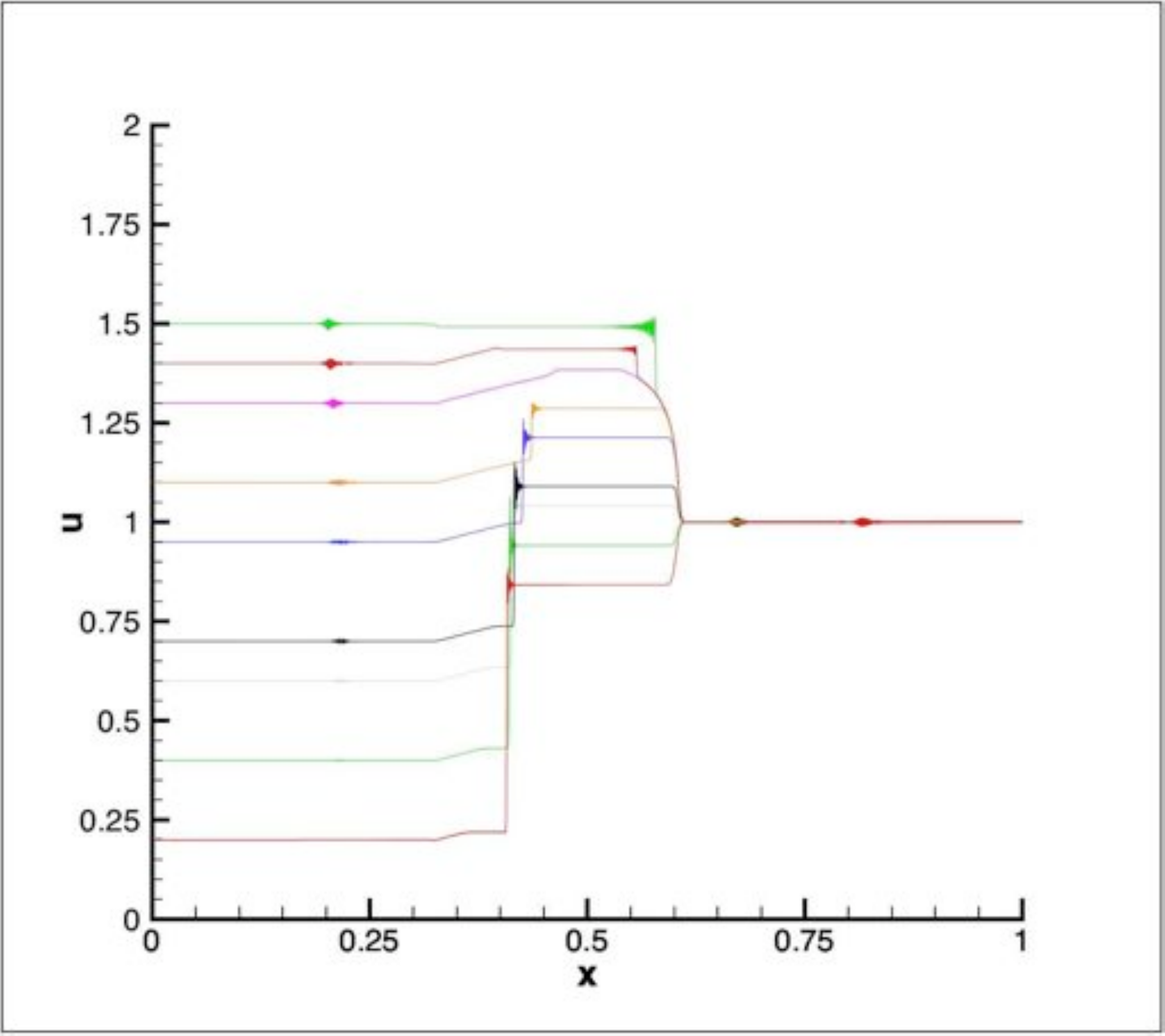}
  \caption{Wave structure ($u$), corresponding to the case $\tau_L=0.8$, $\tau_R=2$, $u_R=1$ and variable $u_L$.}
 \label{u}
\end{center}
\end{figure} %
\begin{figure}
\begin{center}
  \includegraphics[width=8cm,height=8cm]{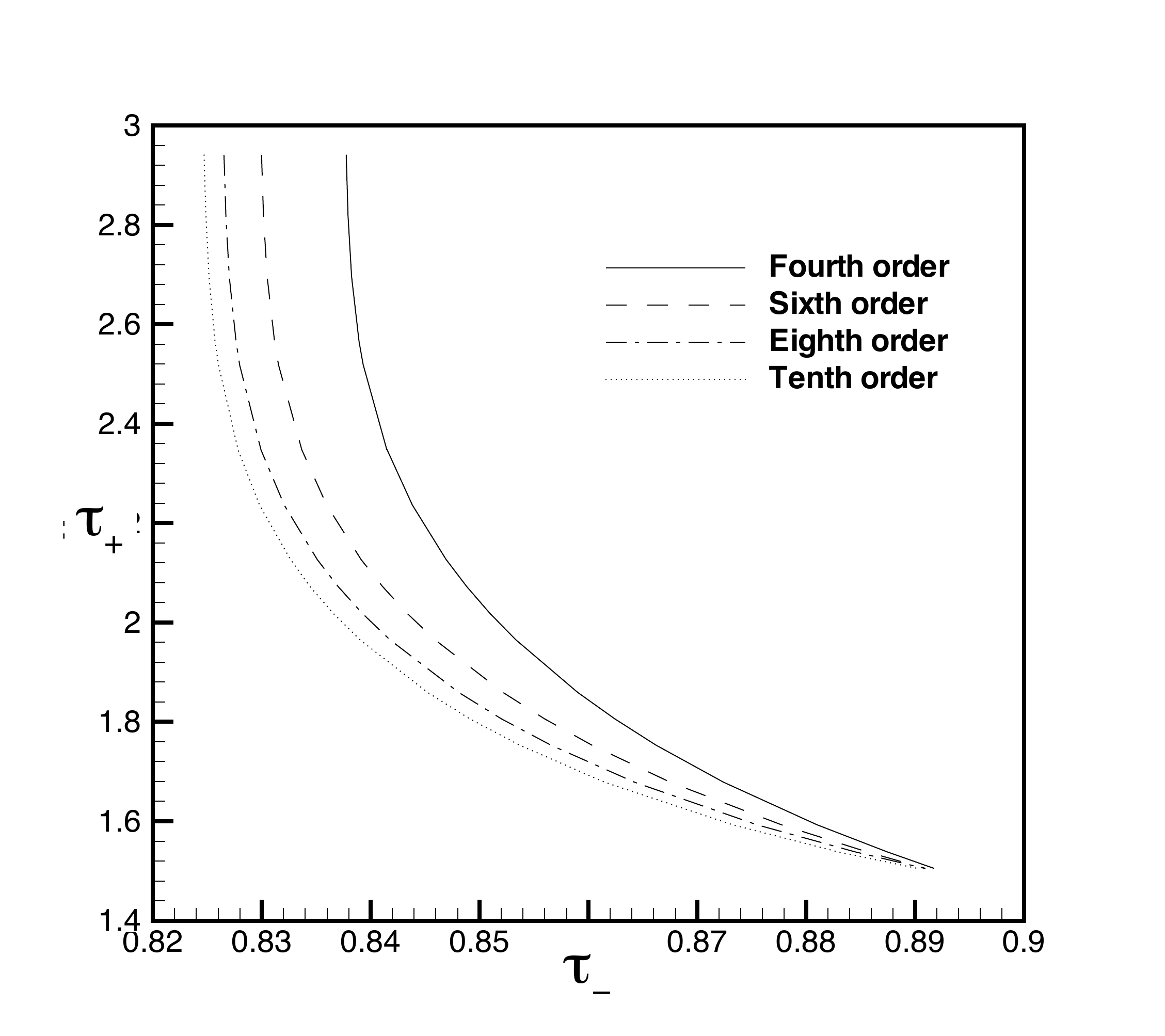}
  \caption{The kinetic function for $\tau$ obtained for $u_R=1$.}
 \label{Kinetiv_v--ur=1}
\end{center}
\end{figure} 
\begin{figure}
\begin{center}
  \includegraphics[width=8cm,height=8cm]{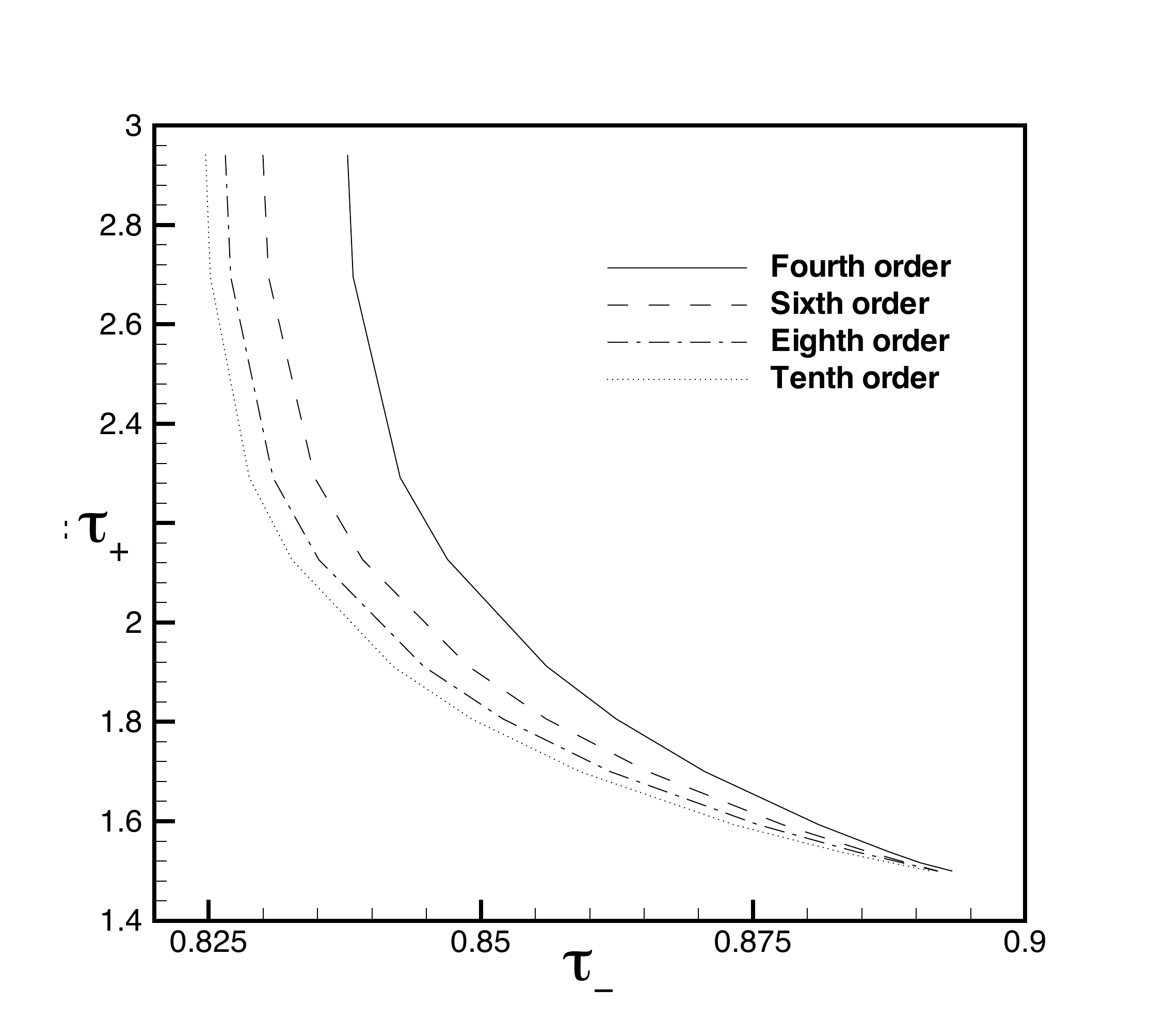}
  \caption{The kinetic function for $\tau$ obtained for $u_R=1.5$.}
 \label{Kinetiv_v--ur=1.5}
\end{center}
\end{figure}
\begin{figure}
\begin{center}
  \includegraphics[width=8cm,height=8cm]{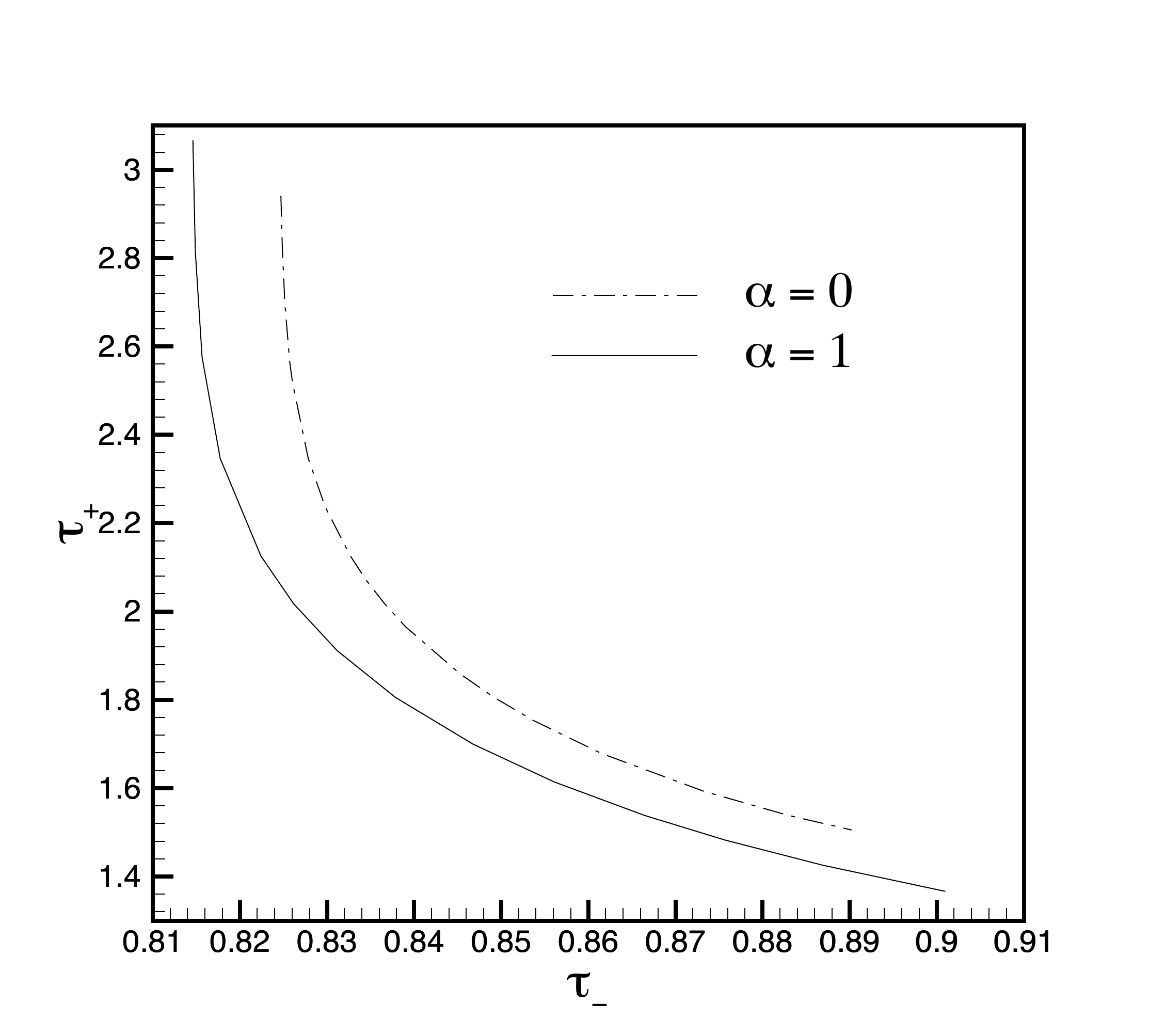}
  \caption{The kinetic function for $\tau$ obtained for $u_R=1$, and capillarity effect with $\alpha=1$, using the tenth order scheme. The dash curve shows the case without the capillarity effect.}
 \label{disp}
\end{center}
\end{figure} 
%
%
\begin{figure}
\begin{center}
  \includegraphics[width=7cm,height=7cm]{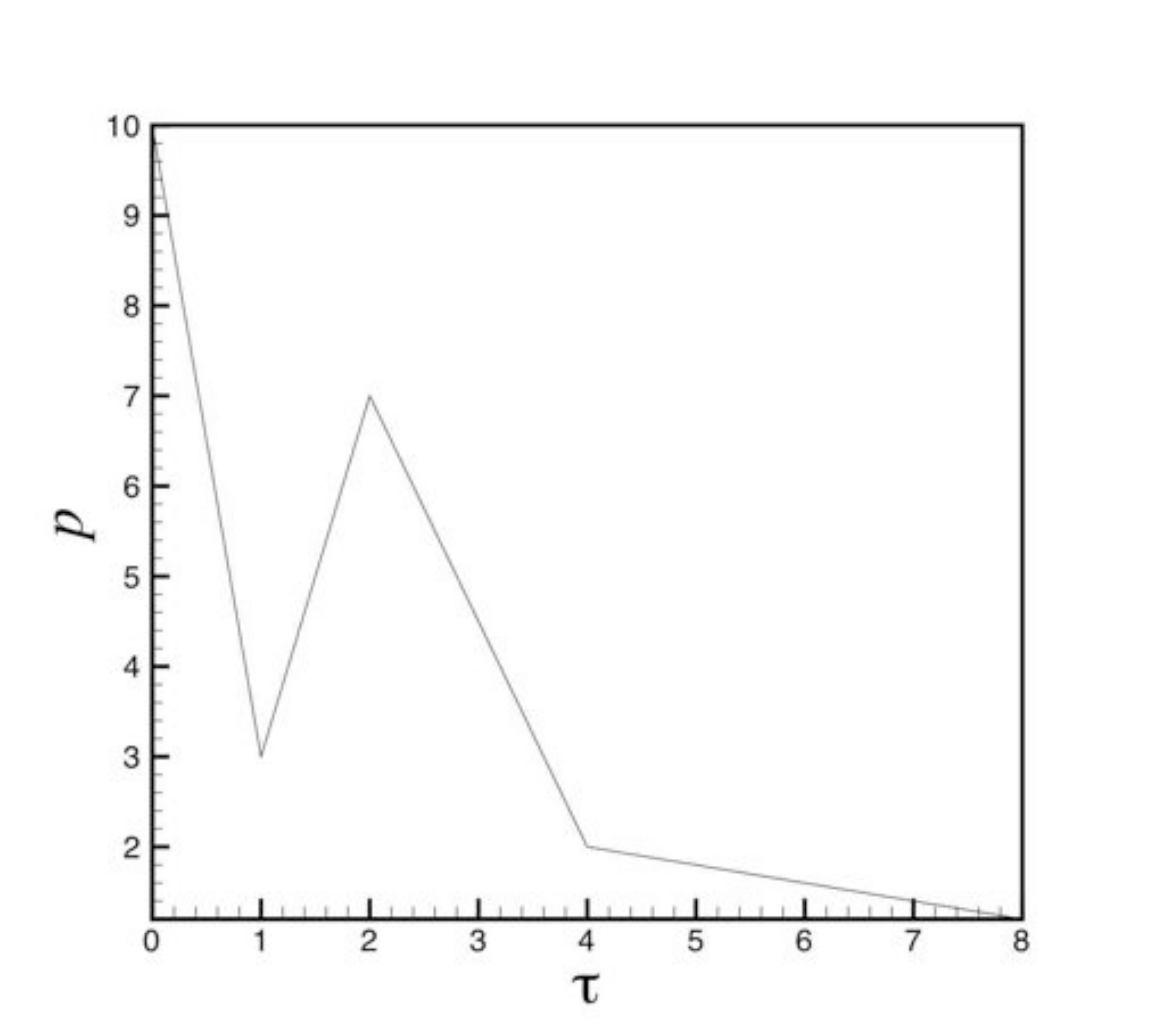}
  \caption{A piecewise linear flux function resembling the van der Waals pressure law.}
 \label{fluxLin}
\end{center}
\end{figure} %
%
\begin{figure} 
\begin{center}
\begin{tabular}{cc}
\begin{minipage}{7.0cm}
\begin{center}
\includegraphics*[height=7.0cm, width=7.0cm]{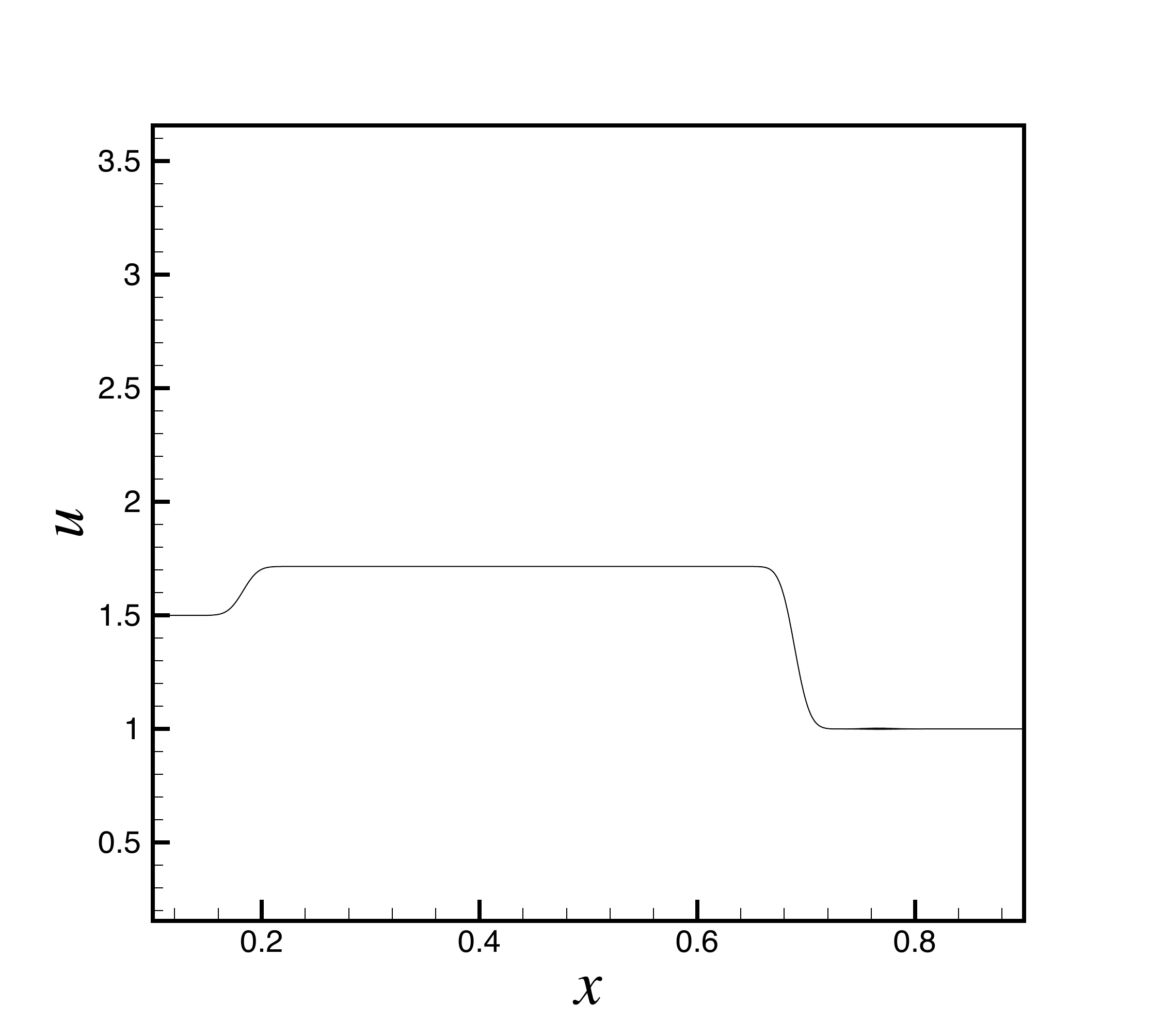}
\end{center}
\end{minipage}
&
\begin{minipage}{7.0cm}
\begin{center}
\includegraphics*[height=7.0cm, width=7.0cm]{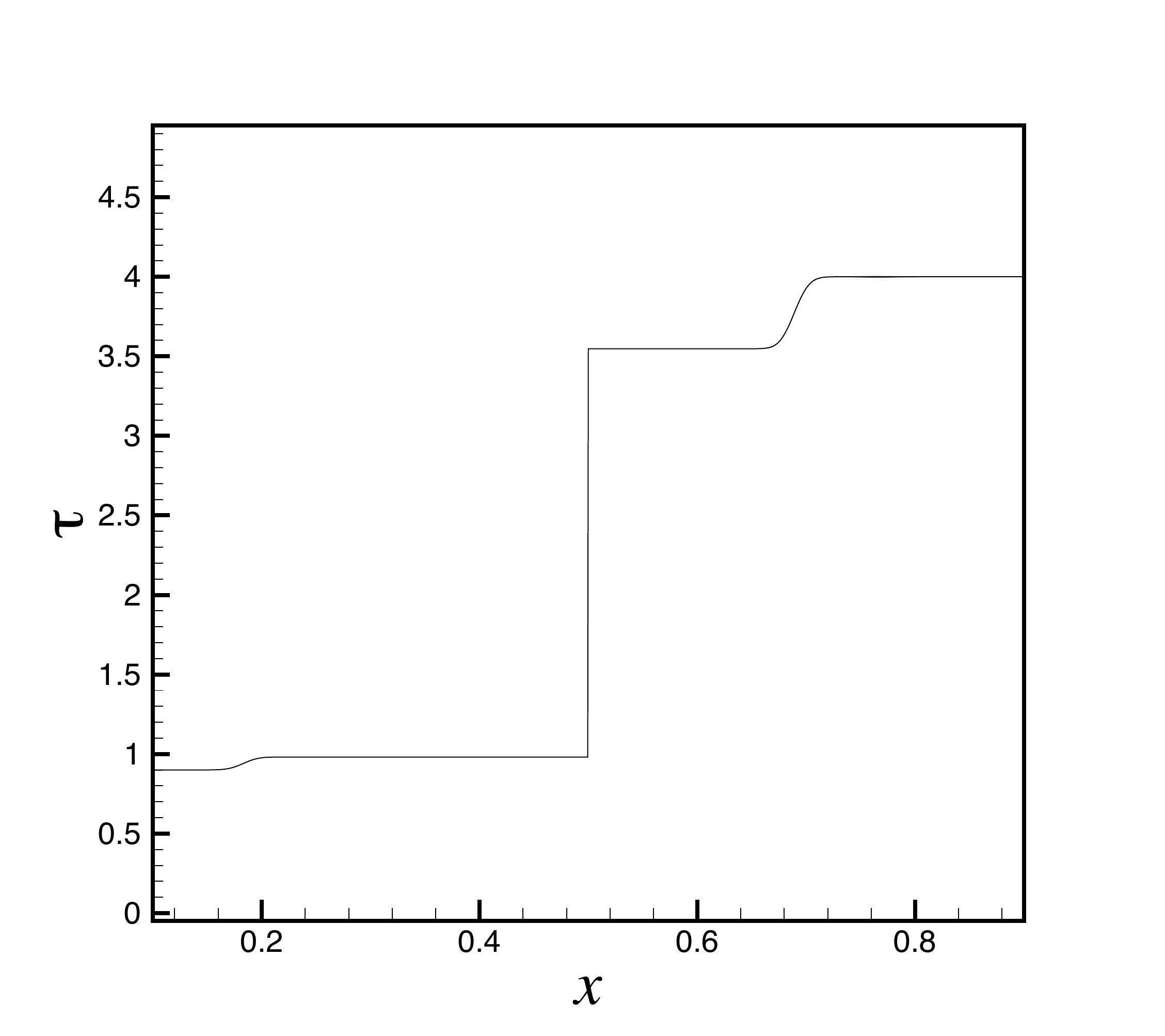}
\end{center}
\end{minipage}
\end{tabular}
\caption{A typical wave structure in Regime A; $u$ (left) and $\tau$ (right) at time $t=0.12$.}
\label{RegimeA}
\end{center}
\end{figure}
\begin{figure}
\begin{center}
  \includegraphics[width=7cm,height=7cm]{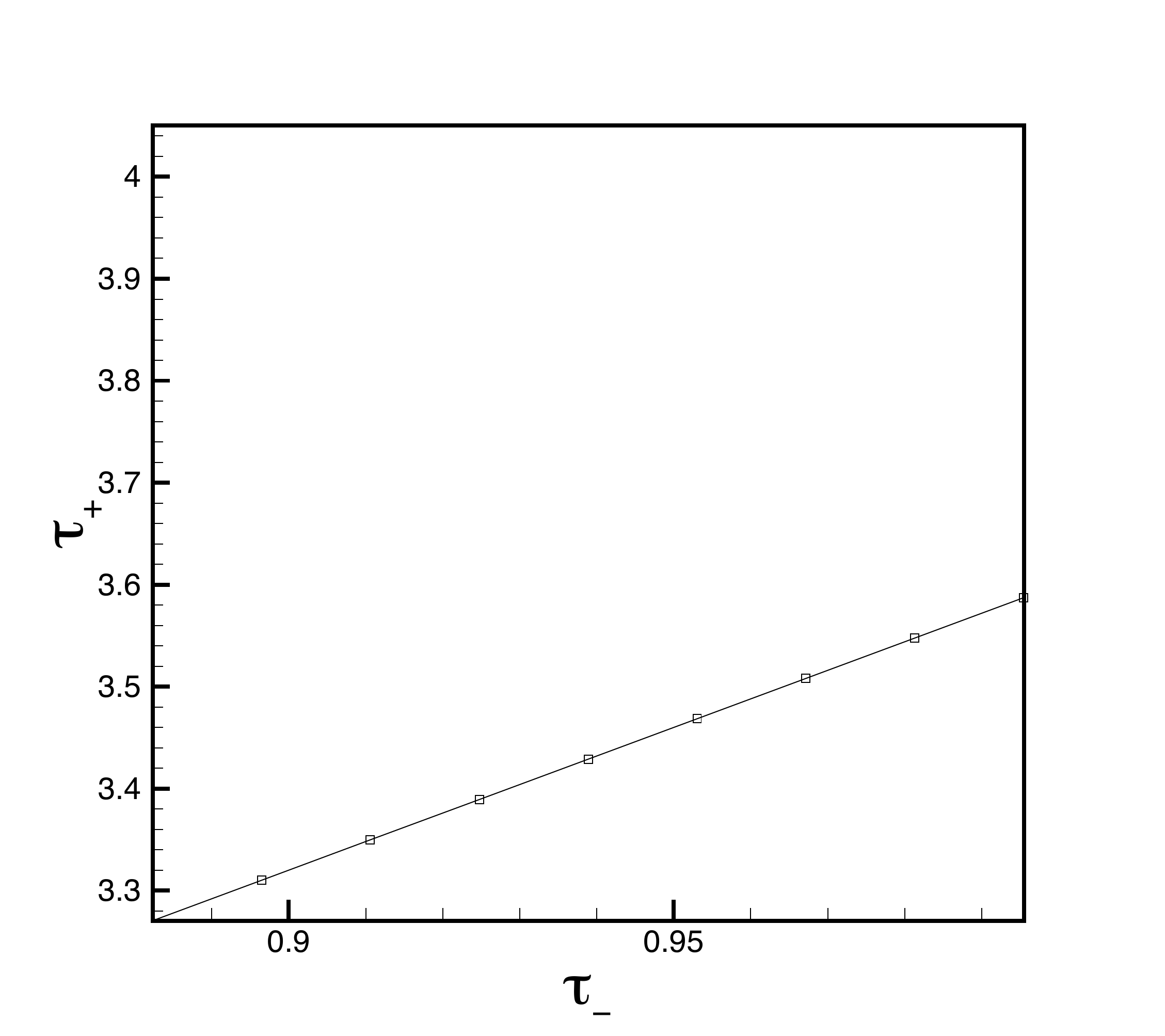}
  \caption{Kinetic function for Regime A.}
 \label{RegimeAKin}
\end{center}
\end{figure} 
%
%
\begin{figure} 
\begin{center}
\begin{tabular}{cc}
\begin{minipage}{7.0cm}
\begin{center}
\includegraphics*[height=7.0cm, width=7.0cm]{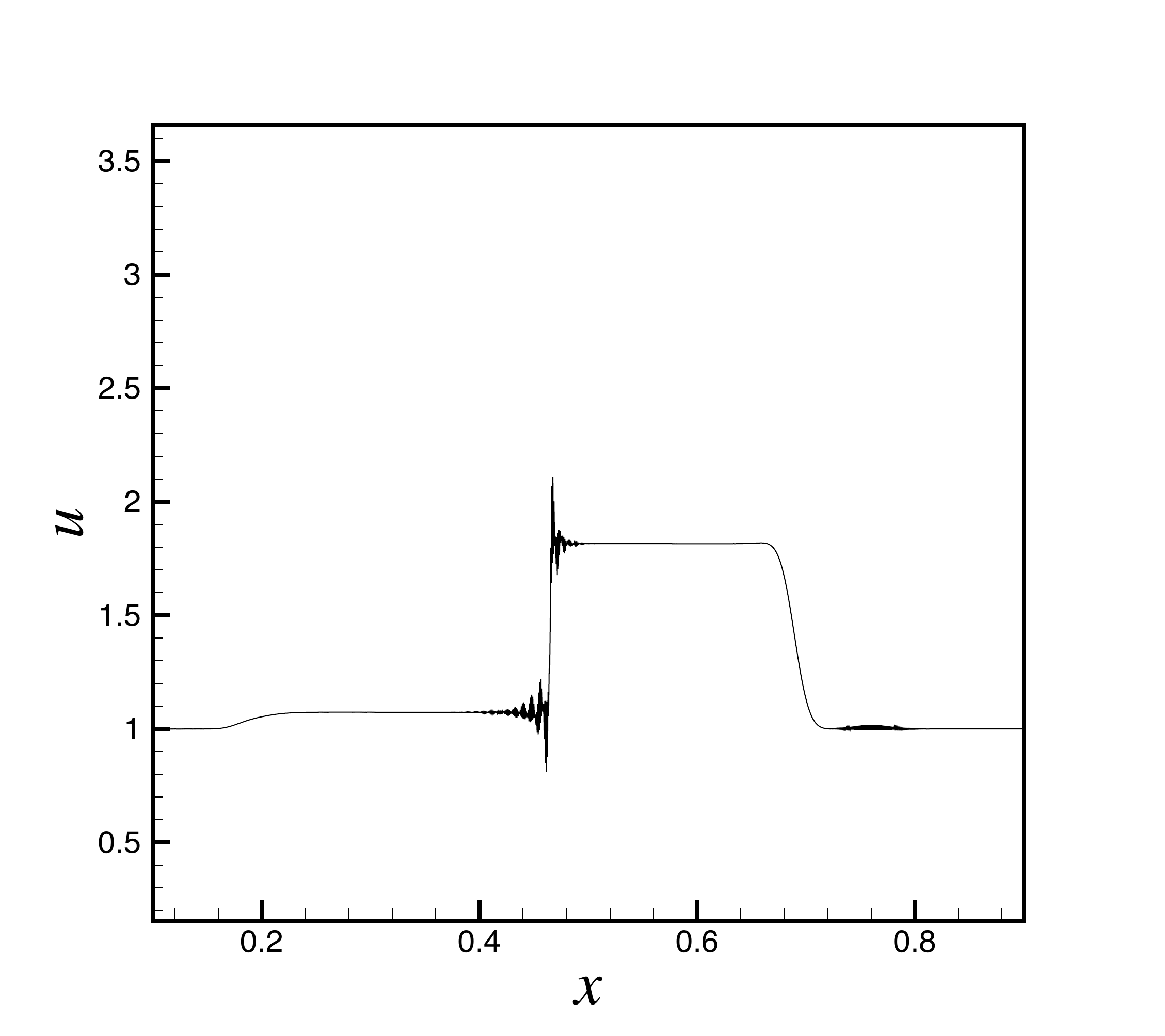}
\end{center}
\end{minipage}
&
\begin{minipage}{7.0cm}
\begin{center}
\includegraphics*[height=7.0cm, width=7.0cm]{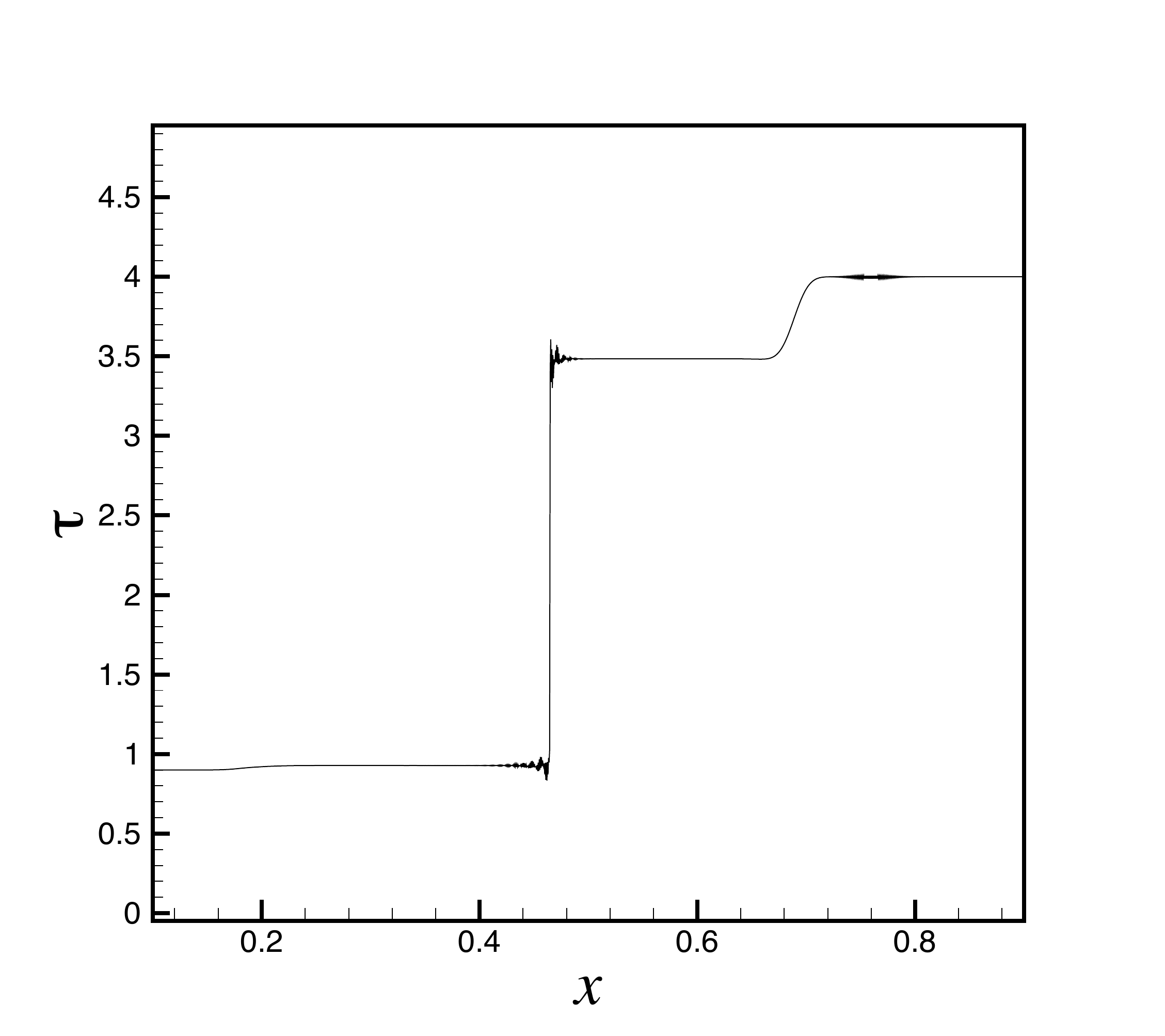}
\end{center}
\end{minipage}
\end{tabular}
\caption{A typical wave structure in Regime B; $u$ (left) and $\tau$ (right) at time $t=0.12$.}
\label{RegimeB}
\end{center}
\end{figure}
\begin{figure}
\begin{center}
  \includegraphics[width=7cm,height=7cm]{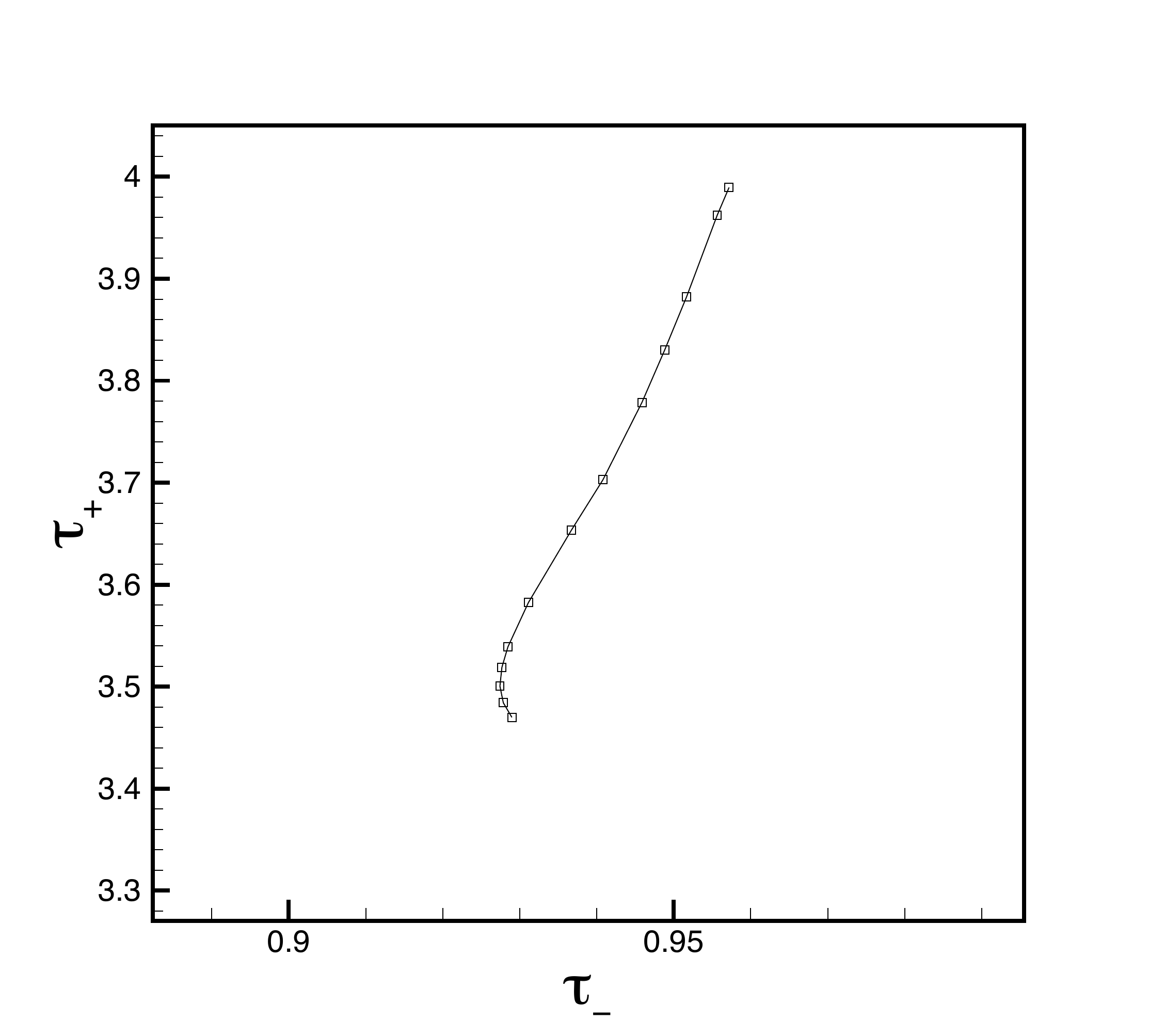}
  \caption{Kinetic function for Regime B.}
 \label{RegimeBKin}
\end{center}
\end{figure} 
%
%
\begin{figure} 
\begin{center}
\begin{tabular}{cc}
\begin{minipage}{7.0cm}
\begin{center}
\includegraphics*[height=7.0cm, width=7.0cm]{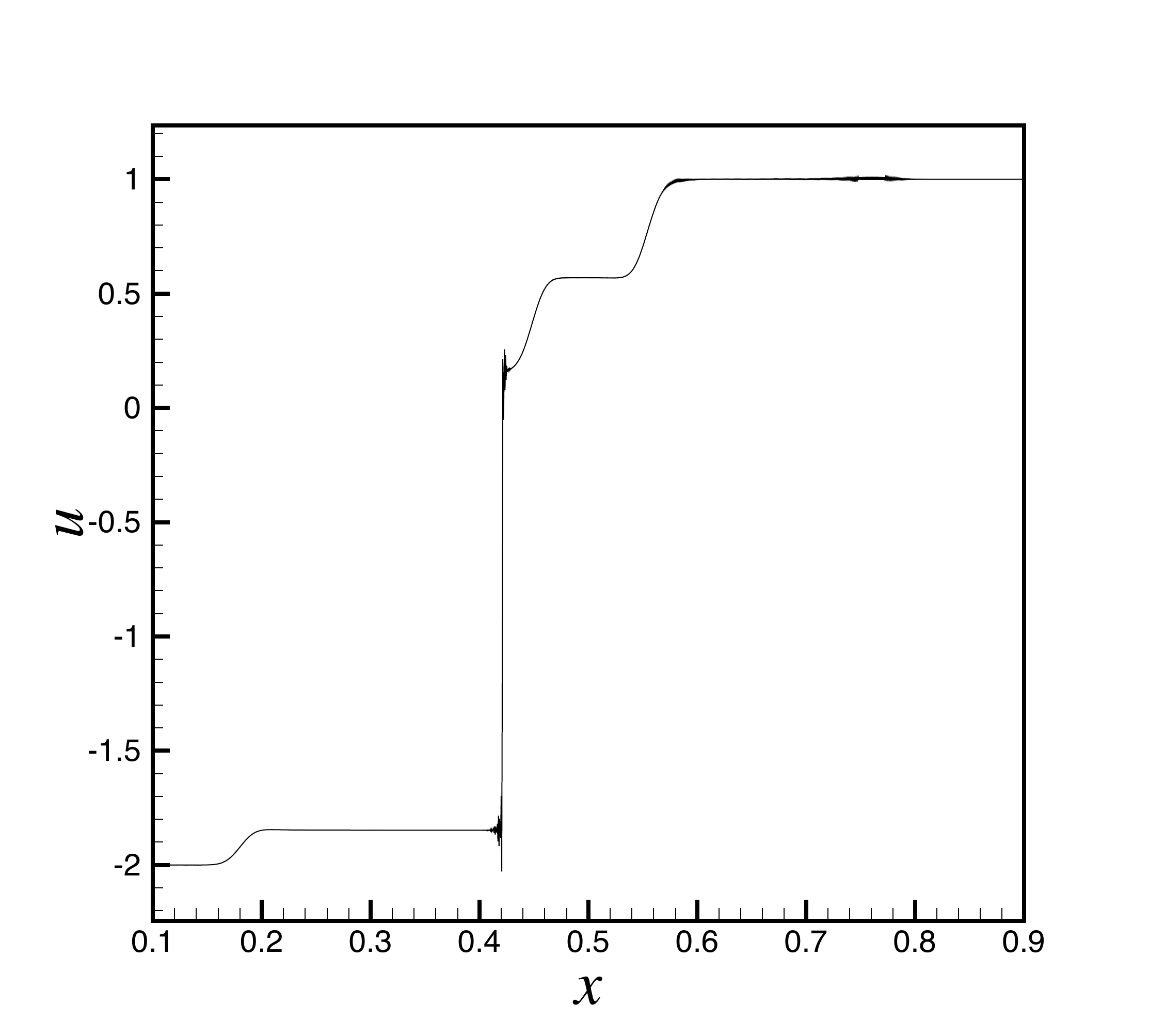}
\end{center}
\end{minipage}
&
\begin{minipage}{7.0cm}
\begin{center}
\includegraphics*[height=7.0cm, width=7.0cm]{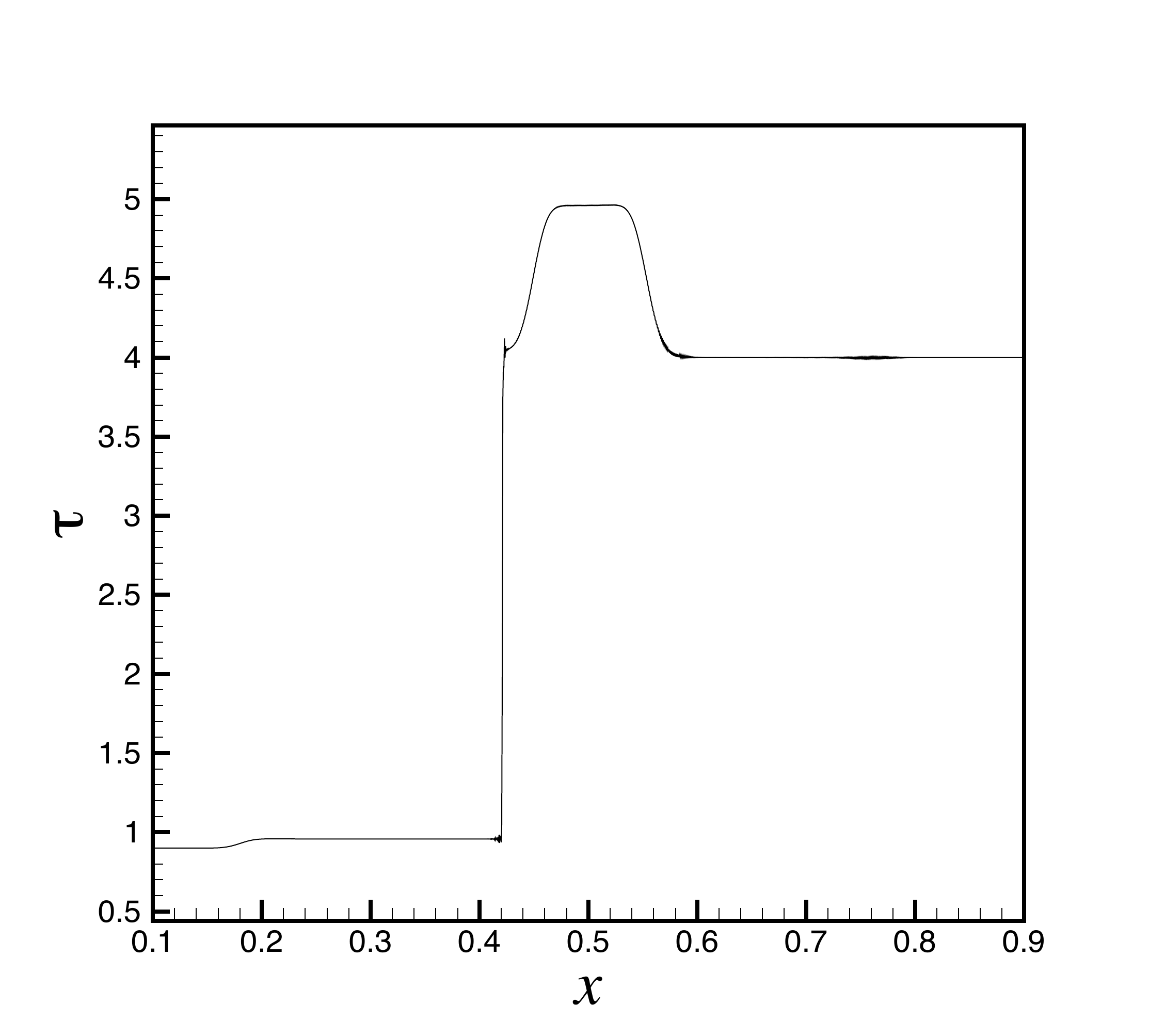}
\end{center}
\end{minipage}
\end{tabular}
\caption{A typical wave structure in Regime C; $u$ (left) and $\tau$ (right) at time $t=0.12$.}
\label{RegimeC}
\end{center}
\end{figure}
\begin{figure}
\begin{center}
  \includegraphics[width=7cm,height=7cm]{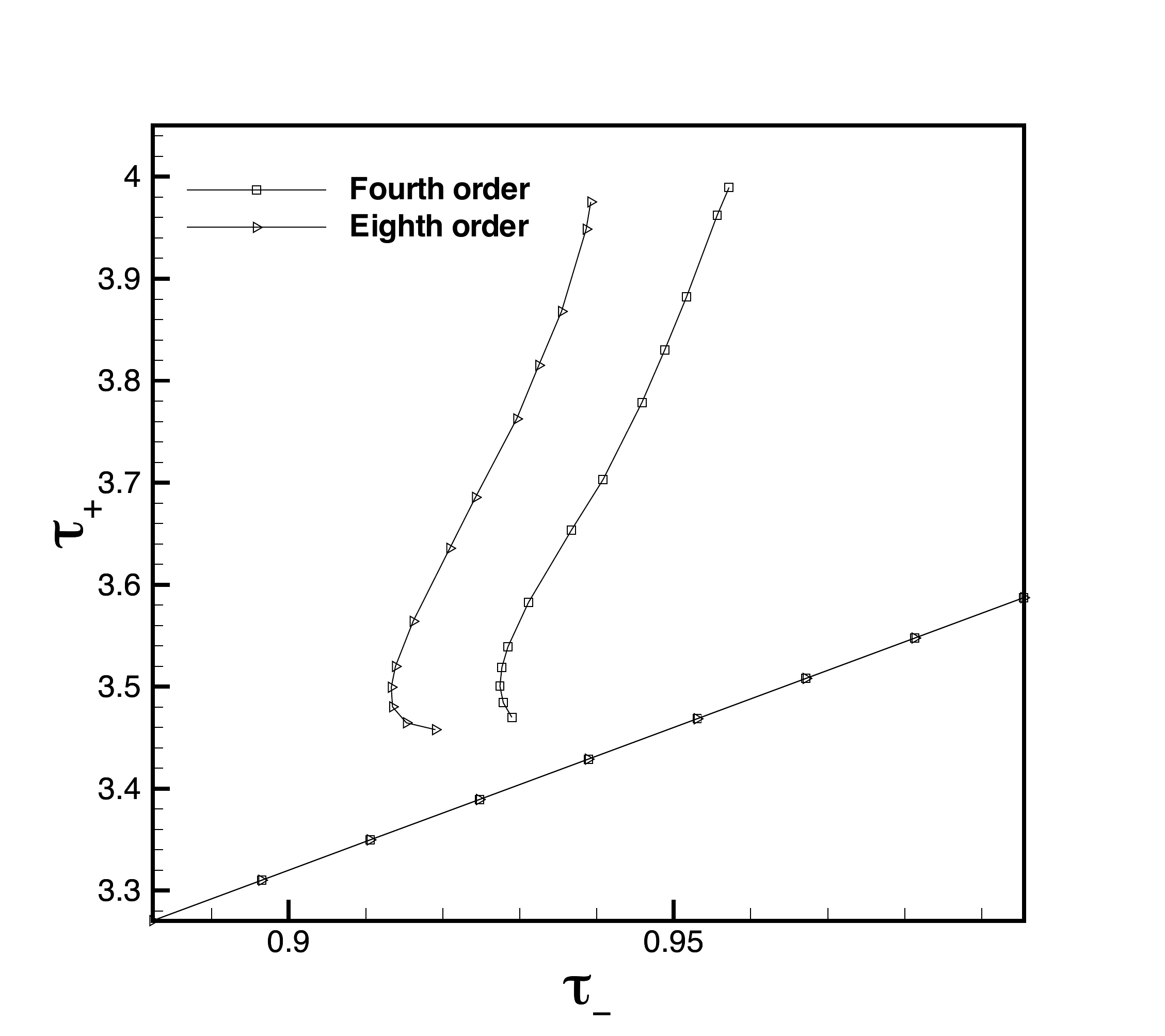}
  \caption{Kinetic function for both Regims A and B, using the fourth and the eighth order schemes.}
 \label{Both}
\end{center}
\end{figure} 
\end{itemize}

\

\


\section{Concluding remarks}
\label{CO-0}
 
{
The specific contributions made in this paper are three-fold.

\begin{itemize}

\item First, we have investigated the role of the equivalent equation associated with a scheme. 
We stated a conjecture and provided numerical evidences demonstrating its validity.
We have shown that the kinetic function associated with a finite difference scheme approaches the
(exact) kinetic function derived from a given (viscosity, capillarity) regularization. 
Interestingly, the accuracy improves as its equivalent equation coincides with the diffusive-dispersive model 
at a higher and higher order of approximation. 
We observed that the spatial accuracy plays a more crucial role than the temporal accuracy. 
In both the continuous model and the discrete scheme,
small scale features are critical to the selection of shocks. Hence, the balance between diffusive and dispersive
features determines which shocks are selected.
These small scale features can not be quite the same at the continuous and at the discrete levels,
since a \emph{continuous} dynamical system of ordinary differential equations 
can not be exactly represented by a \emph{discrete dynamical system} of finite difference equations.  
The effects of the regularization coefficients on the kinetic functions were also investigated.

\item Second, we considered fourth-order models. We demonstrated that a kinetic function can be associated
with the thin liquid film model. We also investigated a generalized Camassa-Holm model, and discovered 
nonclassical shocks for which we could determine a kinetic relation. In both models, the kinetic function 
was found to be monotone decreasing, as required in the general theory \cite{LeFloch2}. 

\item Third, we investigated to what extent a kinetic function can 
be associated with van der Waals fluids, whose flux-function admits two inflection points. 
We established that the Riemann problem admits several solutions whose discontinuities 
have viscous-capillary profiles, and we exhibited non-monotone kinetic functions.  

\end{itemize} 

In summary, the new features of the kinetic functions observed here 
suggests directions for further challenging theoretical work 
on the well-posedness theory for the Cauchy problem associated with nonlinear hyperbolic problems. 

}


\section*{Acknowledgments}
 
The first author (PLF) 
was partially supported by the A.N.R. (Agence Nationale de la Recherche) through the grant 06-2-134423
entitled {\sl ``Mathematical Methods in General Relativity''} (MATH-GR), and by the Centre National de la Recherche Scientifique (CNRS). The second author (MM)
was supported by the NSERC (National Sciences and Engineering Research Council of Canada) 
through the grant PDF-329052-2006. 


\newcommand{\auth}{\textsc}

 
\small 

\section*{Appendix}
\label{A1}
For completeness we list here the high-order discretizations used in this paper which, for instance, 
we state 
for a conservation law with conservative variable $u$ and flux $f=f(u)$. 
\begin{itemize}
\item Fourth order discretization  ($q=4$):
$$
h \, f_x = {1\over 12}f_{i-2}-{2\over 3}f_{i-1}+{2\over 3}f_{i+1}-{1\over 12}f_{i+2},
$$
$$
h^2 \, u_{xx} = -{1\over 12}u_{i-2}+{4\over 3}u_{i-1}-{5\over 2}u_{i}+{4\over 3}u_{i+1}-{1\over 12}u_{i+2},
$$
$$
h^3 \, u_{xxx} = -{1\over 2}u_{i-2}+ u_{i-1}- u_{i+1}+{1\over 2}u_{i+2}.
$$
\item Sixth order discretization ($q=6$):
$$
h \, f_x =  -{1\over 60}f_{i-3}+{3\over 20}f_{i-2}-{3\over 4}f_{i-1}
+{3\over 4}f_{i+1}-{3\over 20}f_{i+2}+{1\over 60}f_{i+3},
$$
$$
{h^2 \over 2} \, u_{xx} = {1\over 180}u_{i-3}-{3\over 40}u_{i-2}+{3\over 4}u_{i-1}
-{49\over 36}u_{i}+{3\over 4}u_{i+1}-{3\over 40}u_{i+2}+{1\over 180}u_{i+3},
$$
$$
{h^3 \over 6} \, u_{xxx} =
{1\over 48}u_{i-3}-{1\over 6}u_{i-2}+{13\over 48}u_{i-1}-{13\over 48}u_{i+1}+{1\over 6}u_{i+2}-{1\over 48}u_{i+3}.
$$
\item Eighth order discretization  ($q=8$):
 $$
 \aligned
 \, f_x =  & {1\over 280}f_{i-4}-{4\over 105}f_{i-3}+{1\over 5}f_{i-2}
      -{4\over 5}f_{i-1}+{4\over 5}f_{i+1}-{1\over 5}f_{i+2}
 \\
 & +{4\over 105}f_{i+3}-{1\over 280}f_{i+4},
 \endaligned
 $$
 $$
 \aligned {h^2 \over 2} \, u_{xx}= & {-1\over 1120}u_{i-4}+{4\over
 315}u_{i-3}-{1\over 10}u_{i-2}+{4\over 5}u_{i-1}-{205\over
 144}u_{i}+{4\over 5}u_{i+1}
 \\
 &  -{1\over 10}u_{i+2}+{4\over 315}u_{i+3}-{1\over 1120}u_{i+4},
 \endaligned
 $$
 $$
 \aligned {h^3 \over 6} \, u_{xxx}= & {-7\over 1440}u_{i-4}+{1\over
 20}u_{i-3}-{169\over 720}u_{i-2}+{61\over 180}u_{i-1}-{61\over
 180}u_{i+1}
 \\
 & +{169\over 720}u_{i+2}-{1\over 20}u_{i+3}+{7\over 1440}u_{i+4}.
 \endaligned
 $$
\item Tenth order discretization ($q=10$):
 $$
 \aligned h \, f_x = & {-1\over 1260}f_{i-5}+{5\over
 504}f_{i-4}-{5\over 84}f_{i-3}+{5\over 21}f_{i-2}-{5\over
 6}f_{i-1}+{5\over 6}f_{i+1}-{5\over 21}f_{i+2}
 \\
 & +{5\over 84}f_{i+3}-{5\over 504}f_{i+4}+{1\over 1260}f_{i+5},
 \endaligned
 $$
 $$
 \aligned {h^2 \over 2} \, u_{xx}= & {1\over 6300}u_{i-5}-{5\over
 2016}u_{i-4}+{5\over 252}u_{i-3}-{5\over 42}u_{i-2}+{5\over
 6}u_{i-1}-{5269\over 3600}u_{i}+{5\over 6}u_{i+1}
 \\
 & -{5\over 42}u_{i+2} +{5\over 252}u_{i+3}-{5\over
 2016}u_{i+4}+{1\over 6300}u_{i+5},
 \endaligned
 $$
 $$
 \aligned {h^3 \over 6} \, u_{xxx}= & {41\over
 36288}u_{i-5}-{1261\over 90720}u_{i-4}+{541\over
 6720}u_{i-3}-{4369\over 15120}u_{i-2} +{1669\over 4320}u_{i-1}
 \\
 & -{1669\over 4320}u_{i+1} +{4369\over 15120}u_{i+2}  -{541\over
 6720}u_{i+3}+{1261\over 90720}u_{i+4}-{41\over 36288}u_{i+5}.
 \endaligned
 $$
\end{itemize}

\end{document}